\documentclass[11pt]{article}
\usepackage{amsmath}
\usepackage{amsthm}
\usepackage{amsfonts}
\usepackage{amssymb}
\usepackage{color}
\usepackage{graphicx}
\usepackage{url}
\usepackage{verbatim}

\usepackage[top=1.25cm, bottom=2cm, left=2.0cm, right=2cm]{geometry} 
\newcommand{\R}{{\mathbb R}}
\newcommand{\cM}{{\mathcal M}}

\usepackage[noend]{algorithmic}
\usepackage{algorithm}
\floatname{algorithm}{Algorithm}

\usepackage{cmbright}
\usepackage[T1]{fontenc}

\title{Topological Microstructure Analysis Using Persistence Landscapes}

\author{Pawe{\l} D{\l}otko \\
        Inria Saclay -- Ile-de-France \\
        1 rue Honor\'e d'Estienne d'Orves \\
        91120 Palaiseau, France
  \and
        Thomas Wanner \\
        Department of Mathematical Sciences \\
        George Mason University \\
        Fairfax, VA 22030, USA }

\date{ \today}

\begin{document}       

\maketitle

\begin{abstract}

Phase separation mechanisms can produce a variety of complicated and
intricate microstructures, which often can be difficult to characterize
in a quantitative way. In recent years, a number of novel topological
metrics for microstructures have been proposed, which measure essential
connectivity information and are based on techniques from algebraic
topology. Such metrics are inherently computable using computational
homology, provided the microstructures are discretized using a 
thresholding process. However, while in many cases the thresholding is
straightforward, noise and measurement errors can lead to misleading
metric values. In such situations, persistence landscapes have been
proposed as a natural topology metric. Common to all of these approaches
is the enormous data reduction, which passes from complicated patterns to
discrete information. It is therefore natural to wonder what type of
information is actually retained by the topology. In the present paper,
we demonstrate that averaged persistence landscapes can be used to
recover central system information in the Cahn-Hilliard theory of phase
separation. More precisely, we show that topological information of
evolving microstructures alone suffices to accurately detect both
concentration information and the actual decomposition stage of a
data snapshot. Considering that persistent homology only measures
discrete connectivity information, regardless of the size of the
topological features, these results indicate that the system parameters
in a phase separation process affect the topology considerably more
than anticipated. We believe that the methods discussed in this paper
could provide a valuable tool for relating experimental data to
model simulations.

\end{abstract}
\section{Introduction}
Complicated patterns which evolve with time occur in a variety of
applied contexts, and quantifying or even just distinguishing such
patterns can pose serious challenges. Over the last decade,
computational topology has emerged as a tool which on the one hand
allows for significant data reduction, while at the same time 
focusing on and retaining essential connectivity information of
the studied patterns. One particularly interesting application
area is materials science, where complex evolving patterns are
frequently created through complicated phase separation processes.
Computational topology encompasses a wide variety of possible
tools~\cite{herbert}, and most of them in one way or another are based
on homology theory. The deeper reason for this can be found in the
inherent computability of homology, and in recent years powerful
algorithms have been devised which enable fast homology computations
for very large data sets, see for example~\cite{dlotko:etal:11a,
kaczynski:etal:04a}, as well as the references therein.

Homology has been used in a number of materials science contexts.
Earlier studies have made use of the easily computable Euler
characteristic, see for example the references in the recent
survey~\cite{wanner:p15a}. However, it has been pointed out that
in some situations the information retained by the Euler
characteristic is not enough to distinguish certain important
pattern features. In contrast, the Betti numbers, which are 
associated with the homology groups and will be described below,
provide a finer metric. They were used for example
in~\cite{gameiro:etal:05a} to relate the pattern complexity
evolution as described by averaged Betti number evolution curves
to the amount of stochasticity inherent in a phase field model
for binary phase separation in metal alloys. This study was the
first to use homology information in the context of model 
validation. Based on available experimental data it was shown
that if the noise in the system is too low, the observed Betti
number evolution curves are qualitatively different from the
experimental ones. In addition, it was demonstrated
in~\cite{gameiro:etal:05a} that while Betti numbers can be
used to separate bulk from boundary behavior, the averaged
Euler characteristic can only describe the boundary effects.
Similar materials science studies in the context of polycrystals
can be found in~\cite{rohrer:miller:10a, fuller:etal:10a}. While
the first of these papers uses homology to study the connectivity
properties of grain-boundary networks in planar sections of
polycrystals, the second paper employs Betti numbers as a means
to describe the thermal-elastic response of calcite-based
polycrystalline materials such as marble. More precisely,
in~\cite{fuller:etal:10a} homological techniques are used 
to characterize not only the elastic energy density and maximum
principal stress response fields in a polycrystal, but also the
respective grain-boundary misorientation distributions which 
generated these response fields. It was shown that this topological
analysis can quantitatively distinguish between different types of
grain-boundary misorientations, and relate them to differences in
the resulting response fields.

In all of the applications described so far, the numerical or
experimental data is given in the form of a field, or in other
words, a real-valued function~$u : \Omega \to \R$ defined on
some domain~$\Omega \subset \R^d$. The associated patterns 
are subsets of the domain~$\Omega$, and usually created through
a thresholding process. For example, after selecting a suitable
threshold~$\theta$, one can consider the sub- and super-level sets
\begin{displaymath}
  \cM^\pm \; = \;
  \left\{ x \in \Omega \; : \;
  \pm\left( u(x) - \theta \right) \ge 0 \right\} \; ,
\end{displaymath}
which in some contexts are called nodal domains of~$u$. As
subsets of~$\Omega \subset \R^d$, the sets~$\cM^\pm$ have 
well-defined singular homology groups, and if the function~$u$
is sufficiently regular, these groups can be computed
precisely; see for example~\cite{cochran:etal:13a, day:etal:09a} 
for the two-dimensional case $d = 2$. However, if the function~$u$
is not smooth, or if the thresholding process involves a field
created from experimental or noisy data, then the above
thresholding process may not capture the correct topology
of the actual underlying pattern. Moreover, in certain
applications the thresholding approach itself might not 
be appropriate, for example if there is no obvious or
physically relevant choice of threshold.

An extension of the concept of homology to situations
involving noise or the lack of a clear thresholding process
has been proposed some fifteen years ago and is called
persistent homology~\cite{herbert}. While this extension
is described in more detail in the next section, persistent
homology is a dimension reduction technique which provides
a topological description of the evolution of sublevel sets
of a mapping~$u$ as a function of the thresholding
level~$\theta$. It gives rise to intervals~$[\theta_1,
\theta_2]$ over which certain topological features persist,
and the length~$\theta_2 - \theta_1$ can in some sense be
viewed as a measure of importance of the specific feature.
Thus, topological features with small interval lengths are
usually considered as ``noise'' or ``unimportant,'' while 
features with long intervals are deemed ``significant.''
Needless to say, the precise meaning of these notions will change
from application to application. The concept of persistent homology
has already been used in a number of materials science contexts,
such as for example in the analysis of granular
media~\cite{granularMedia2, granularMedia1}, in the study of
protein compressibility~\cite{proteinComplexity}, as well as in
the classification of amorphous structures~\cite{amorphousStructure}
and glass~\cite{glass}. Efficient algorithms for computing 
persistent homology are described in~\cite{herbert, mrozek:wanner:10a}.

Despite the success of the above uses of computational topology
in applications, one immediate question is the extent of the resulting
dimension reduction. Through homology, patterns or microstructures are
basically reduced to a finite set of integers, and it is therefore
natural to wonder what information is still encoded in this reduced
measurement, beyond the obvious connectivity information. It was pointed
out in~\cite{desi:etal:11a} that even small Betti number counts,
which taken in isolation do not provide much in terms of pattern
differentiation, can provide significant information when viewed
in a stochastic, i.e., averaged setting. More precisely, it was shown
in~\cite{desi:etal:11a} that during the phase separation process 
called nucleation, averaged droplet counts on small domains give
extremely precise projections for the observed averaged droplet
counts on large domains, where the droplet count is just the
zero-dimensional Betti number.
\begin{figure} \centering
  \setlength{\unitlength}{1 cm}
  \begin{picture}(14.9,7.3)
    \put(0.0,3.8){
      \includegraphics[width=3.5cm]{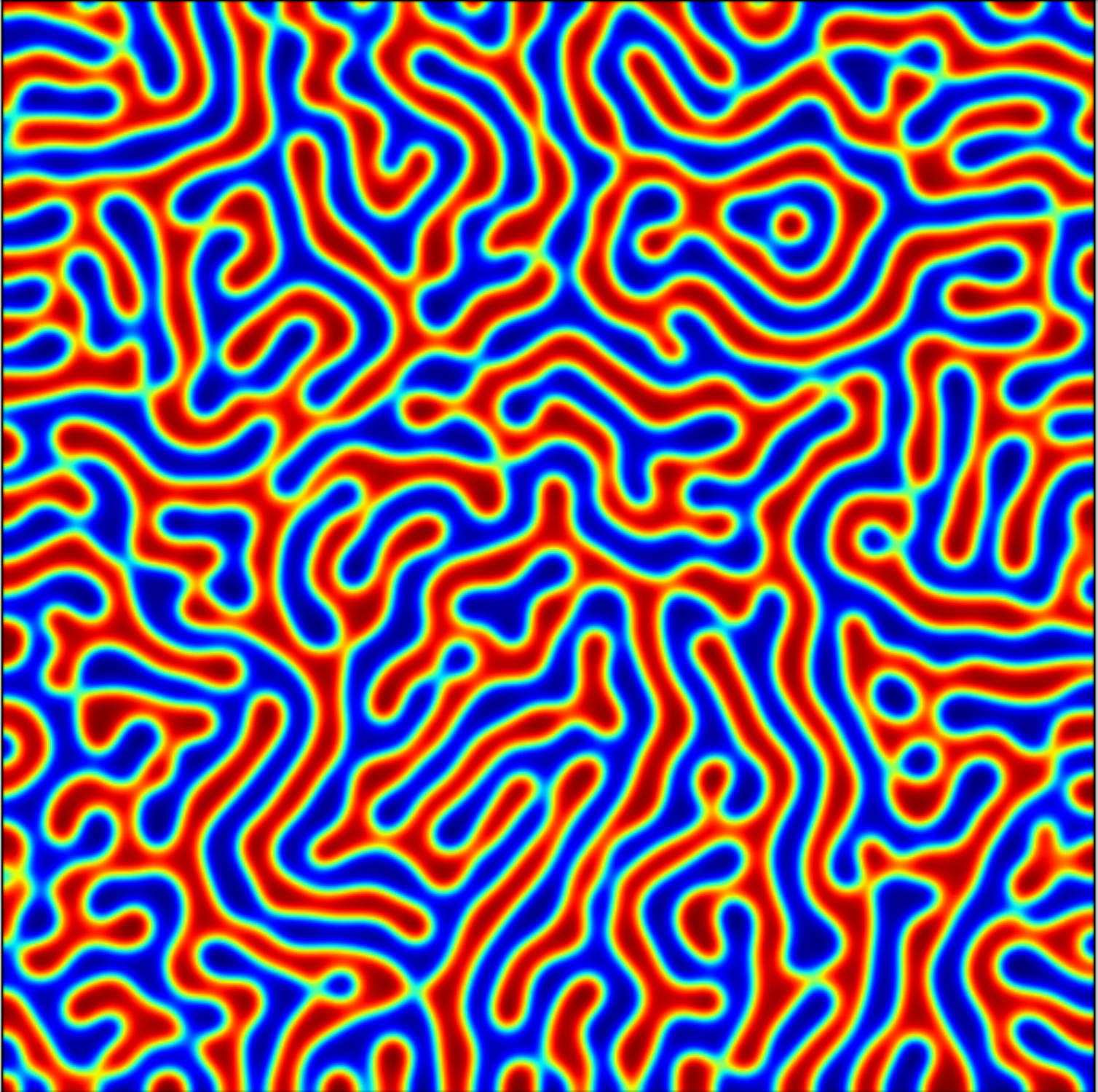}}
    \put(3.8,3.8){
      \includegraphics[width=3.5cm]{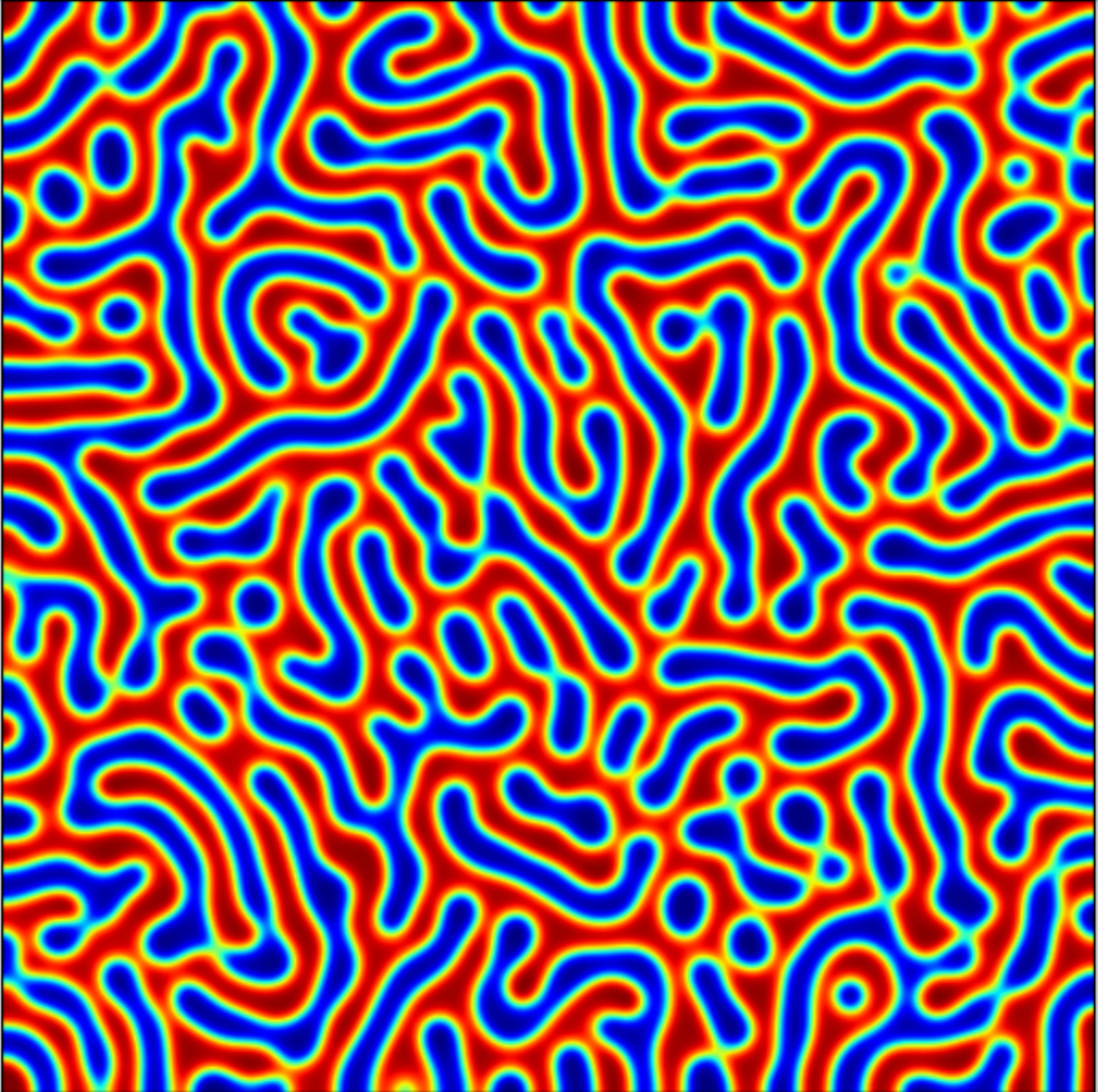}}
    \put(7.6,3.8){
      \includegraphics[width=3.5cm]{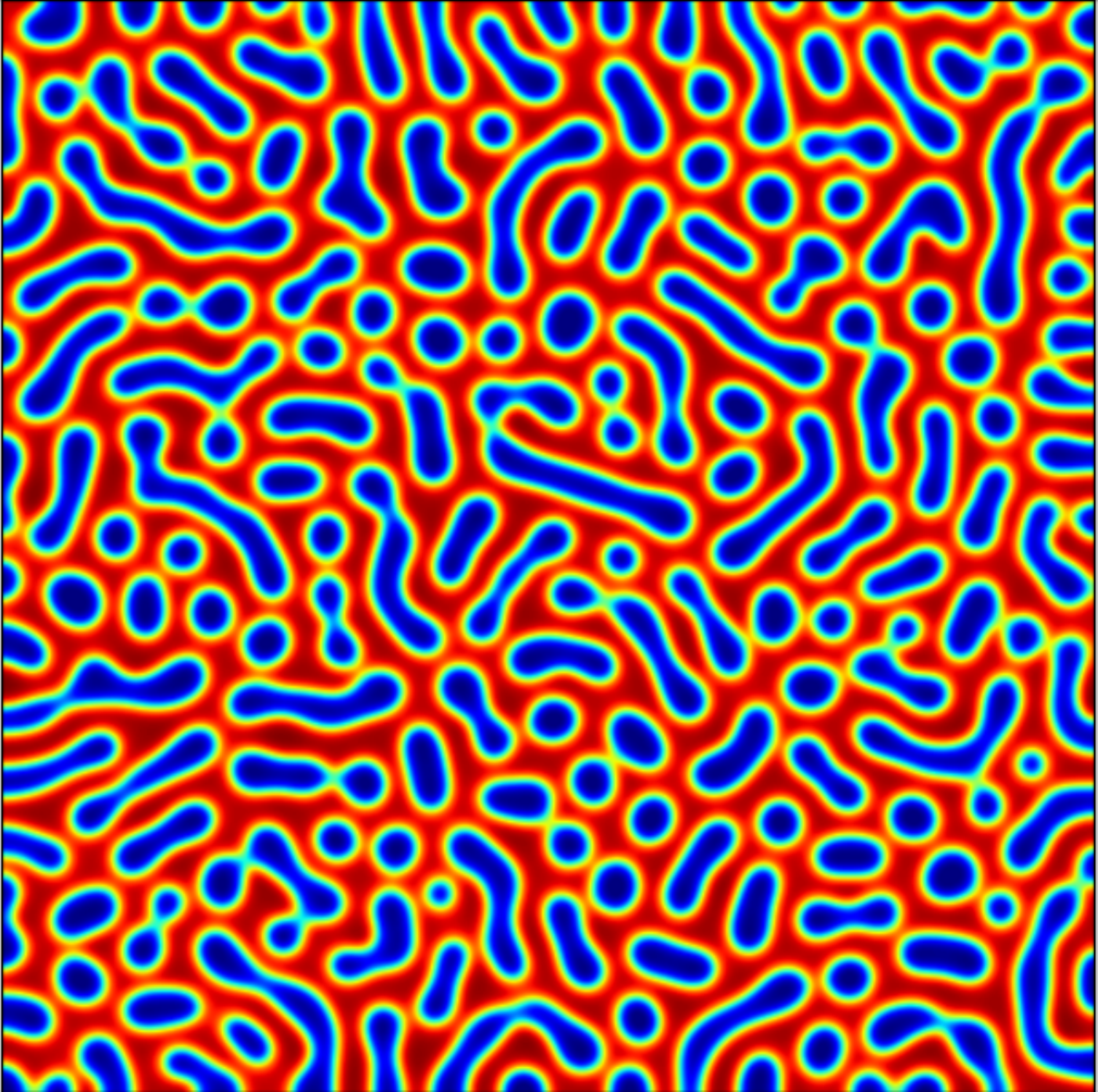}}
    \put(11.4,3.8){
      \includegraphics[width=3.5cm]{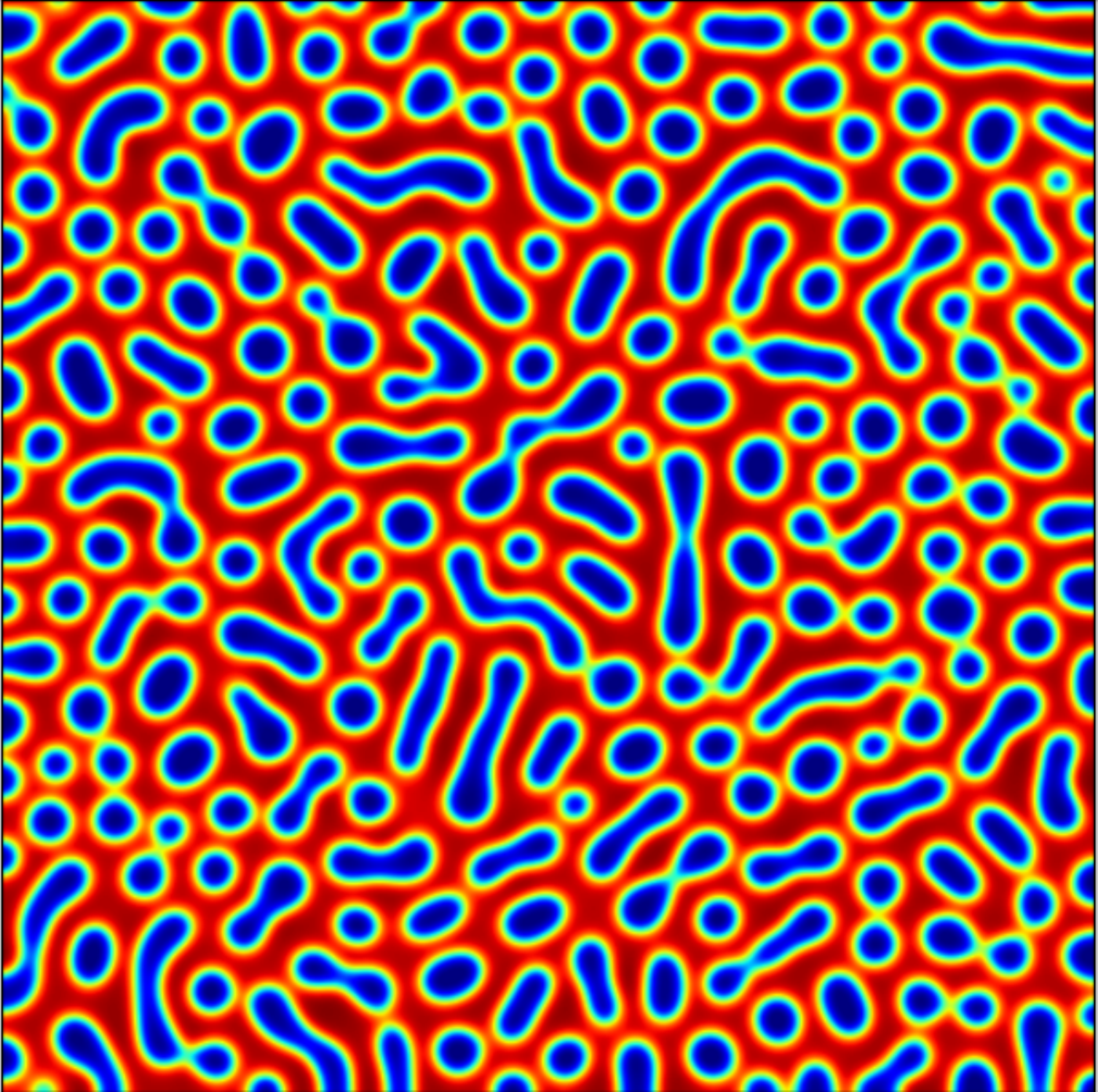}}
    \put(0.0,0.0){
      \includegraphics[width=3.5cm]{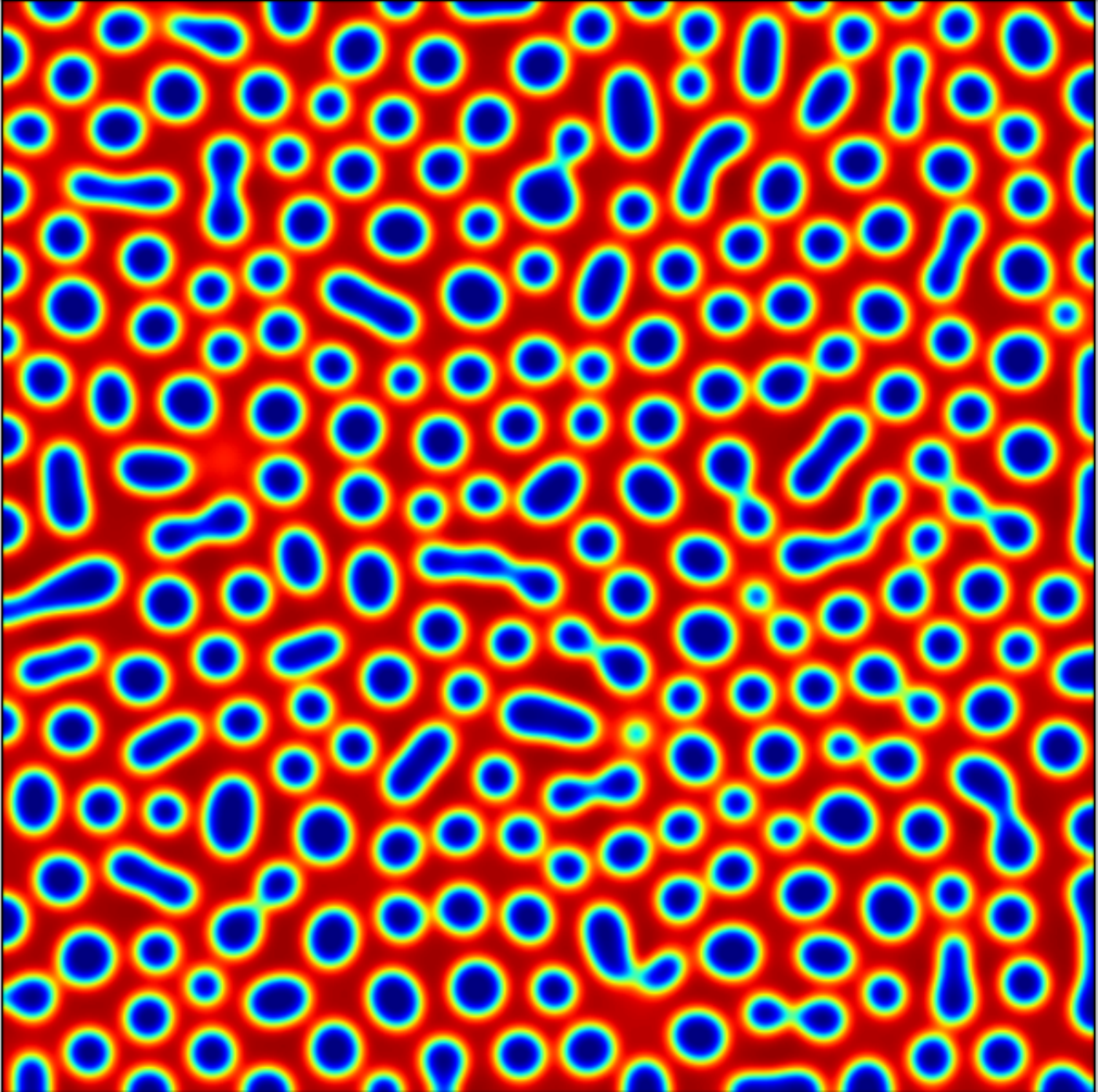}}
    \put(3.8,0.0){
      \includegraphics[width=3.5cm]{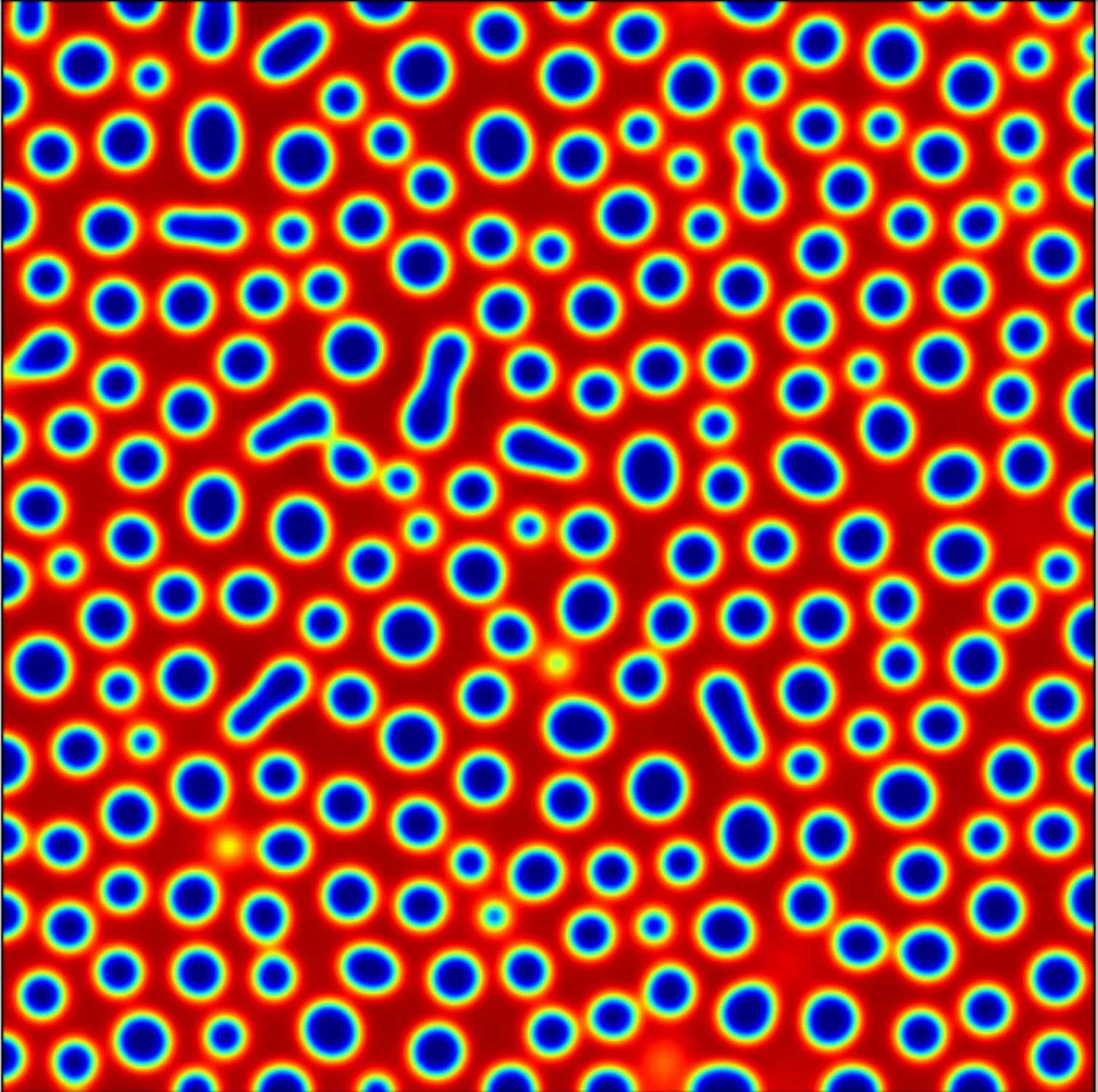}}
    \put(7.6,0.0){
      \includegraphics[width=3.5cm]{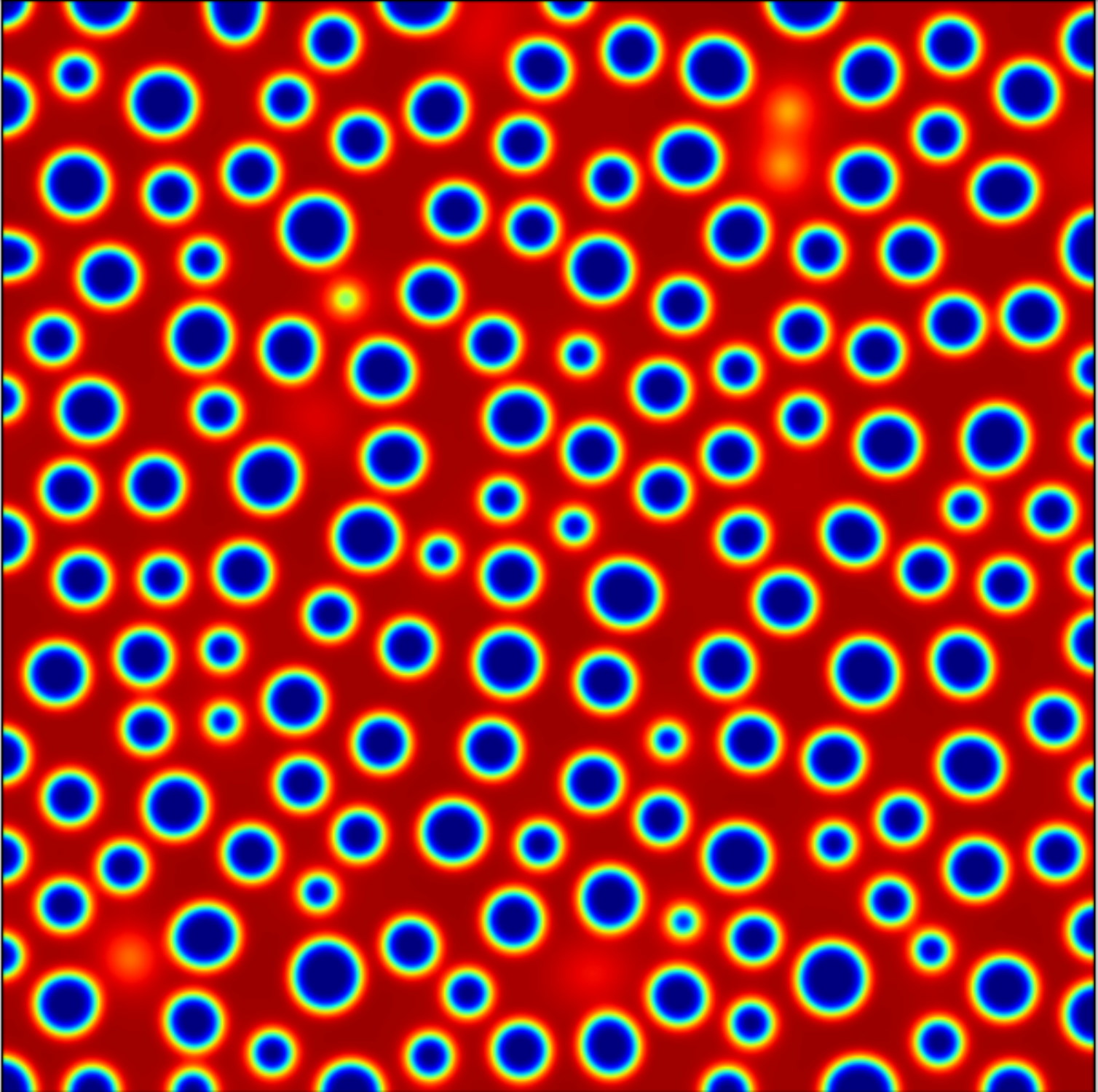}}
    \put(11.4,0.0){
      \includegraphics[width=3.5cm]{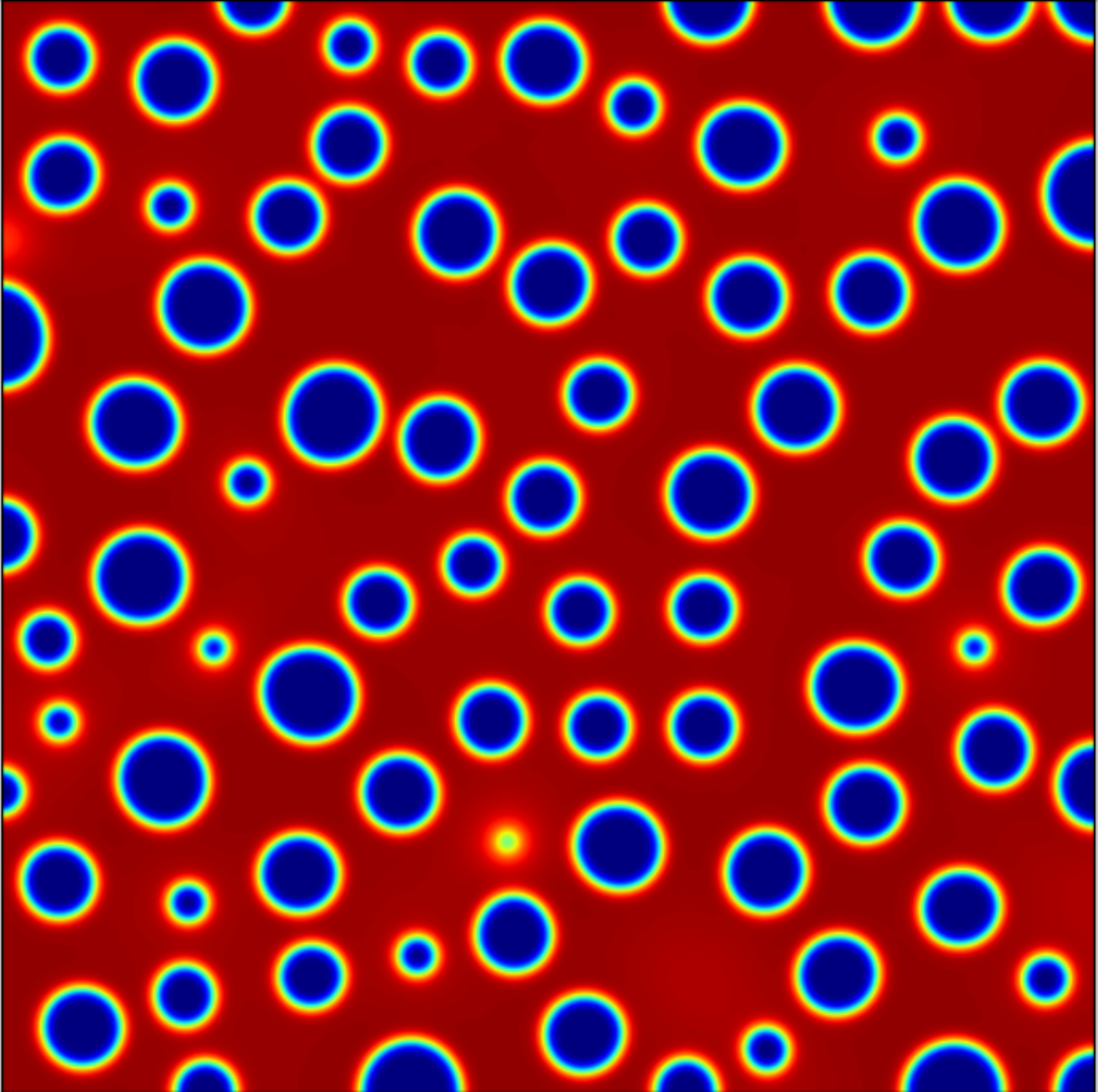}}
  \end{picture} 
  \caption{Phase separation through spinodal decomposition in the
           Cahn-Hilliard-Cook model~(\ref{chc}) with white noise forcing.
           The images show snapshots of solutions which originate at the
           homogeneous state~$u \equiv \mu$, and for parameter values
           $\epsilon=0.005$ and $\sigma = 0.001$. From top left to bottom
           right the images are for mass~$\mu = 0.07 \cdot k$,
           where $k = 0,\ldots,7$.}
  \label{fig:chcpatt}
\end{figure}

In the present paper, we demonstrate that when viewed in a stochastic
and time evolving framework, topological information encodes considerably
more than anticipated. This will be accomplished in the setting of
instantaneous phase separation in binary metal alloys, as modelled
by the Cahn-Hilliard theory of spinodal decomposition~\cite{bloemker:etal:05a,
cahn:68a}. This phase separation phenomenon is initiated immediately
after a high-temperature melt of two uniformly mixed metal components
is quenched, i.e., rapidly cooled. Depending on the concentrations of
the involved components, they will quickly separate to form complicated
microstructures which contain some apparent element of randomness. Some
of the resulting patterns in two space dimensions are shown in
Figure~\ref{fig:chcpatt}. These patterns evolve with time, and
after the initial phase separation process a coarsening stage sets
in, during which the characteristic length scale of the patterns
increases.

The first mathematical model for spinodal decomposition 
was introduced by Cahn and Hilliard~\cite{cahn:59a, cahn:hilliard:58a},
who proposed a nonlinear evolution equation for the relative concentration
difference $u = \rho_A - \rho_B$, where~$\rho_A$ and~$\rho_B$ denote
the relative concentrations of the two components, i.e., $\rho_A +
\rho_B = 1$. Their model is based on the Ginzburg-Landau free energy
given by
\begin{equation} \label{defenergy}
  E_\epsilon(u) = \int_\Omega \left( \frac{\epsilon^2}{2} |\nabla u|^2
                  + F(u) \right) \, dx \; ,
\end{equation}
where~$\Omega \subset \R^d$ is a bounded domain, and the positive
parameter~$\epsilon$ models interaction distance. The bulk free energy~$F$
is a double well potential, which for the purposes of this paper is
taken as
\begin{equation} \label{defF}
  F(u) = \frac{1}{4} \left( u^2 - 1 \right)^2 \; .
\end{equation}
Taking the variational derivative $\delta E_\epsilon / \delta u$ of
the Ginzburg-Landau free energy~(\ref{defenergy}) with respect to the
concentration variable~$u$, one then obtains first the chemical
potential $w = -\epsilon^2 \Delta u + F'(u)$, and then the
associated Cahn-Hilliard equation $\partial u / \partial t =
\Delta w$, i.e., the fourth-order partial differential equation
\begin{equation} \label{ch}
  \frac{\partial u}{\partial t} =
  -\Delta \left( \epsilon^2 \Delta u - F'(u) \right) \; ,
\end{equation}
subject to homogeneous Neumann boundary conditions for both~$w$
and~$u$. Due to these boundary conditions, any mass flux through the
boundary is prohibited, and therefore mass is conserved. We generally
consider initial conditions for~(\ref{ch}) which are small-amplitude
random perturbations of a spatially homogeneous state, i.e., we
assume that $u(0,x) \approx \mu$ for all $x \in \Omega$, as well as
$\int_\Omega u(0,x) \, dx / |\Omega| = \mu$, where we use the standard
abbreviation $|\Omega| = \int_\Omega 1 \, dx$. One can easily see
that such initial conditions lead to instantaneous phase separation
in the Cahn-Hilliard equation as long as the conserved total mass~$\mu$
satisfies $F''(\mu) < 0$. In the context of~(\ref{defF}), this
condition is equivalent to~$|\mu| < 1 / \sqrt{3}$.

While the deterministic Cahn-Hilliard model~(\ref{ch}) provides
a basic qualitative model for spinodal decomposition, it cannot
produce microstructures whose evolution during the initial phase
separation stage agrees with experiments~\cite{gameiro:etal:05a}.
As it turns out, this deficiency is due to the fact that the
original model~(\ref{ch}) completely ignores thermal fluctuations.
To remedy this, Cook~\cite{cook:70a} extended the model by adding
a random fluctuation term~$\xi$, i.e., he considered the
stochastic Cahn-Hilliard-Cook model
\begin{equation} \label{chc}
  \frac{\partial u}{\partial t} =
  -\Delta \left( \epsilon^2 \Delta u - F'(u) \right) +
    \sigma \cdot \xi \; ,
\end{equation}
where the expected value and the correlation of the noise
process satisfy both
\begin{displaymath}
  \langle \xi(t,x) \rangle = 0
  \qquad\mbox{ and }\qquad
  \langle \xi(t_1,x_1) \xi(t_2,x_2) \rangle = 
    \delta(t_1 - t_2) q(x_1 - x_2) \; .
\end{displaymath}
Here $\sigma > 0$ is a measure for the intensity of the fluctuation
and~$q$ describes the spatial correlation of the noise. Thus, the
noise is always uncorrelated in time, and for the special case~$q =
\delta$ we obtain space-time white noise. In general, however, the
noise process will exhibit spatial correlations. For more details we
refer the reader to the discussions in~\cite{bloemker:etal:05a, langer:71a}.
As mentioned before, while both the deterministic and the stochastic
model generate microstructures which are qualitatively similar to the
ones observed during spinodal decomposition~\cite{cahn:65a}, only the 
stochastic model can in principle lead to a more quantitative
agreement~\cite{gameiro:etal:05a}. Recent mathematical results
for the models~(\ref{ch}) and~(\ref{chc}) have identified the
observed microstructures as certain random superpositions of
eigenfunctions of the Laplacian, and were able to explain the
dynamics of the decomposition process in more
detail~\cite{bloemker:etal:01b, bloemker:etal:08a,
sander:wanner:99a, sander:wanner:00a, wanner:04a}.

Based on the insight obtained in~\cite{gameiro:etal:05a}, our
studies in the present paper concentrate exclusively on the
stochastic Cahn-Hilliard-Cook model~(\ref{chc}). For this model,
we show that persistent homology averages suffice to accurately
deduce the total mass~$\mu$, and therefore the total alloy component
concentrations, as well as the actual stage of the decomposition
process. Both of these observations are somewhat surprising, since
there is no obvious link between the connectivity information retained
by homology and either of these quantities. Our results demonstrate,
however, that during the phase separation process, the evolution of
averaged persistent homology information in some sense is a very
detailed encoding of these system parameters. Expressed differently,
while each mass value~$\mu$ uniquely determines an averaged persistent
homology evolution for randomly selected solutions originating
close to the homogeneous state~$\mu$, the topology encoded in these
averaged evolutions is sufficiently different to allow for the solution
of an inverse problem. A similar statement can be made for the stage
of the decomposition. Moreover, we wold like to point out that
unlike the Betti numbers, persistent homology information cannot
easily be averaged directly. We therefore make use of the recent concept
of persistence landscapes~\cite{peterLandscapes, plt}, which will
be described in more detail below. In our studies, we are also using the
idea of topological process, which will be introduced later in this paper.
It allows us to deal with a processes that produces a sequence of topological
spaces.

The remainder of this paper is organized as follows. In
Section~\ref{sec:topology} we survey the topological background
for this paper, to keep our presentation self-contained and
accessible for applied scientists. After introducing cubical
homology in Section~\ref{sec:introcubhom}, we discuss the concept
of persistent homology in Section~\ref{sec:intropers}. Since for
our study the averaging of persistence information is crucial,
we then turn our attention to persistence landscapes in
Section~\ref{sec:intropl}. After these preparations, the topological
analysis of microstructures created by the Cahn-Hilliard-Cook
model~(\ref{chc}) is the subject of Section~\ref{sec:chc}. We
begin by presenting the basic methodology in Section~\ref{sec:basicmethod},
turn our attention to the topological determination of the total
mass~$\mu$ in Section~\ref{sec:chcmass}, and show in Section~\ref{sec:chctime}
that even the time of a solution snapshot can be recovered accurately
from the persistence information during the spinodal decomposition and
the initial coarsening regime. Finally, Section~\ref{sec:future} contains
a brief summary and draws some conclusions.
\section{Topological Background}
\label{sec:topology}
In order to keep this paper as self-contained as possible, this section
provides background information on homology and persistence. More precisely,
we recall the basic definition of cubical homology, introduce persistent
homology, and finally describe the concept of persistence landscapes.
\subsection{Cubical Homology}
\label{sec:introcubhom}
Homology theory associates discrete objects with topological spaces in
such a way that continuous deformations of the space do not change
the objects. In other words, homology is a topological invariant of the
underlying topological space. We will see below that the discrete objects
are Abelian groups, and depending on the form in which the topological
space is given, a number of different homology theories have been developed
over the years to compute them. If a given space can be treated within
several different homology theories, then the homology groups assigned
by these different homology theories will be isomorphic. For the purposes
of this paper, we restrict our attention to one specific homology theory,
namely \emph{cubical homology}. For this, we assume that our topological
space is given as a specific union of cubical sets of possibly
varying dimensions, leading to the concept of a \emph{cubical complex}.

\begin{figure} \centering
  \setlength{\unitlength}{1 cm}
  \begin{picture}(14.5,5.0)
    \put(0.0,4.5){
      \parbox{1cm}{\emph{(a)}}}
    \put(1.5,0.0){
      \includegraphics[height=5cm]{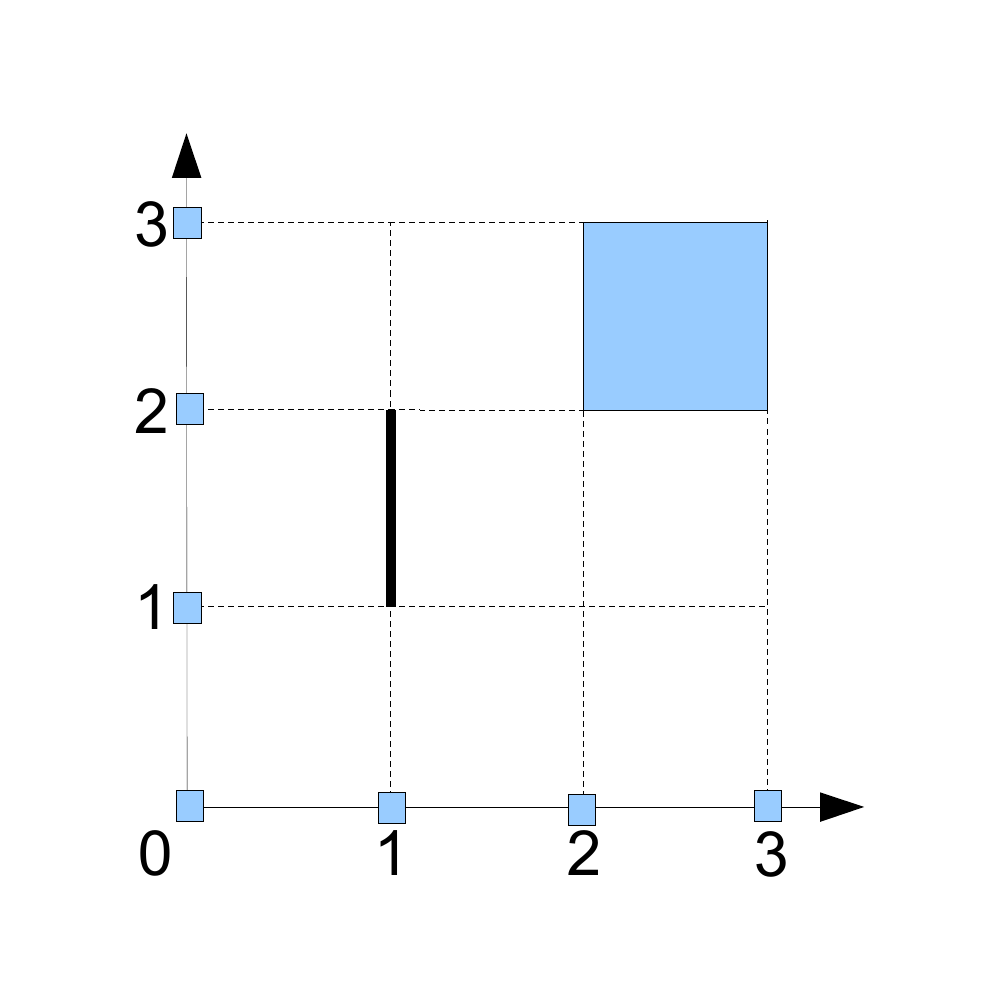}}
    \put(8.0,4.5){
      \parbox{1cm}{\emph{(b)}}}
    \put(9.5,0.0){
      \includegraphics[height=5cm]{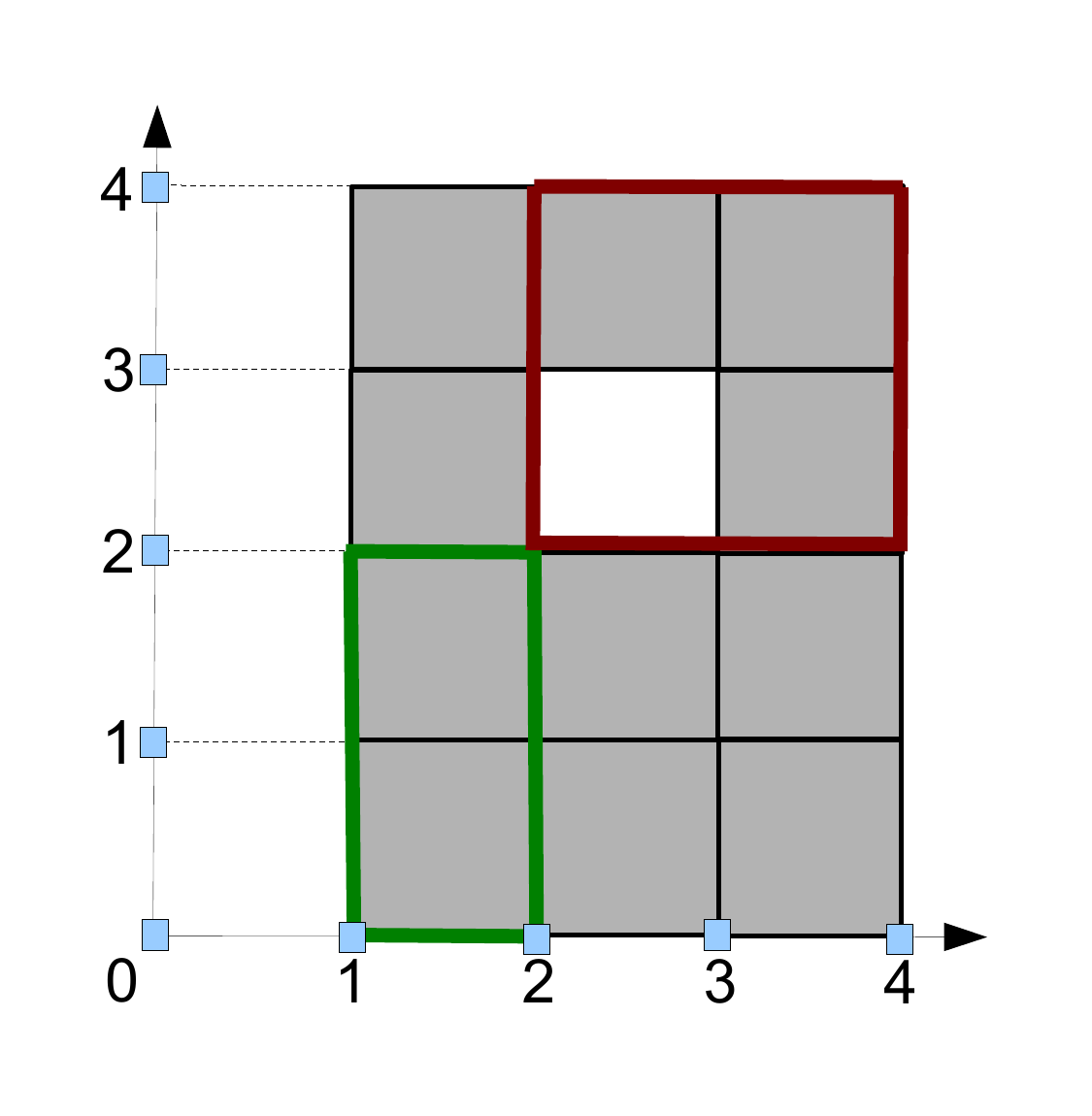}}
  \end{picture}
  \caption{Illustration of two fundamental concepts from cubical homology.
           \emph{(a) Boundary Operator}: The left image shows the two cubes~$C_1
           = [1,1] \times [1,2]$ and~$C_2 = [2,3] \times [2,3]$, whose dimensions
           are one and two, respectively. Their boundaries can easily be
           computed as $\partial C_1 = [1,1] \times [1,1] + [1,1] \times [2,2]$,
           as well as $\partial C_2 = [2,2] \times [2,3] + [3,3] \times [2,3]
           + [2,3] \times [2,2] + [2,3] \times [3,3]$. Note that in both cases,
           the formal notion of boundary is consistent with the intuitive
           boundary of that cube.
           \emph{(b) Cycles}: Consider the cubical complex shown in gray in 
           the right image, together with all contained edges and vertices. Then
           both the collection of edges shown in red, and the collection shown
           in green form one-dimensional chains which are cycles, i.e., their
           boundaries vanish. This is due to the fact that every cell of
           dimension zero in the support of each chain occurs in the boundary
           of exactly two cells, and therefore cancels out in the chain's
           boundary.}
  \label{fig:cubhom}
\end{figure}

To be more precise, define an \emph{elementary interval} as a compact
real interval of the type~$[\ell,\ell+1]$ or~$[\ell,\ell]$, where~$\ell$
denotes an integer. While the first interval type is referred to as
\emph{non-degenerate}, the second one is called a \emph{degenerate
interval}. An \emph{elementary cube}~$C$ is a product of elementary
intervals $I_1 \times \ldots \times I_n$, and its \emph{dimension}
equals the number of non-degenerate intervals~$I_k$ in this product.
Furthermore, the total number~$n$ of intervals in the product is
called the \emph{embedding dimension} of the cube~$C$. For example,
Figure~\ref{fig:cubhom}\emph{(a)} depicts an elementary cube of
dimension one in black, and a two-dimensional elementary cube in
blue, and both cubes have embedding dimension two.

While elementary cubes form the building blocks of the cubical complexes
introduced below, we first need to introduce a few algebraic concepts
based on them. A \emph{chain}~$c$ is a formal sum
\begin{equation} \label{def:chain}
  c = \alpha_1 C_1 + \ldots + \alpha_m C_m \; ,
  \quad\mbox{where }
  C_1, \ldots, C_m
  \mbox{ are elementary cubes and }
  \alpha_1, \ldots, \alpha_m \in F
  \mbox{ are scalars.}
\end{equation}
In this formula, $F$ denotes an arbitrary field. However, for the purposes
of this paper, we always consider the field~$F = \mathbb{Z}_2$ which
consists of the two elements~$0$ and~$1$, as well as addition and
multiplication modulo two. Thus, for the chain~$c$ in~(\ref{def:chain})
a cube~$C_k$ is either present or absent, depending on whether
$\alpha_k = 1$ or $\alpha_k = 0$, respectively. Furthermore, the
\emph{support} of the chain~$c$ is the collection of all elementary
cubes~$C_k$ which have nonzero coefficient~$\alpha_k$ in~$c$. Two
one-dimensional chains consisting of four one-dimensional elementary
cubes each are shown in green and red in Figure~\ref{fig:cubhom}\emph{(b)}.

Next we have to introduce the \emph{algebraic boundary} of a chain.
For an elementary interval $I = [n,m]$, its boundary is the chain
supported by the two degenerate intervals~$[m,m]$ and~$[n,n]$,
and this immediately implies both
\begin{displaymath}
  \partial [n,n+1] = [n+1,n+1] + [n,n]
  \qquad\mbox{ as well as }\qquad
  \partial [n,n] = 0 \; ,
\end{displaymath}
depending on whether $m = n+1$ or $m = n$, respectively. Note
that if the interval~$I$ is degenerate, we have used the fact that
$\partial [n,n] = [n,n] + [n,n] = (1+1) [n,n] = 0$ in the
field~$F = \mathbb{Z}_2$. Combining the definition of the boundary
for elementary intervals with the definition of elementary cubes, 
one can then lift the boundary definition to the latter using 
linearity. More precisely, let $C = I_1 \times \ldots \times I_n$
denote an elementary cube. Then we define
\begin{displaymath}
  \partial C = \partial \left( I_1 \times \ldots \times I_n \right) =
  \sum_{i=1}^n I_1 \times \ldots \times \partial I_i \times
    \ldots \times I_n \; ,
\end{displaymath}
where $I_1 \times \ldots \times \partial [m,m+1] \times \ldots \times I_n
= I_1 \times \ldots \times [m+1,m+1] \times \ldots \times I_n + I_1 \times
\ldots \times [m,m] \times \ldots \times I_n$ in the nondegenerate case,
and $I_1 \times \ldots \times \partial [m,m] \times \ldots \times I_n = 0$
in the degenerate one. Finally, the boundary operator extends linearly
from elementary cubes to chains, i.e., for a chain $c = \alpha_1 C_1 +
\ldots + \alpha_\ell C_\ell$ one defines its boundary via $\partial c =
\alpha_1 \partial C_1 + \ldots + \alpha_\ell \partial C_\ell$. These
concepts are illustrated further in Figure~\ref{fig:cubhom}. 

The boundary operator formalizes the intuitive notion of the boundary
of a $k$-dimensional elementary cube. For $k = 0$, such a cube is a point,
and its boundary is zero, while for $k = 1$ the cube is an interval with
the boundary consisting of its endpoints. Note, however, that this notion
does not depend on the embedding dimension of the cube, in contrast to the
topological boundary. In other words, the boundary~$\partial C$ is
determined from the intrinsic dimension of the cube. This latter fact
also lies at the heart of the fundamental property of the boundary 
operator, which states that
\begin{displaymath}
  \partial\partial c = 0
  \qquad\mbox{ for every chain } c \; .
\end{displaymath}
One can easily see why this statement is true for elementary cubes.
If~$C$ is an elementary cube of dimension~$k$, then its boundary is
a union of elementary cubes of dimension~$k-1$. This implies that the
chain~$\partial\partial C$ consists of cubes of dimension~$k-2$.
Every one of these $k-2$-dimensional cubes is contained in the boundary
of exactly two $k-1$-dimensional cubes in~$\partial C$, i.e., the
chain~$\partial\partial C$ contains every elementary cube exactly
twice. Since we are considering the field~$F = \mathbb{Z}_2$, these
cubes cancel out and yield a zero boundary. This is illustrated in
Figure~\ref{fig:cubhom}\emph{(b)}.

As we mentioned earlier, elementary cubes are the fundamental building
blocks for the topological spaces which can be treated using cubical
homology. Loosely speaking, we will consider subsets of~$\mathbb{R}^n$
which can be obtained as unions of elementary cubes of embedding dimension~$n$.
More precisely, we say that a finite collection~$\mathcal{K}$ of elementary
cubes with embedding dimension~$n$ is a \emph{cubical complex}, if for every
cube~$C \in \mathcal{K}$ all elementary cubes in its boundary~$\partial C$
are also contained in~$\mathcal{K}$. For a given cubical complex, we can
then talk about its associated chains. For this, let~$0 \le k \le n$ be
fixed. Then chains which consist exclusively of elementary cubes
in~$\mathcal{K}$ with dimension~$k$ are called \emph{$k$-chains}
of~$\mathcal{K}$, and the set of all $k$-chains is denoted
by~$C_k(\mathcal{K})$. One can easily see that~$C_k(\mathcal{K})$
is an additive group, where again we use coefficients in~$F = \mathbb{Z}_2$,
and we will therefore refer to~$C_k(\mathcal{K})$ as the \emph{$k^{\rm th}$ chain
group}. Furthermore, if the boundary operator~$\partial$ is applied to a
$k$-chain~$c$, then~$\partial c$ is clearly a $k-1$-chain. In fact,
one can show that the map~$\partial : C_k(\mathcal{K}) \to
C_{k-1}(\mathcal{K})$ is a group homomorphism. Since we consider
coefficients in the field~$F = \mathbb{Z}_2$, even more is true:
The chain groups~$C_k(\mathcal{K})$ are vector spaces over~$\mathbb{Z}_2$,
and the maps~$\partial : C_k(\mathcal{K}) \to C_{k-1}(\mathcal{K})$ are
all linear.

We are now finally in a position to describe \emph{cubical homology}.
For this, let~$\mathcal{K}$ denote a cubical complex. Intuitively 
speaking, homology is used to count the number of $k$-dimensional
holes in the cubical complex, for every $0 \le k \le n$. As a first
step towards this goal, we call a $k$-chain $c \in C_k(\mathcal{K})$
a \emph{$k$-cycle}, if its boundary vanishes, i.e., if we have $\partial
c = 0$. The set of all $k$-cycles is denoted by~$Z_k(\mathcal{K})$,
and one can easily see that it is a subgroup of the $k^{\rm th}$ chain group.
For example, in Figure~\ref{fig:cubhom}\emph{(b)}, two $1$-cycles are
shown in red and green. It is clear from this image, that cycles alone
cannot be used to detect holes in the cubical complex~$\mathcal{K}$. While
the red $1$-cycle in the figure encloses a hole, the green one does not.

In order to distinguish between cycles which enclose holes and cycles
which do not, we need to introduce one final concept. For this, assume
as before that~$\mathcal{K}$ is a cubical complex and that $0 \le k \le n$.
Then a $k$-cycle is called a \emph{$k$-boundary}, if it is equal to the
boundary of some $k+1$-chain. We denote the set of all $k$-boundaries
by~$B_k(\mathcal{K})$, and since this set forms a subgroup of the $k^{\rm th}$
chain group, it is called the \emph{$k^{\rm th}$ boundary group}. In fact, 
due to the fundamental property~$\partial\partial = 0$ of the boundary
operator, one has the inclusion $B_k(\mathcal{K}) \subset Z_k(\mathcal{K})$,
i.e., the $k$-boundaries form a subgroup of the $k$-cycles.
We now return to the two $1$-cycles shown in
Figure~\ref{fig:cubhom}\emph{(b)}. Notice that the green cycle is 
the boundary of the $2$-chain which consists of all elementary cubes
in the interior of the green loop. In other words, the green cycle is
bounding a region of $2$-dimensional cubes in the cubical complex, 
which completely fill its interior. On the other hand, one can show
that the red $1$-cycle cannot be represented as the boundary of a
$2$-chain, and this cycle does enclose a hole in the cubical set.

The last example captures the essence of homology. In order to detect
holes in a cubical complex, one needs to consider $k$-cycles, as they are
the $k$-chains without boundary. On the other hand, one needs to make sure
that such a $k$-cycle actually encloses a hole, and for this one has to make
sure that the cycle is not the boundary of a $k+1$-chain in~$\mathcal{K}$.
Otherwise, this $k+1$-chain would fill up the holes created by the cycle.
Notice, however, that just counting the $k$-cycles which are not
$k$-boundaries does not give the correct number of holes. In fact, in
Figure~\ref{fig:cubhom}\emph{(b)} one can easily find $1$-cycles which
are different from the red one, but which still enclose the same
hole. While not quite as intuitive, one can show that two cycles enclose
the same holes, if their difference happens to be a $k$-boundary.

From a mathematical point of view, counting the number of holes seems
to be related to counting equivalence classes. Call two $k$-cycles
{\em homologous\/}, if their difference is contained in the boundary
group~$B_k(\mathcal{K})$. This defines an equivalence relation on the
set of all $k$-cycles. Furthermore, a cycle is is homologous to zero
if and only if it is a $k$-boundary, i.e., if it does not enclose a
hole. Based on this discussion, it should not come as a surprise that
the quotient groups
\begin{displaymath}
  H_k(\mathcal{K}) = Z_k(\mathcal{K}) / B_k(\mathcal{K})
  \qquad\mbox{ for all }\qquad
  0 \le k \le n
\end{displaymath}
capture information about $k$-dimensional holes in the cubical
complex~$\mathcal{K}$. These groups are called \emph{homology
groups}, and since we consider coefficients in the field
$F = \mathbb{Z}_2$, they are in fact quotient vector spaces.

The homology groups defined above are abstract algebraic 
objects which can be assigned to a cubical complex. But what
exactly is their relation to the number of $k$-dimensional
holes in~$\mathcal{K}$. To see this, notice that if we have
two $k$-cycles~$c_1$ and~$c_2$ which enclose to different holes
in the cubical complex, then the equivalence class
in~$H_k(\mathcal{K})$ determined by their sum~$c_1 + c_2$
encloses both holes. In other words, two different holes 
in~$\mathcal{K}$ give rise to three different equivalence 
classes in the $k^{\rm th}$ homology group. In fact, what we need
are the \emph{independent} classes that can be found
in~$H_k(\mathcal{K})$, i.e., we need to determine the size
of a basis in the quotient vector space~$H_k(\mathcal{K})$,
which is the same as the rank of the quotient group. If we
now define
\begin{displaymath}
  \beta_k(\mathcal{K}) =
  \mbox{rank } H_k(\mathcal{K}) \; ,
  \qquad\mbox{ called the \emph{$k^{\rm th}$ Betti number} of }
  \mathcal{K} \; ,
\end{displaymath}
then the number of $k$-dimensional holes in the cubical 
complex is given by~$\beta_k(\mathcal{K})$. Furthermore,
the $k^{\rm th}$ homology group is isomorphic
to~$\mathbb{Z}_2^{\beta_k(\mathcal{K})}$. Intuitively,
the zero-dimensional Betti number counts the number of
connected components of the cubical complex, while the
one-dimensional Betti number counts the number of tunnels
or one-dimensional holes. Similarly, the two-dimensional Betti
number counts the number of voids or cavities in the complex,
and this classification of holes continues for higher dimensions.
For example, the cubical complex in Figure~\ref{fig:cubhom} is
connected and has one hole, but as a flat object no cavities.
This implies both~$\beta_0(\mathcal{K}) = 1$
and~$\beta_1(\mathcal{K}) = 1$, while~$\beta_2(\mathcal{K}) = 0$.

The cubical homology theory introduced above provides a very
rough quantitative measurement for the topological complexity of
a complex~$\mathcal{K}$ and its underlying topological space.
To make this more precise, let~$|\mathcal{K}| \subset \R^n$
denote the union of all elementary cubes in~$\mathcal{K}$.
Then one can show that the information provided by the homology
groups is \emph{invariant} under continuous transformations
of~$|\mathcal{K}|$. Specifically, if two cubical complexes~$\mathcal{K}_1$ 
and~$\mathcal{K}_2$ are such that their underlying topological
spaces~$|\mathcal{K}_1|$ and~$|\mathcal{K}_2|$ are homeomorphic,
then the homology groups~$H_k(\mathcal{K}_1)$
and~$H_k(\mathcal{K}_2)$ are isomorphic. In fact, this is 
still true even if the spaces~$|\mathcal{K}_1|$ and~$|\mathcal{K}_2|$
are only homotopy equivalent. In other words, homology captures 
connectivity and hole information that is invariant under
continuous deformations.
\subsection{Persistent Homology}
\label{sec:intropers}
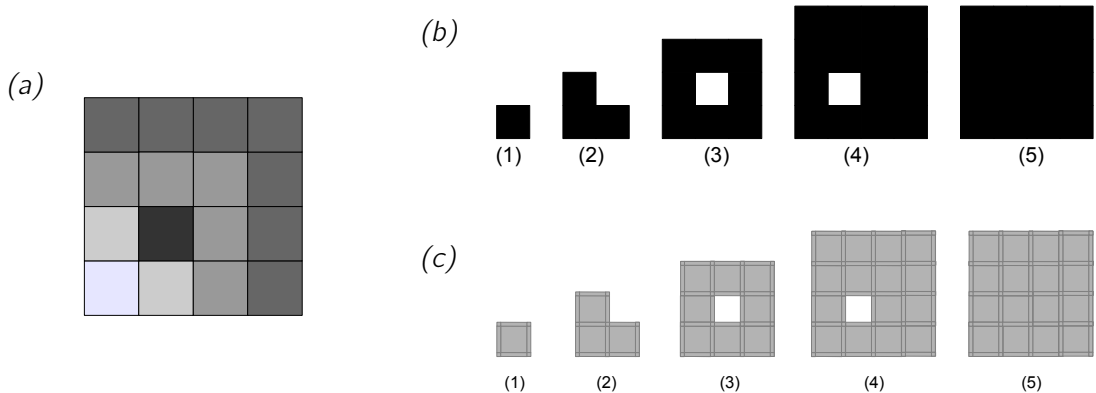
\begin{figure} \centering
  \setlength{\unitlength}{1 cm}
  \begin{picture}(14.5,5.5)
    \put(0.0,4.0){
      \parbox{1cm}{\emph{(a)}}}
    \put(1.0,1.0){
      \includegraphics[width=3cm]{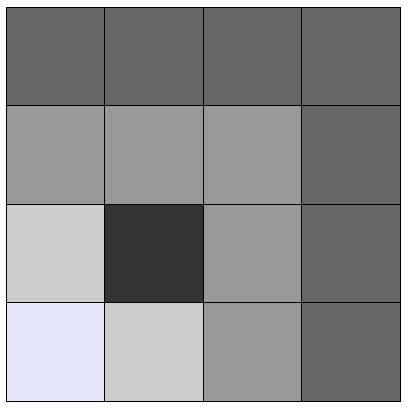}}
    \put(5.5,4.7){
      \parbox{1cm}{\emph{(b)}}}
    \put(6.5,3.0){
      \includegraphics[width=8cm]{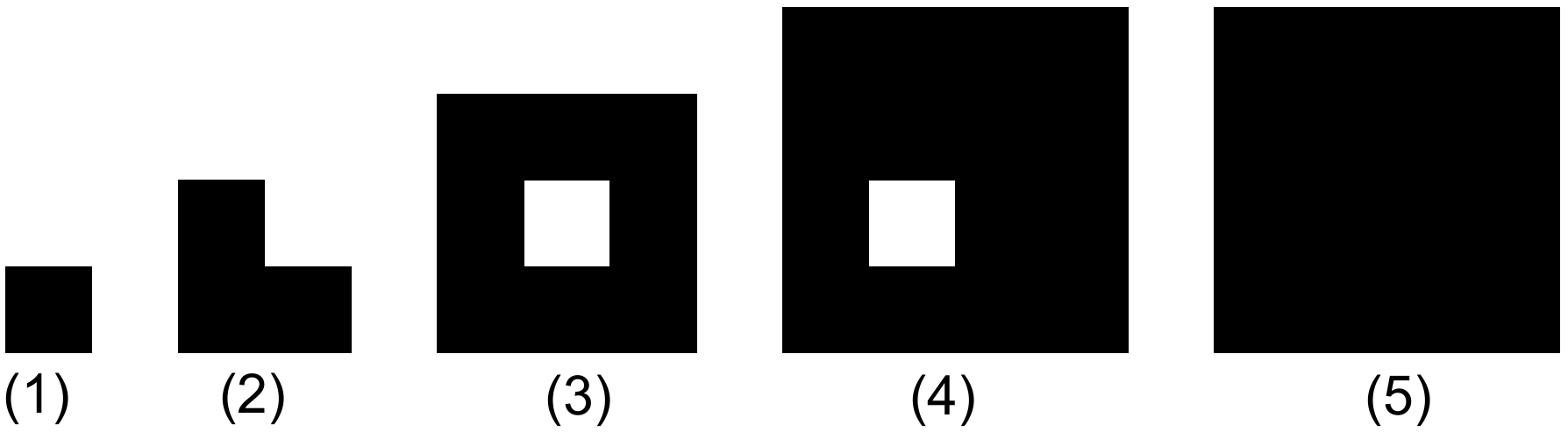}}
    \put(5.5,1.7){
      \parbox{1cm}{\emph{(c)}}}
    \put(6.5,0.0){
      \includegraphics[width=8cm]{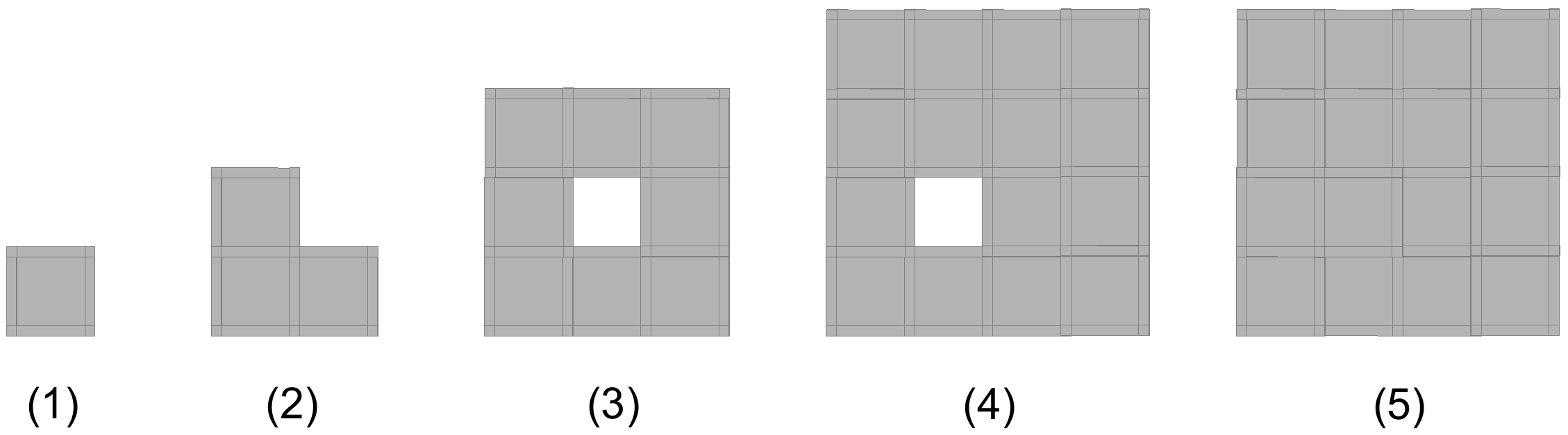}}
  \end{picture}
  \caption{Cubical complexes and filtration assigned to a gray-scale
           image. The image is shown in~\emph{(a)} and consists of sixteen
           pixels. From this image, one can extract cubical subcomplexes
           by considering only pixels with gray-scale intensity less than
           a given threshold. Depending on the threshold, this leads to one
           of the five cubical complexes shown in~\emph{(b)}. Alternatively,
           one can assign one filtration of cubical complexes to the gray-scale
           image, as shown in~\emph{(c)}. This filtration consists of five
           nested cubical complexes.}
  \label{fig:filtrationOfAnImage}
\end{figure}

The homology groups introduced in the last section are topological
invariants which can be assigned to a fixed topological space~$|\mathcal{K}|$.
In applications, finding the ``correct'' space to study is often a
straightforward task, but sometimes the situation is less clear.
Consider for example a gray-scale image, and suppose we are trying
to use homology to detect holes in the image. One simple instance
of this is shown in Figure~\ref{fig:filtrationOfAnImage}\emph{(a)}.
In this $4 \times 4$-pixel image, the lighter-colored pixels seem
to enclose a hole given by the solid black pixel. If we are
trying to use homology to detect this hole, we need to create
a cubical complex which consists of the ``correct pixels'' of the
image. This can be accomplished for example by choosing a gray-scale
threshold, and then defining~$\mathcal{K}$ to consist of all pixels
with intensity less than that threshold, together with all their 
edges and vertices. Depending on the threshold, one obtains the 
cubical complexes shown in Figure~\ref{fig:filtrationOfAnImage}\emph{(b)},
and the third and fourth of these complexes do correctly identify the
hole visible in the image. Note, however, that the remaining three
complexes fail to detect the hole.

But what about a situation in which the image contains two holes,
whose interiors occur at different gray-scale levels? Depending on
the intensities of the respective surrounding pixels, it is conceivable
that one cannot find a threshold value which detects both holes at the
same time. In situations like this, finding one threshold level for the
whole image clearly is the wrong goal. Rather, it would be nice to have
an automatic way which allows for the detection of holes at different 
levels --- and using homology in its original form alone makes this a
difficult task.

To address problems of the above type, the concept of \emph{persistent
homology} has been proposed, and it will be described in more detail
in the following. For a literature review on the subject, we refer the
reader to~\cite{herbert} and the references therein. The fundamental
idea behind persistent homology can easily be motivated by the complexes
shown in Figure~\ref{fig:filtrationOfAnImage}. Rather than singling out
a specific cubical complex in~\emph{(b)}, persistent homology deals
with all of these complexes using the concept of a \emph{filtration}.
In the situation of the figure, the underlying filtration is the 
increasing (with respect to set inclusion) sequence of the five cubical
complexes shown in~\emph{(c)}, indexed by the threshold levels indicated
by~(1) through~(5). In other words, a filtration encodes how the full
image is assembled by adding pixels according to their gray-scale level,
where pixels of the same level will be added simultaneously. Even though
we only mention the pixels, which are two-dimensional cubes, we implicitly
assume that at every stage of the filtration we do have cubical complexes
as defined in the last section. This means that the edges and vertices of the
pixels have to be added at the same time as the pixel. This is indicated
in Figure~\ref{fig:filtrationOfAnImage}\emph{(c)} by highlighting the
one- and zero-dimensional cubes. This implies for example that at stage~(1)
the shown cubical complex has four vertices, four edges, and one square,
while at stage~(2) we have eight vertices, ten edges, and three squares.
The goal of persistent homology is to detect and keep track of holes in
this sequence of cubical complexes, as they grow from complex~(1) to
complex~(5). In the simple example of Figure~\ref{fig:filtrationOfAnImage},
persistent homology will allow us to deduce that a single one-dimensional
hole can be observed starting with the third complex, and that it will
disappear again at level~(5).

The example involving gray-scale images can easily be generalized. Formally,
a \emph{filtration of cubical complexes} is a nested sequence of cubical
complexes
\begin{equation} \label{def:filtration}
  \mathcal{K}_1 \subset \mathcal{K}_2 \subset \ldots \subset
  \mathcal{K}_m \; .
\end{equation}
Notice that since we require every~$\mathcal{K}_i$ to be a cubical complex,
if a cube of dimension~$k$ appears for the first time in~$\mathcal{K}_i$,
all of its boundary elements have to appear at the latest in~$\mathcal{K}_i$,
but they could appear earlier. While filtrations could be specified explicitly,
one often convenient, and in fact equivalent, way of defining them is based
on the concept of a \emph{filtering function}. This is a function~$f :
\mathcal{K}_m \to \{ 1, \ldots, m \}$ which assigns every cube~$C$ in the final
cubical complex an integer between~$1$ and~$m$ in such a way that $f(C) = i$
if and only if $C \in \mathcal{K}_i$ and $C \not\in \mathcal{K}_{i-1}$,
where we define~$\mathcal{K}_0 = \emptyset$. One can readily see that for
this filtering function we have
\begin{equation} \label{def:filteringfcomplex}
  \mathcal{K}_i = \left\{ C \in \mathcal{K}_m \; : \;
    f(C) \le i \right\}
  \qquad\mbox{ for all }\qquad
  1 \le i \le m \; .
\end{equation}
Equivalently, but more useful for practical applications, assume
that~$\mathcal{K}$ denotes a cubical complex and that we are given a
function~$f : \mathcal{K} \to \{ 1, \ldots, m \}$ such that the following
holds: If~$A,B \in \mathcal{K}$ are any two cubes such that~$B$ is contained
in the boundary of~$A$, then the function~$f$ satisfies $f(A) \ge f(B)$. In
this case, we call~$f$ a \emph{filtering function} on~$\mathcal{K}$,
and~(\ref{def:filteringfcomplex}) defines a filtration of cubical complexes
with~$\mathcal{K}_m = \mathcal{K}$. We would like to point out that the
above monotonicity property of~$f$ from a cube to its boundary is all that
is needed to ensure that every~$\mathcal{K}_i$ is indeed a cubical complex.

The above concept of a filtering function can easily be related to our image
example from Figure~\ref{fig:filtrationOfAnImage}. In this case, let~$\mathcal{K}$
denote the cubical complex consisting of all sixteen pixels, together with their
edges and vertices. Furthermore, suppose that the gray-scale intensities
are given by the five values~$a_1 < a_2 < \ldots < a_5$. Then for every
two-dimensional pixel~$A$ in~$\mathcal{K}$ we define~$f(A) = i$, if its
gray-scale intensity is given by~$a_i$. For an edge~$B \in \mathcal{K}$,
we let~$f(B)$ be the smallest value of~$f$ on a pixel which contains~$B$ as
an edge, and for a vertex~$C \in \mathcal{K}$, the value~$f(C)$ is the
smallest value of~$f$ on an edge which contains~$C$. Then~(\ref{def:filteringfcomplex})
leads to the filtration shown in Figure~\ref{fig:filtrationOfAnImage}\emph{(c)}.
This example also shows that choosing~$\{ 1, \ldots, m \}$ as the index
set for a filtration of cubical complexes can be done without loss of generality
--- through a simple transformation, any finite and discrete subset of~$\mathbb{R}$
can be mapped onto this set. Alternatively, we could keep the values of the
filtering function in the set~$\{ a_1, \ldots, a_5 \}$, and index the 
cubical complexes by these distinct real numbers. This extension will be
useful later on.

Having introduced the concept of a filtration, we now turn our attention
to persistent homology. Suppose we are given a filtration of cubical
complexes as in~(\ref{def:filtration}). Then for each dimension~$k$ and
each complex in the sequence we obtain a homology group~$H_k(\mathcal{K}_i)$.
Note that these homology groups are in general not ordered by inclusion,
despite the fact that the underlying cubical complexes are. However, if~$c$
denotes an arbitrary $k$-cycle in~$\mathcal{K}_i$, then one can readily see 
that~$c$ is also a $k$-cycle in~$\mathcal{K}_{i+1}$. This implies that we
can define a map from~$H_k(\mathcal{K}_i)$ to~$H_k(\mathcal{K}_{i+1})$ by
mapping the equivalence class of~$c$ in the first homology group to the
equivalence class of~$c$ in the second one. In order words, the inclusions
in~(\ref{def:filtration}) lead to a sequence of maps
\begin{displaymath}
  H_k(\mathcal{K}_1) \rightarrow H_k(\mathcal{K}_2) \rightarrow
  \ldots \rightarrow H_k(\mathcal{K}_m) \; ,
\end{displaymath}
and these maps are in fact group homomorphisms. 

It was shown in the last section that if~$\gamma \in H_k(\mathcal{K}_i)$
is a nontrivial homology class, and if~$c$ is a $k$-chain representing
this class, then~$c$ encloses a $k$-dimensional hole in~$\mathcal{K}_i$.
Now consider the class~$\tilde{\gamma}$ induced by~$c$
in~$H_k(\mathcal{K}_{i+1})$. If this equivalence class is nontrivial,
then~$c$ still encloses a hole in~$\mathcal{K}_{i+1}$. On the other hand,
if the equivalence class~$\tilde{\gamma}$ is trivial, then the hole has
been filled in by cubes that were added to~$\mathcal{K}_i$ during the
formation of~$\mathcal{K}_{i+1}$. In other words, the homomorphisms
induced by inclusion allow us to detect the disappearance of holes.
Note, however, that holes can also disappear in a second way. Rather
than being mapped to a trivial homology class, it is also possible
that two different homology classes~$\gamma_1, \gamma_2 \in
H_k(\mathcal{K}_i)$ are mapped to the same nontrivial class
in~$H_k(\mathcal{K}_{i+1})$, and this case will be discussed in more
detail later.

What about the appearance of holes? Suppose again that~$\gamma \in
H_k(\mathcal{K}_i)$ is a nontrivial homology class representing a
$k$-dimensional hole in~$\mathcal{K}_i$. One can easily see that if
this hole was already present in~$\mathcal{K}_{i-1}$, then~$\gamma$
has to be in the image of the inclusion-induced homomorphism
from~$H_k(\mathcal{K}_{i-1})$ to~$H_k(\mathcal{K}_i)$. In other
words, the hole~$\gamma$ appears for the first time in~$\mathcal{K}_i$,
if and only if it is not in the range of this homomorphism.

Being able to track the appearance and disappearance of holes
in the filtration is the main idea behind \emph{persistent homology}.
While it is usually defined via the so-called \emph{persistent homology
groups}, which are just the images of the above-mentioned maps between
homology groups induced by inclusion, we proceed right away to an
equivalent formulation which is based on the \emph{lifespan of
homology classes}. We say that a $k$-dimensional homology class~$\alpha$
is \emph{born} in the complex~$\mathcal{K}_j$, if we have both $\alpha
\in H_k(\mathcal{K}_j)$ and $\alpha \not\in \iota(H_k(\mathcal{K}_{j-1}))$,
where~$\iota : H_k(\mathcal{K}_{j-1}) \to H_k(\mathcal{K}_j)$ is induced
by inclusion. Moreover, we say that a $k$-dimensional homology class~$\alpha$
\emph{dies} in the complex~$\mathcal{K}_\ell$, if~$H_k(\mathcal{K}_\ell)$
is the first homology group where~$\alpha$ either became trivial, or where
it became identical to another class which was born earlier in the 
filtration. Note that the first of these cases, corresponds to the 
disappearance of a hole discussed above, while the second case addresses
the issue of two nontrivial homology classes merging --- in this case,
the class that was created later in the filtration dies, while the one
which was created earlier survives. In other words, merging effects can be
treated in a unique way, through a process called the \emph{elder rule}.

Using the definitions of the last paragraph, persistent homology provides
us with a list of nontrivial and independent homology classes, together 
with their birth times and death times. If~$\gamma$ is one of these classes
and if its birth and death times are given by~$j$ and~$\ell$, respectively,
then the interval~$[j,\ell)$ is called the \emph{persistence interval} 
associated with~$\gamma$, and the integer~$\ell - j$ is called its
\emph{lifespan}. If there are classes which are born, but never die,
then we set $\ell = \infty$. In typical applications, persistent homology
classes with relatively long lifespan are considered to be more important
than the ones that have shorter lifespan.
\begin{figure} \centering
  \setlength{\unitlength}{1 cm}
  \begin{picture}(12,7)
    \put(0.0,0.0){
      \includegraphics[height=7cm]{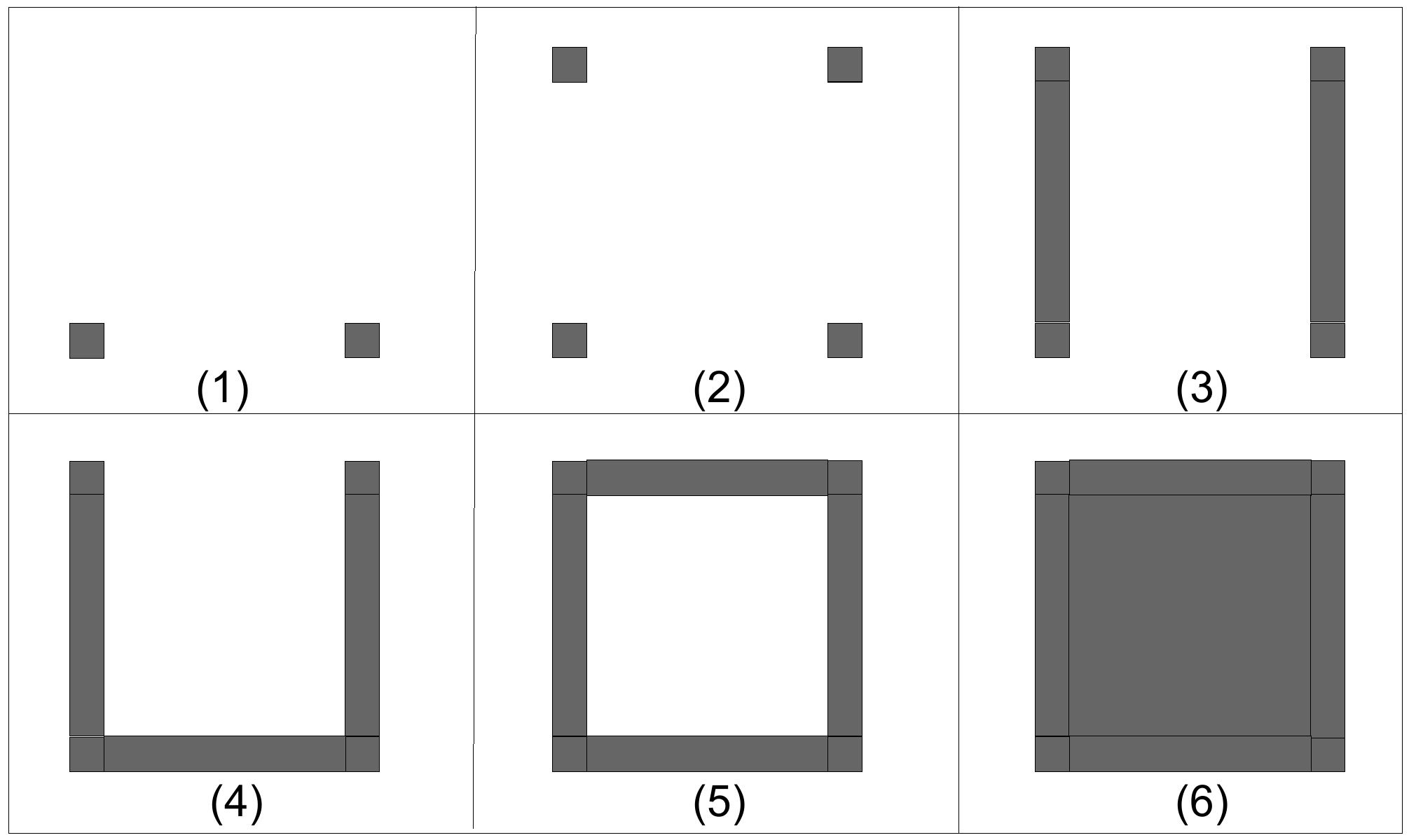}}
  \end{picture}
  \caption{Illustration of persistent homology. The six images show a
           filtration of cubical complexes in the plane, which leads
           to the following homology classes and lifespans.
           (1)~In the complex~$\mathcal{K}_1$, two connected components
           are born, i.e., two nontrivial zero-dimensional homology classes.
           (2)~In~$\mathcal{K}_2$, another two connected components are born.
           (3)~In the complex~$\mathcal{K}_3$, the two connected components
           born in~$\mathcal{K}_2$ die, since each of them merges with a 
           component which was born in~$\mathcal{K}_1$ --- these latter two 
           components still live on.
           (4)~In~$\mathcal{K}_4$, however, the remaining two components
           merge, so one dies and the other survives. Since both were born 
           at the same time, one can randomly choose which one survives.
           (5)~In the complex~$\mathcal{K}_5$ a $1$-dimensional cycle is
           created, which gives rise to a nontrivial homology class.
           (6)~This $1$-dimensional class dies in the final
           complex~$\mathcal{K}_6$.}
  \label{fig:persistence}
\end{figure}

The lifespans of homology classes are illustrated in
Figure~\ref{fig:persistence}. As described in the caption, this figure
shows a filtration of six cubical complexes in the plane. Persistent
homology yields five different homology classes. Four of these
are in dimension zero, and their persistence intervals are given
by~$[1,\infty)$, $[1,4)$, $[2,3)$, and again~$[2,3)$, respectively.
The one nontrivial class in dimension one has persistence interval~$[5,6)$.

The persistence intervals provide information about when homology
classes are born, evolve, and die in the filtration. But they also
can be used to recover the homologies of the involved cubical complexes.
It can be shown that the $k^{\rm th}$ Betti number of the complex~$\mathcal{K}_i$
in the filtration is given by the number of persistence intervals in
dimension~$k$ which contain~$i$. For example, in the situation of
Figure~\ref{fig:persistence}, the complex~$\mathcal{K}_3$ has Betti
number~$\beta_0(\mathcal{K}_3) = 2$, since only the intervals~$[1,\infty)$
and~$[1,4)$ contain the index~$3$. We would like to pint out again, however,
that persistent homology provides considerably more information than the
sequence of Betti numbers of the complexes~$\mathcal{K}_i$, as it allows
us to track for how long specific homology classes survive.
\begin{figure} \centering
  \setlength{\unitlength}{1 cm}
  \begin{picture}(14,7)
    \put(0.0,0.0){
      \includegraphics[height=7cm]{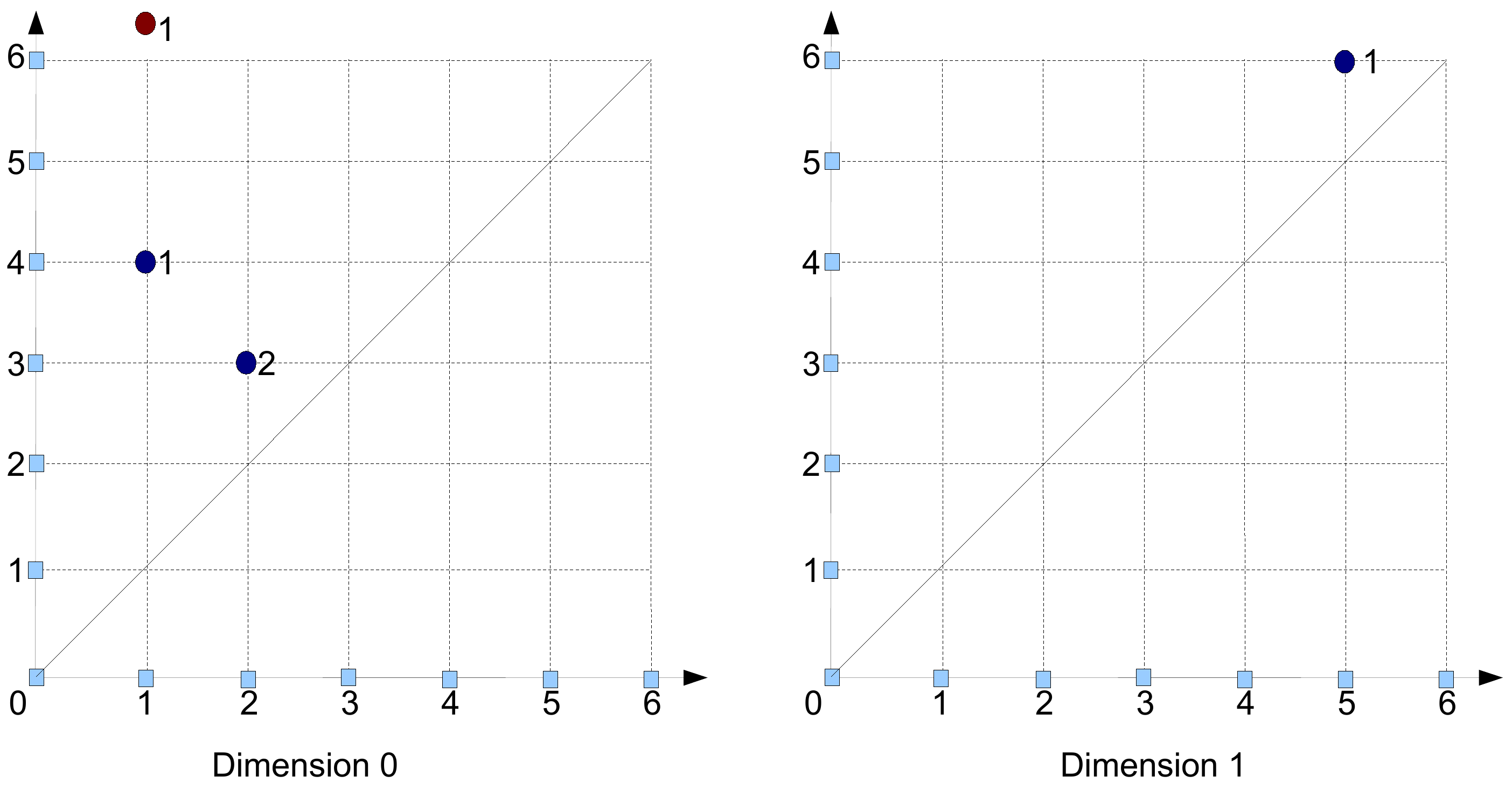}}
  \end{picture}
  \caption{Persistence diagrams for the filtration shown in
           Figure~\ref{fig:persistence}. The left image shows the 
           persistence diagram in dimension zero, which contains the
           points~$(1,\infty)$ and $(1,4)$ once, and the
           point~$(2,3)$ twice. These multiplicities are indicated
           by the numbers adjacent to the points. The right image contains
           the diagram for dimension one, which only contains the single
           point~$(5,6)$.}
  \label{fig:persistenceDiagram}
\end{figure}

While the persistence intervals encode all of the information
provided by persistent homology, dealing with the intervals directly
is fairly cumbersome. In practice, these intervals are visualized
in a \emph{persistence diagram} as follows. Fix a dimension~$k$. Then
the persistence diagram in dimension~$k$ is a finite collection of 
points (possibly of multiplicity larger than one) in the plane such
that~$(a,b) \in \mathbb{R}^2$ is a point in the $k^{\rm th}$ diagram if
and only if persistent homology provides a nontrivial $k$-dimensional
homology class with interval~$[a,b)$. In other words, a persistence
diagram is a multiset of points in~$\mathbb{R}^2$ which corresponds
to the persistence intervals. Notice that all of these points lie
above the diagonal in the first quadrant of~$\mathbb{R}^2$. For the
filtration in Figure~\ref{fig:persistence} the persistence diagrams
in dimensions zero and one are shown in Figure~\ref{fig:persistenceDiagram}.

For our special situation of a filtration indexed by integers as
in~(\ref{def:filtration}), the points in a persistence diagram always
have integer coordinates. However, if one considers a cubical
complex~$\mathcal{K}$ together with a filtering function $f : \mathcal{K}
\to \{ a_1, \ldots, a_m \}$, one can consider the associated filtration
\begin{equation} \label{def:filteringfcomplex2}
  \mathcal{K}_{a_1} \subset \mathcal{K}_{a_2} \subset \ldots \subset
    \mathcal{K}_{a_m}
  \qquad\mbox{ defined via }\qquad
  \mathcal{K}_{a_i} = \left\{ C \in \mathcal{K} \; : \; f(C) \le a_i
    \right\} \; .
\end{equation}
In this more general setting, persistent homology produces persistence
diagrams in a completely analogous way, but now with points of the
form~$(a_j, a_\ell) \in \R^2$ which lie above the main diagonal. As 
mentioned before, the values~$a_j$ could be the gray-scale intensities
in an underlying image. What happens if these values are slightly 
perturbed? One would expect that new persistence diagram is only a
slight perturbation of the old one, provided the perturbations are 
small enough.

\begin{figure} \centering
  \setlength{\unitlength}{1 cm}
  \begin{picture}(17.0,5.0)
    \put(0.0,0.0){
      \includegraphics[height=5cm]{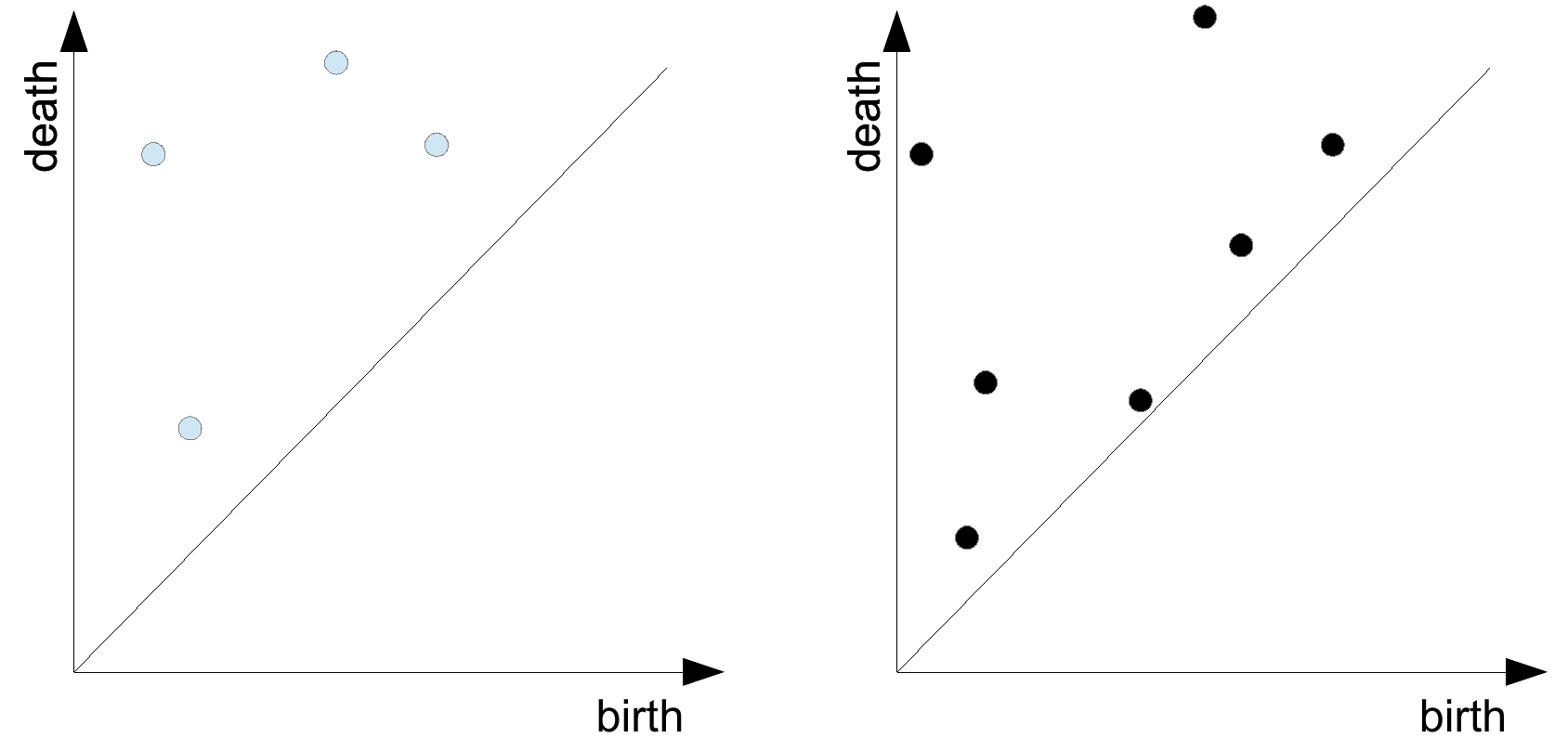}}
    \put(12.0,0.0){
      \includegraphics[height=5cm]{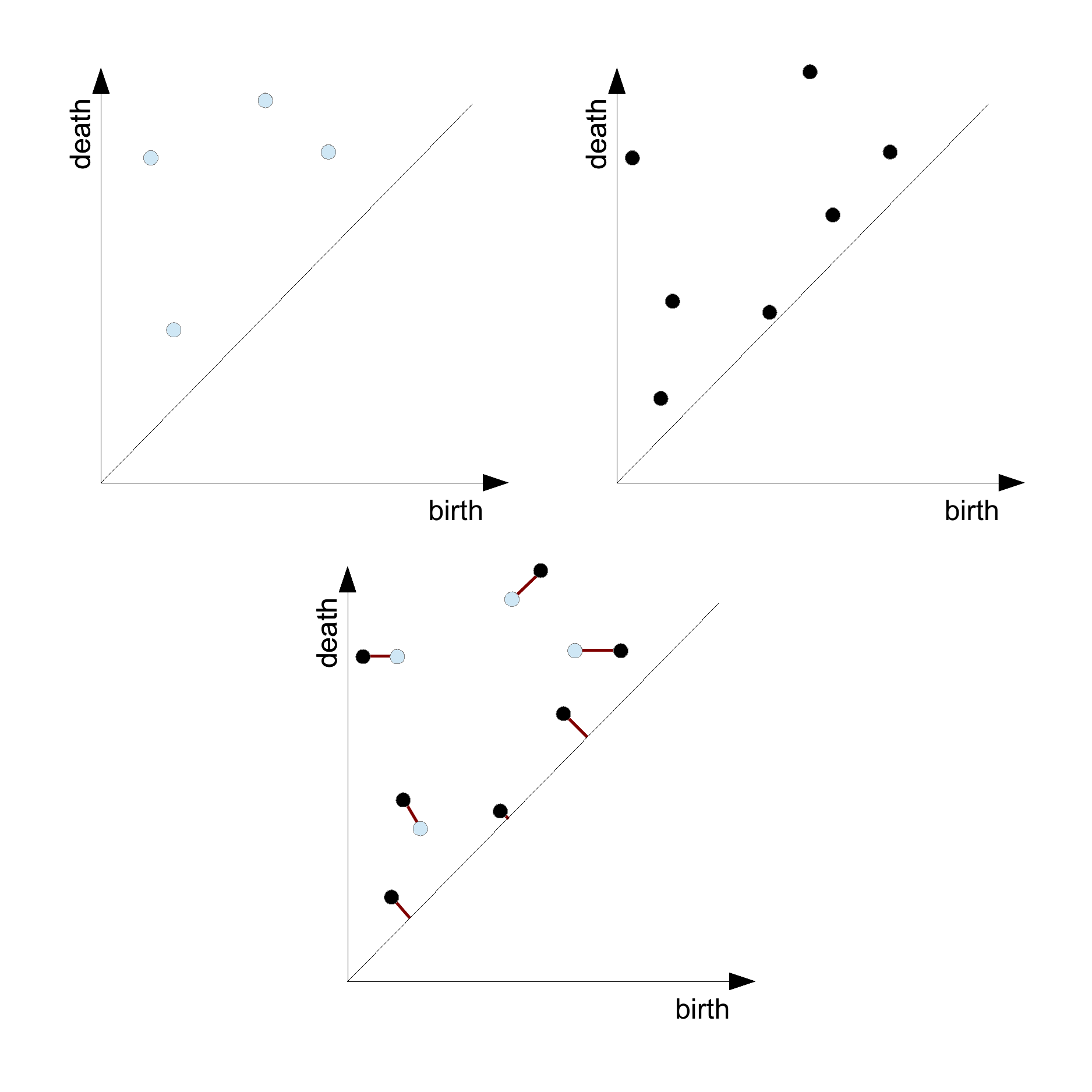}}
  \end{picture}
  \caption{For the persistence diagrams~$D_1$ and~$D_2$ shown in the first
           two images, the third image contains a matching $b : \tilde{D}_1
           \to \tilde{D}_2$ in which points have to be moved only over
           fairly short distances.}
  \label{fig:bottelneckDistance}
\end{figure}

One of the main features of persistent homology is that this 
perturbation idea can be quantified in a most satisfying way. 
In order to do this, one needs to introduce a metric on the space
of persistence diagrams, and there are a number of ways for doing
this~\cite{herbert}. For the purposes of this paper, we focus on the
so-called \emph{bottleneck distance}, which is closely related to the
earth mover's distance discussed in~\cite{earthMoverDistance}. To
define this distance, consider two persistence diagrams~$D_1$ and~$D_2$,
both in the same dimension. Note, however, that these persistence diagrams
could be created from two completely different filtrations. Let us now
assume that in addition to the finitely many points corresponding to
persistence intervals, we add all points of the form~$(x,x)$ with~$x
\in \mathbb{R}$ to both diagrams, each of these points with infinite
multiplicity. While this might seem strange at first glance, it is
necessary for comparing persistence diagrams with different number 
of points. We denote the resulting new diagrams by~$\tilde{D}_1$
and~$\tilde{D}_2$, respectively. After these preparations, we can
finally define the bottleneck distance of the diagrams~$D_1$ and~$D_2$,
which will be denoted by~$W_{\infty}(D_1,D_2)$ and is defined as
\begin{displaymath}
  W_{\infty}(D_1,D_2) \; = \;
  \inf_{\mbox{\footnotesize bijections } b : \tilde{D}_1 \rightarrow \tilde{D}_2}
    \left( \sup_{x \in \tilde{D}_1} \| x - b(x) \| \right) \; .
\end{displaymath}
In other words, the bottleneck distance tries to find a matching
$b : \tilde{D}_1 \rightarrow \tilde{D}_2$ which minimizes the largest
distance between a point~$x \in \tilde{D}_1$ and its image~$b(x) \in \tilde{D}_2$.
Note that obtaining matchings of this type is made possible only through
extending the persistence diagrams by the points on the diagonals, since~$D_1$
and~$D_2$ could have very different numbers of points. Since points on the
diagonal would correspond to homology classes which die immediately when
born, this should not alter the information content of the persistence
diagrams. Intuitively, the bottleneck distance between two diagrams~$D_1$
and~$D_2$ is of size~$\epsilon$, if it is possible to move every point
of~$D_1$ to a point in~$D_2$ or the diagonal by less than distance~$\epsilon$,
as well as vice versa. This is illustrated in Figure~\ref{fig:bottelneckDistance}.
As we mentioned before, typically, persistent homology classes with relatively
long lifespan are considered to be more important than the ones that have shorter
lifespan --- and the latter ones give rise to points in the persistence diagram
which are closer to the diagonal. This implies that these are the points that
can easily be moved onto the diagonal in a matching.

We now return to filtrations which are generated through a filtering
function~$f : \mathcal{K} \to \{ a_1, \ldots, a_m \}$ as defined
in~(\ref{def:filteringfcomplex2}). How do perturbations of~$f$ affect
the resulting persistence diagrams? To describe this, let us consider
two filtering functions~$f_1$ and~$f_2$ on the same cubical
complex~$\mathcal{K}$, but we do not assume that the finite ranges
of~$f_1$ and~$f_2$ are the same. If we further denote the persistence
diagrams in dimension~$k$ associated with~$f_1$ and~$f_2$
by~$D_k(\mathcal{K},f_1)$ and~$D_k(\mathcal{K},f_2)$, respectively,
then one can prove a \emph{stability theorem} in the form of the
estimate
\begin{equation} \label{th:stabilityTheorem}
  W_{\infty} \left( D_k(\mathcal{K},f_1),
    D_k(\mathcal{K},f_2) \right)
  \; \leq \;
  \left\| f_1 - f_2 \right\|_{\infty} \; ,
\end{equation}
where~$\|f_1 - f_2\|_{\infty}$ denotes the usual maximum norm.
See~\cite{herbert} for more details. In other words, the bottleneck
distance of the two persistence diagrams is bounded by the maximal
function value difference of the two filtering functions.

The stability theorem encapsulated by~(\ref{th:stabilityTheorem})
is one of the central properties of persistence diagrams, which is
invaluable in applications. It implies, for example, that if observational
data is perturbed by noise which is bounded by a constant~$\epsilon$, then
the persistence diagram of the noisy data will be at most~$\epsilon$ away
from the persistence diagram of a true underlying data. In other words,
persistence is a robust topological metric. In the context of this paper,
this robustness property allows us to compute persistence diagrams of
finite-dimensional spectral approximations to the true solution~$u$
of the Cahn-Hilliard-Cook equation~(\ref{chc}), which are basically
indistinguishable from the true underlying persistence diagrams due
to spectral accuracy.

In closing, and for the sake of completeness, we would like to point out
that the bottleneck distance is only one possibility for measuring the
distance between persistence diagrams. In fact, in some situations its
insistence on only registering the largest move distance~$\| x - b(x) \|$
is too stringent. In such cases one can make use of the
\emph{Wasserstein distance}~\cite{herbert}, which is defined as
\begin{displaymath}
  W_p(D_1,D_2) \; = \;
  \inf_{\mbox{\footnotesize bijections } b : \tilde{D}_1 \rightarrow \tilde{D}_2}
    \left( \sum_{x \in \tilde{D}_1} \| x - b(x) \|^p \right)^{1/p}
\end{displaymath}
for $p \ge 1$. This distance has the advantage of taking all line segments in
a matching into account, but it comes at the price of no longer satisfying
the stability result~(\ref{th:stabilityTheorem}). However, more elaborate
stability estimates which are valid for the Wasserstein distance are available
in a number of situations, and we refer the reader to~\cite{herbert} for more
information.
\subsection{Persistence Landscapes}
\label{sec:intropl}
By moving from the concept of homology groups to the concept of persistent
homology, we have made the transition from thresholding a gray-scale image
at one specific level to considering the fine structure of the image as
we pass through all possible shades of gray. As we mentioned in the introduction,
the goal of this paper is to show that the topology of material microstructures
such as the ones in Figure~\ref{fig:chcpatt} encodes interesting facts about the
underlying dynamics. After suitable discretization these microstructures
can be represented as gray-scale images, with the value of the solution~$u$
of~(\ref{chc}) serving as the filtration parameter. Recall, however, that the
Cahn-Hilliard-Cook model is a stochastic partial differential equation, and
therefore with probability one no two realizations of the dynamics will lead
to identical patterns. In order to uncover this ``typical'' behavior, one would
therefore like to consider statistical properties, in particular, averages
of the topological information.

While averages of Betti numbers are easily accessible, and have in fact been
used for example in~\cite{gameiro:etal:05a}, it is not immediately clear
how the average of persistence diagrams can be computed. One of the first 
attempts at defining the mean of two persistence diagrams is the notion
of \emph{Fr\'echet mean} which was explored in~\cite{frechet1}. Without
going into too much detail, assume we are given two persistence diagrams~$D_1$
and~$D_2$, and assume that we have found a matching~$b : \tilde{D}_1 \to
\tilde{D}_2$ which realizes the $p$-Wasserstein distance~$W_p(D_1, D_2)$,
for some $p \ge 1$. Then the Fr\'echet mean of the diagrams~$D_1$ and~$D_2$
is defined as the collection of midpoints of the line segments between~$x$
and~$b(x)$, for all $x \in \tilde{D}_1$. Despite the fact that the Wasserstein
distance between~$D_1$ and~$D_2$ takes the length of all matching distances
$\|x - b(x)\|$ into account, one can easily see that the above notion of
Fr\'echet mean is not well-defined. Consider for example a diagram~$D_1$ 
which consists of the points~$(j,k)$ and~$(j+1,k+1)$, while the diagram~$D_2$
has the points~$(j,k+1)$ and~$(j+1,k)$, and where we assume that $j+1 < k$. Then
there are two minimal matchings which involve these points: One matching maps
along two horizontal line segments of length one, while the other one uses
two vertical line segments of length one. This implies that in the former
case the Fr\'echet mean consists of the points~$(j+1/2, k)$ and~$(j+1/2,
k+1)$, whereas in the latter case it is given by the points~$(j, k+1/2)$
and~$(j+1, k+1/2)$.

It is not surprising that this obvious non-uniqueness of the Fr\'echet
mean makes a statistical analysis of persistence diagrams which involves 
this concept nontrivial. From a mathematical point of view, these problems
can be resolved, see for example~\cite{frechet2}. Unfortunately, however,
these ideas have not yet been implemented algorithmically, as they pose
serious technical and computational difficulties. For the purposes of the
present paper we therefore pursue a different approach, which uses the
concept of \emph{persistence landscapes}. These objects were introduced
in~\cite{peterLandscapes}, and they prove to be an algorithmically feasible
tool for the statistical analysis of persistence diagrams.

In order to motivate the definition of persistence landscapes, we return
one final time to the Fr\'echet mean. Its non-uniqueness stems from the
fact that it operates on a discrete set of points. In some sense, this
situation strongly resembles the problem of computing the mean of two
integers within the ring of integers. Also in this case, the mean is
often not unique. In the integer case, however, the mean problem can be
solved easily by introducing rational numbers. We will see below that the
viewpoint behind persistence landscapes is very similar --- it is based on
the idea of embedding persistence diagrams into the larger space of piecewise
linear functions, where the mean is uniquely defined.
For the remainder of this section, we present an intuitive introduction
to persistence landscapes, rather than a formal definition. For a more
detailed discussion we refer the reader to~\cite{peterLandscapes}. Important
for our applications is the fact that persistence landscapes can be computed 
and manipulated efficiently, and algorithms for this can be found in~\cite{plt}.

\begin{figure} \centering
  \setlength{\unitlength}{1 cm}
  \begin{picture}(15.0,10.0)
    \put(0.0,0.0){
      \includegraphics[height=10cm]{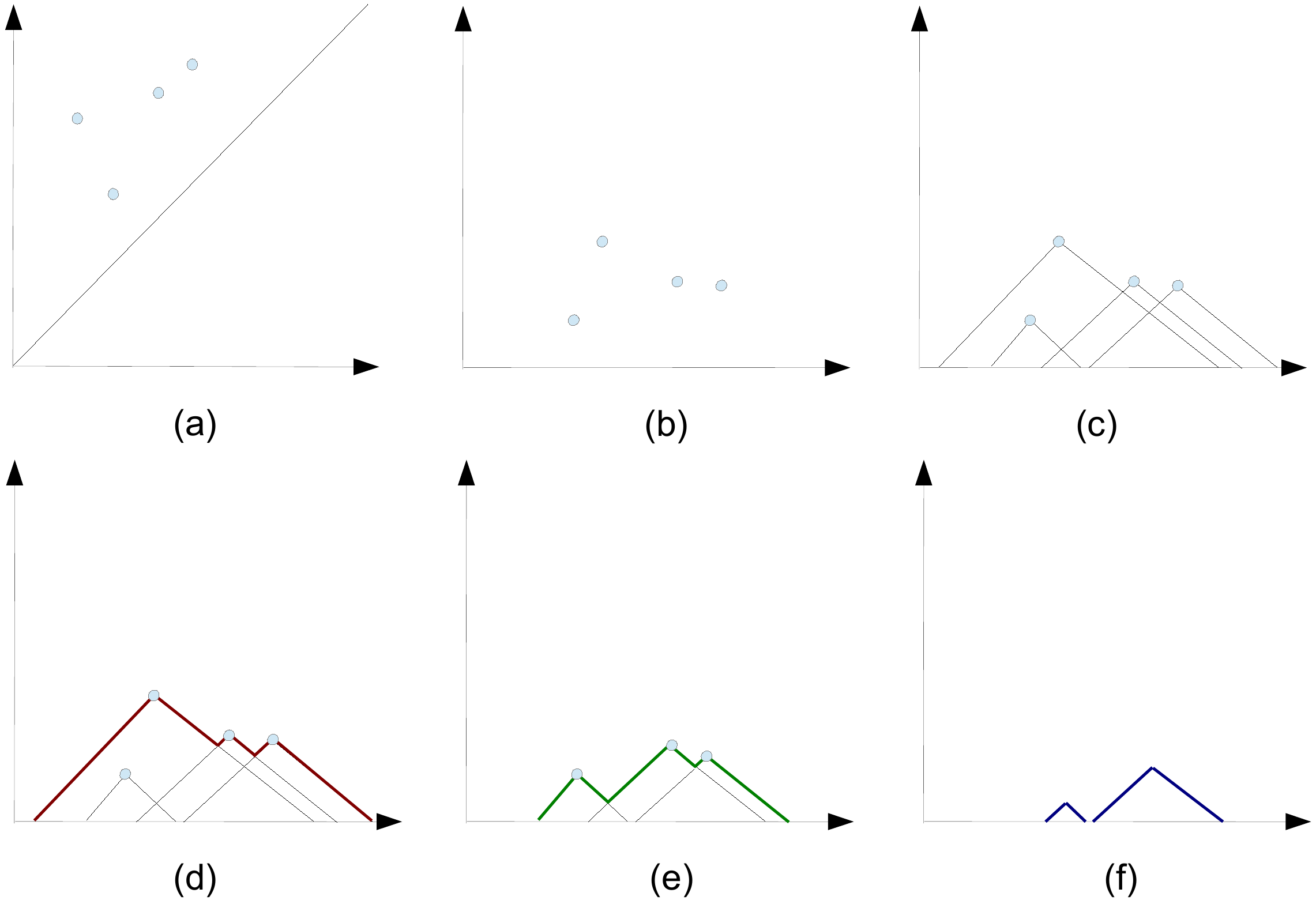}}
  \end{picture}
  \caption{The idea behind persistence landscapes:
           (a)~contains the initial persistence diagram, while
           (b)~shows the transformed diagram after mapping all points and the
               diagonal with the transformation~(\ref{def:pltransformation}).
           (c)~depicts the piecewise linear functions~$f_{(b,d)}$ defined
               in~(\ref{eq:basicLand}) for every one of the transformed
               points. Finally, the graph in
           (d)~illustrates the first persistence landscape function~$\lambda_1$
               in red, which is the upper envelope of the graphs from the
               previous picture.
           The second and third persistence landscape functions~$\lambda_2$
           and~$\lambda_3$ are shown in green and blue in~(e) and~(f),
           respectively, and they follow the second- and third-largest
           function values of the functions~$f_{(b,d)}$ in image~(c). All
           subsequent persistence landscape functions~$\lambda_k$ are
           identically zero.}
  \label{fig:persistenceLandscapeIllustration}
\end{figure}

To begin with, suppose we are given a point~$(b,d) \in \mathbb{R}^2$ in a
persistence diagram, i.e., the point corresponds to a homology class with
birth time~$b$ and death time~$d$. With this point, we associate a piecewise
linear function~$f_{(b,d)} : \mathbb{R} \rightarrow [0,\infty)$, which is
defined as
\begin{equation} \label{eq:basicLand}
  f_{(b,d)}(x) =
  \left\{ \begin{array}{ccl}
            0     & \mbox{ if } & x \not\in (b, d) \; , \\[2ex]
            x - b & \mbox{ if } & x \in \left( b, \frac{b+d}{2}
              \right] \; , \\[2ex]
            d - x & \mbox{ if } & x \in \left(\frac{b+d}{2},
              d \right) \; .
  \end{array} \right.
\end{equation}
In other words, the function~$f_{(b,d)}$ is a hat function, whose
maximal value is~$(d-b)/2$, i.e., half the lifespan of the homology
class represented by the point~$(b,d)$. Moreover, this global maximum
of~$f_{(b,d)}$ occurs at the center of the associated persistence interval,
and on either side the function decays linearly towards zero with slopes~$\pm 1$.
In fact, there is another, more geometric way of describing~$f_{(b,d)}$. For this,
map the point~$(b,d)$ in the persistence diagram via the linear map
\begin{equation} \label{def:pltransformation}
  \left( \begin{array}{c} b \\ d \end{array} \right)
  \; \mapsto \;
  \left( \begin{array}{rr} 1/2 & 1/2 \\ -1/2 & 1/2 \end{array} \right)
  \left( \begin{array}{c} b \\ d \end{array} \right) \; ,
\end{equation}
and then draw lines with slopes~$\pm 1$ from the image of the point
towards the horizontal axis. Clearly, the function~$f_{(b,d)}$ encodes
the original data provided by the point~$(b,d)$, but in a more convenient
way: The support of the function equals the closure of the persistence
interval~$[b,d)$, while the height of the function is half the lifespan.
Note also that the map~(\ref{def:pltransformation}) transforms the main
diagonal into the horizontal axis. This procedure is illustrated in
Figure~\ref{fig:persistenceLandscapeIllustration}\emph{(a)},
\emph{(b)}, and~\emph{(c)}.

After these preparations, we can now introduce the concept of persistence
landscape. As was mentioned earlier, in many applications a homology class
is considered important, if its lifespan~$d-b$ is large. Since the heights
of the functions~$f_{(b,d)}$ shown in
Figure~\ref{fig:persistenceLandscapeIllustration}\emph{(c)} are proportional
to the lifespan, the functions close to the top are therefore considered
more significant than the ones closer to the horizontal axis. Motivated
by this observation, the \emph{persistence landscape} of the birth-death
pairs~$(b_i , d_i)$, where $i = 1,\ldots,m$, which constitute the given
persistence diagram is the sequence of functions $\lambda_k : \mathbb{R}
\rightarrow [0,\infty)$ for $k \in \mathbb{N}$, where~$\lambda_k(x)$
denotes the $k^{\rm th}$ largest value of the numbers~$f_{(b_i,d_i)}(x)$,
for $i = 1, \ldots, m$, and we define $\lambda_k(x) = 0$ if $k > m$.
Equivalently, this sequence of functions can be combined into a single
function $L : \mathbb{N} \times \mathbb{R} \to [0,\infty)$ of two
variables, if we define $L(k,t) = \lambda_k(t)$. See again
Figure~\ref{fig:persistenceLandscapeIllustration} for an illustration
of these concepts.

One can easily see that the persistence landscape~$L$ encodes exactly the
same information as the underlying persistence diagram. Yet, in some 
sense, this information is presented in a way which is mathematically
more accessible. First of all, it creates a filtered representation of
the persistence diagram which is premised on the assumption that homology
classes with large lifespans are more important than ones with shorter
lifespans. This is reflected in part by the monotonicity property
\begin{equation} \label{Lmonotone}
  L(k, x) \ge L(k+1, x)
  \qquad\mbox{ for all }\qquad
  k \in \mathbb{N}
  \quad\mbox{ and }\quad
  x \in \mathbb{R} \; .
\end{equation}
Secondly, and most importantly, it represents the persistence data
as a piecewise linear function~$L$, and averaging
functions of this type results again in a piecewise linear function.
Moreover, the monotonicity property~(\ref{Lmonotone}) is transferred
from the involved landscapes to the average.

Thus, the introduction of persistence landscapes allows us to easily talk
about the average of persistence information. Assume we are given persistence
landscapes~$L_1, \ldots, L_\ell$. Then the \emph{averaged persistence landscape}
is defined by~$\sum_{i=1}^\ell L_i / \ell$, and it is again a piecewise linear
function which satisfies~(\ref{Lmonotone}). In addition, it is clearly uniquely
defined. But considering the persistence information in the form of a real-valued
function defined on~$\mathbb{N} \times \mathbb{R}$ has another advantage.
For such functions, we can use the standard $L^p$-metric to quantify how
``different'' two landscape functions are. In fact, it was shown
in~\cite{peterLandscapes} that persistence landscapes are stable with
respect to the $L^p$-metric, i.e., small changes in the filtering function
which generates the underlying persistence diagram leads to small changes
in the resulting persistence landscape with respect to the $L^p$-norm.

Also from a computational point of view persistence landscapes are
useful. The construction of persistence landscapes has been implemented
by the first author of this paper in a library called {\tt Persistence
Landscape Toolbox}, see~\cite{plt} for more details. The library is publicly
available and can be downloaded from the first author's webpage. It implements
various manipulations of persistence landscapes, including the statistical
classifier which will be heavily used in the analysis of the present paper,
and which will be described in more detail in the next section.
\section{Topological Analysis of the Cahn-Hilliard-Cook model}
\label{sec:chc}
We now turn to the analysis of the Cahn-Hilliard-Cook model using
persistence landscapes. We begin in Section~\ref{sec:basicmethod} by
introducing the basic methodology, including a description of the
persistence landscapes based statistical classifier and of the numerical
simulation techniques. After that, Section~\ref{sec:chcmass} is concerned
with showing that the topology evolution of~(\ref{chc}) encodes the central
mass parameter~$\mu$, while Section~\ref{sec:chctime} illustrates that
even the decomposition stage can be recovered accurately from the
topological information.
\subsection{Basic Methodology}
\label{sec:basicmethod}
In order to generate the microstructures for the persistence landscape
analysis, we simulated the Cahn-Hilliard-Cook model~(\ref{chc}) at
the parameter values
\begin{displaymath}
  \epsilon = 0.005 \; , \qquad
  \sigma = 0.001 \; , \qquad
  \mbox{ with space-time white noise~$\xi$} \; ,
\end{displaymath}
and on the square domain~$\Omega = (0,1)^2 \subset \mathbb{R}^2$.
All simulations were performed in~{\tt C++} using a linearly implicit
spectral method. More precisely, the solution~$u(t,x)$ for $t \ge 0$
and $x = (x_1,x_2) \in \Omega$ is approximated using the Fourier
expansion
\begin{equation} \label{def:spectralapprox}
  u(t,x) \; \approx \;
  \sum_{k_1,k_2=0}^{K-1}
    \hat{u}_{k_1,k_2}(t) \cdot c_{k_1,k_2} \cdot
    \cos k_1 \pi x_1 \cdot \cos k_2 \pi x_2 \; ,
\end{equation}
where the constants~$c_{k_1,k_2}$ are chosen in such a way that the
functions $c_{k_1,k_2} \cos k_1 \pi x_1 \cos k_2 \pi x_2$ form a complete
orthonormal set in~$L^2(\Omega)$, for $k_1,k_2 \in \mathbb{N}_0$. Also the
white noise process~$\xi$ is approximated in this way, i.e., we use cut-off
noise which acts independently on all~$K^2 - 1$ Fourier modes in the
representation~(\ref{def:spectralapprox}). Recall that the noise is
assumed to be mass-conserving, so it does not act on the constant mode.
Throughout this paper, the simulations use $K = 256$, and to avoid aliasing
effects, the necessary Fourier coefficients of the nonlinearity~$F'(u)$ are
computed using the two-dimensional discrete cosine transform on a grid of
size~$512^2$.

For a variety of different mass values~$\mu$, we then performed the
following simulations. Let~$u_0 \in L^2(\Omega)$ denote a randomly chosen
initial condition with average~$\mu$ which satisfies
\begin{displaymath}
  \left\| u_0 - \mu \right\|_{L^\infty(\Omega)} \; = \; 10^{-4} \; .
\end{displaymath}
Then the above numerical method is used to compute the solution~$u$ of
the Cahn-Hilliard-Cook model~(\ref{chc}) with $u(0,\cdot) = u_0$ for times
\begin{displaymath}
  0 \; < \; t \; \le \; T_e \; = \;
  \frac{80 \epsilon^2}{F''(\mu)^2} \; ,
\end{displaymath}
where we recall that $F''(\mu) = 3\mu^2 - 1$. While the choice of
simulation endtime might seem strange at first sight, it takes into
account the different instability strengths of the homogeneous state~$\mu$
as a function of the mass. In fact, by choosing the endtime~$T_e$ as 
above every simulation leads to basically the same phase separation
horizon. The time interval is discretized using~$100000$ steps, and
solution snapshots are recorded at a variety of times, as will be
described in more detail below.
\begin{figure} \centering
  \setlength{\unitlength}{1 cm}
  \begin{picture}(16.0,9.9)
    \put(0.0,6.6){
      \includegraphics[width=5cm]{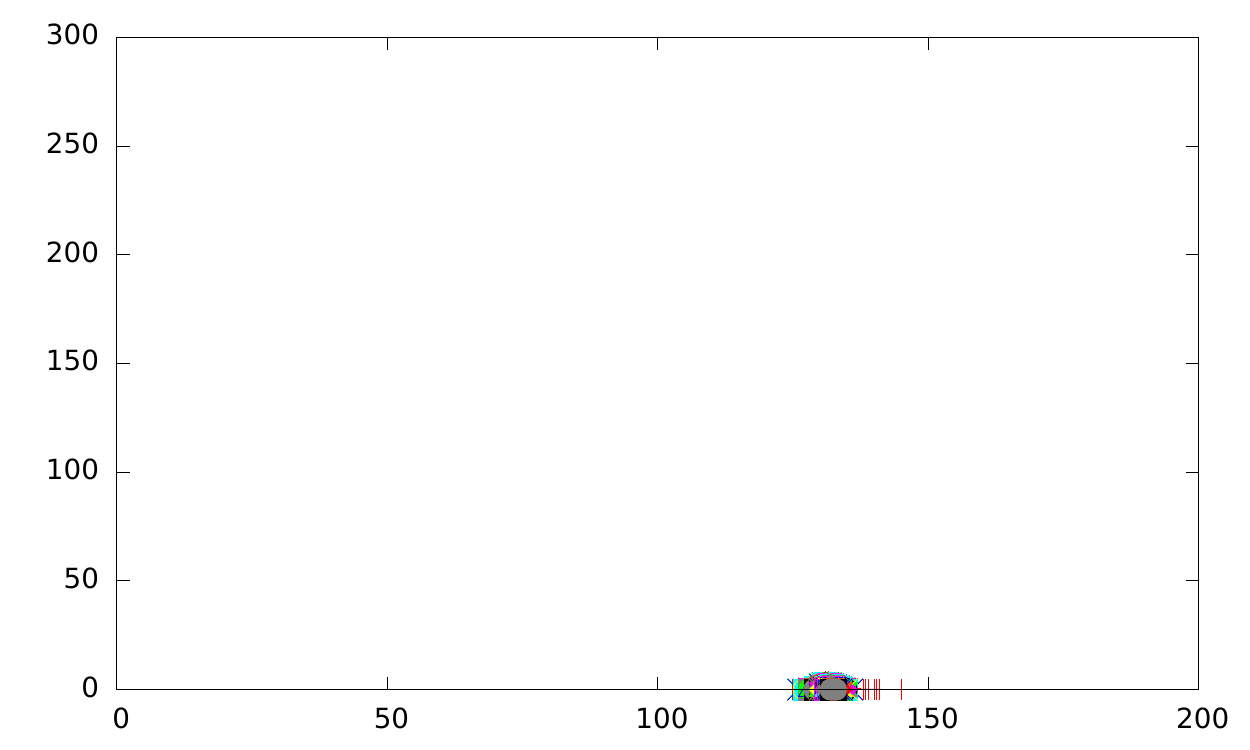}}
    \put(5.5,6.6){
      \includegraphics[width=5cm]{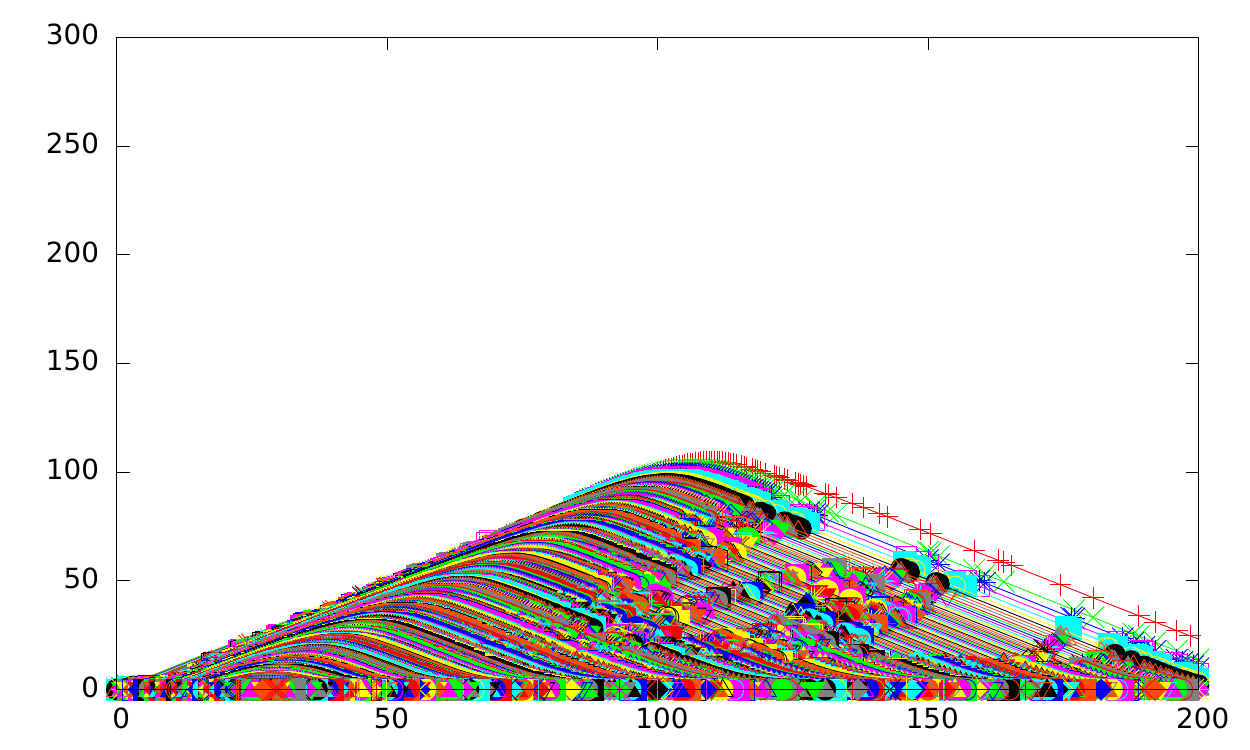}}
    \put(11.0,6.6){
      \includegraphics[width=5cm]{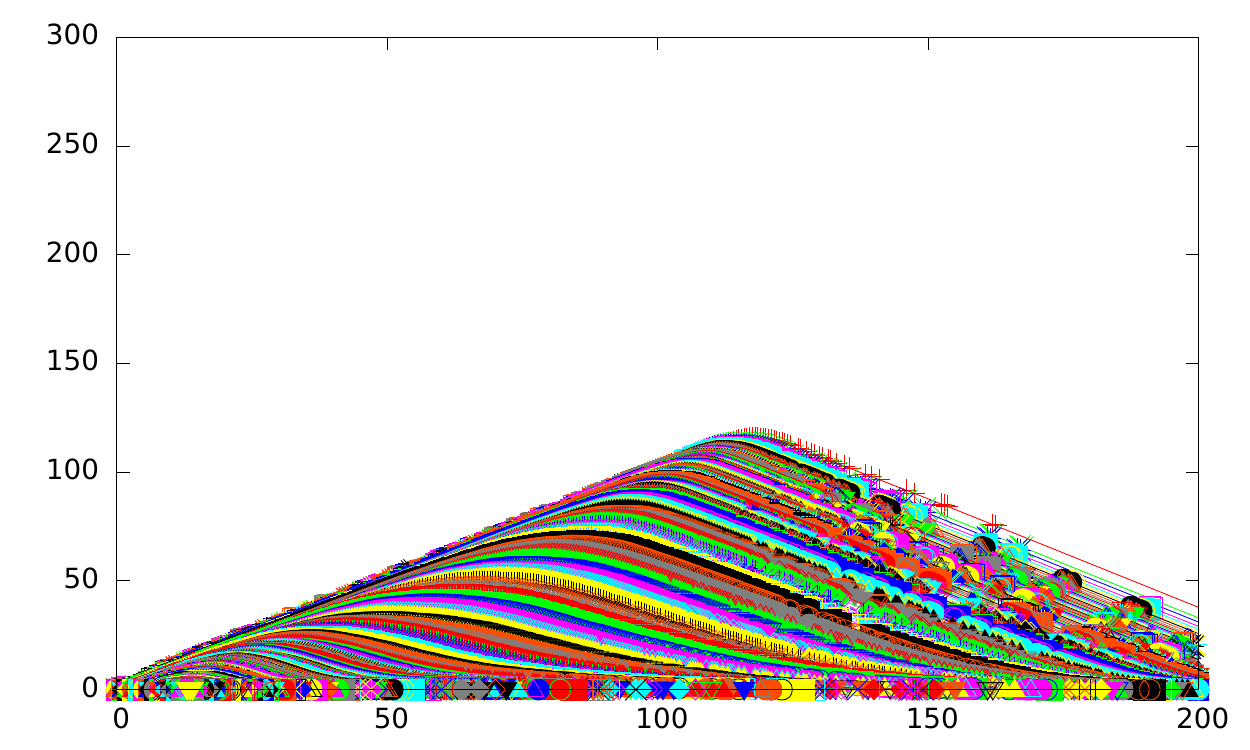}}
    \put(0.0,3.3){
      \includegraphics[width=5cm]{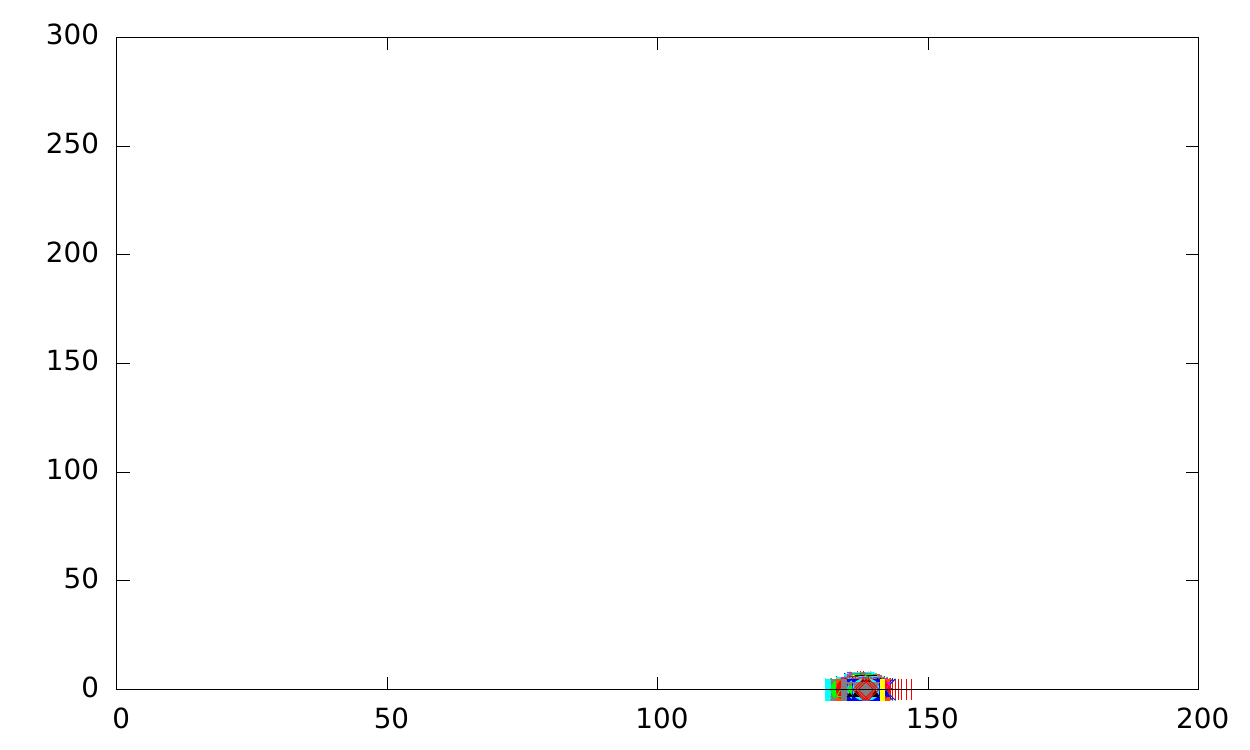}}
    \put(5.5,3.3){
      \includegraphics[width=5cm]{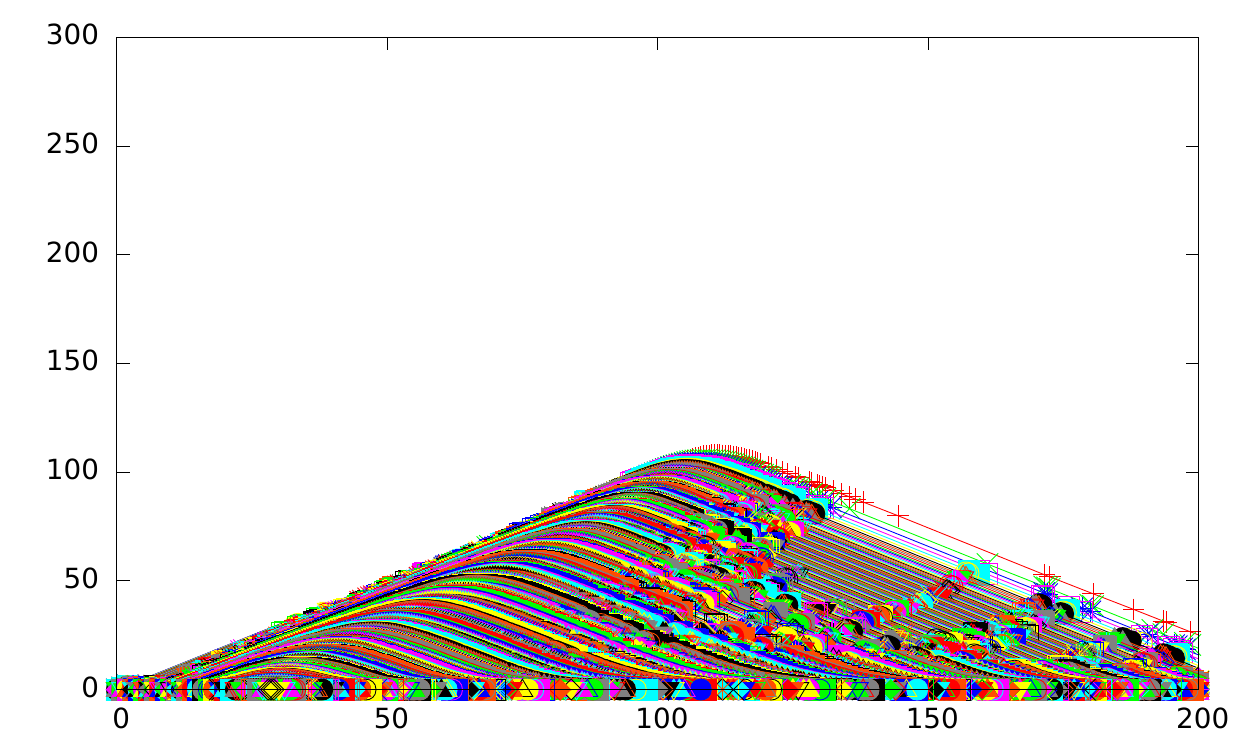}}
    \put(11.0,3.3){
      \includegraphics[width=5cm]{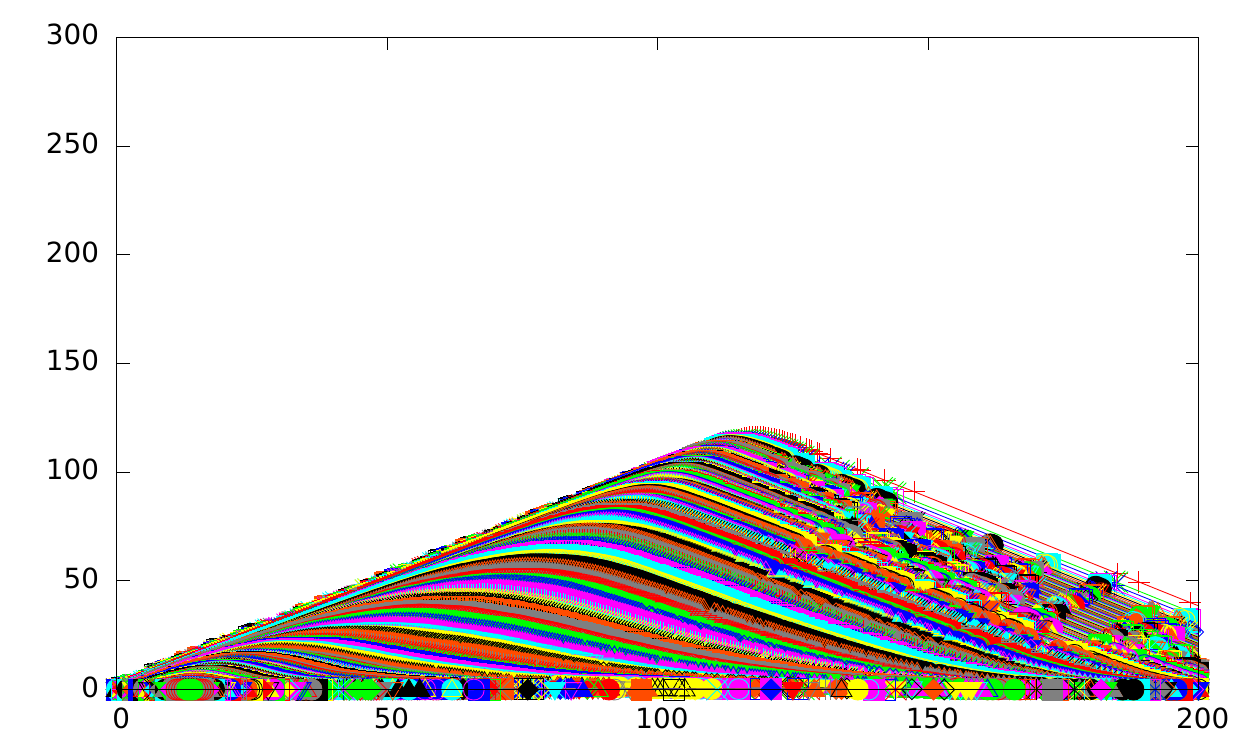}}
    \put(0.0,0.0){
      \includegraphics[width=5cm]{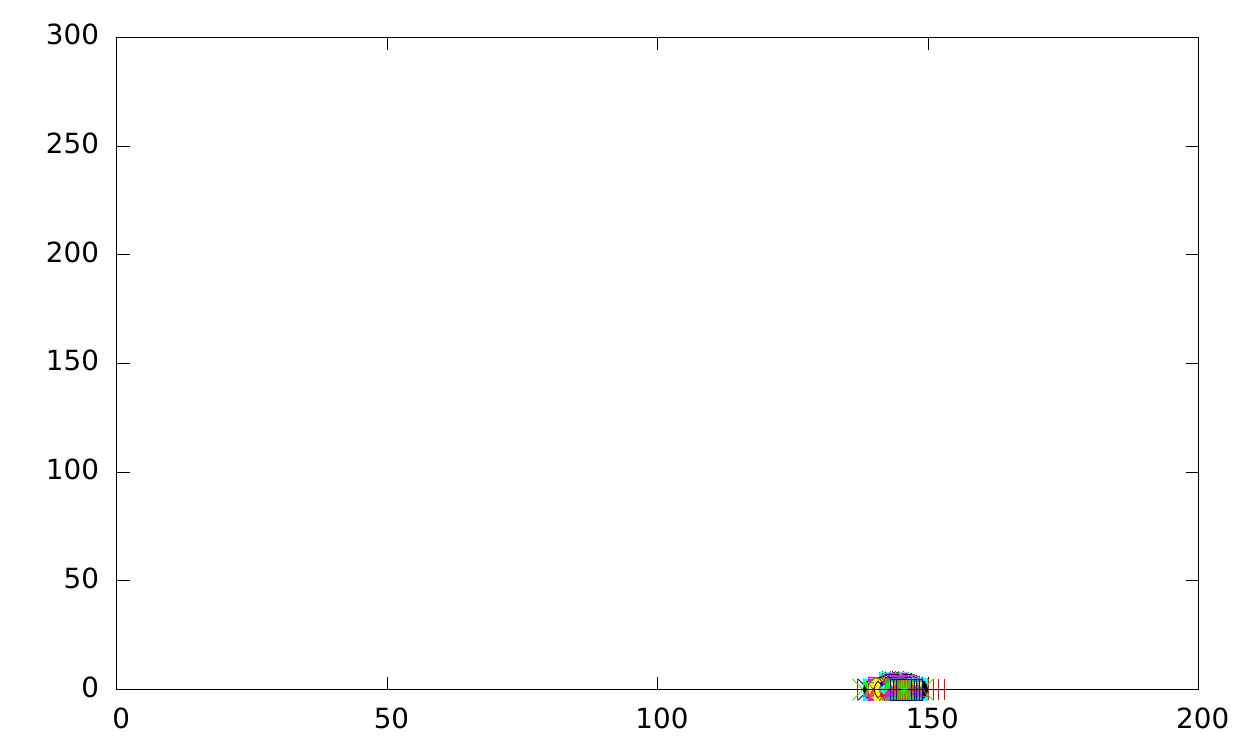}}
    \put(5.5,0.0){
      \includegraphics[width=5cm]{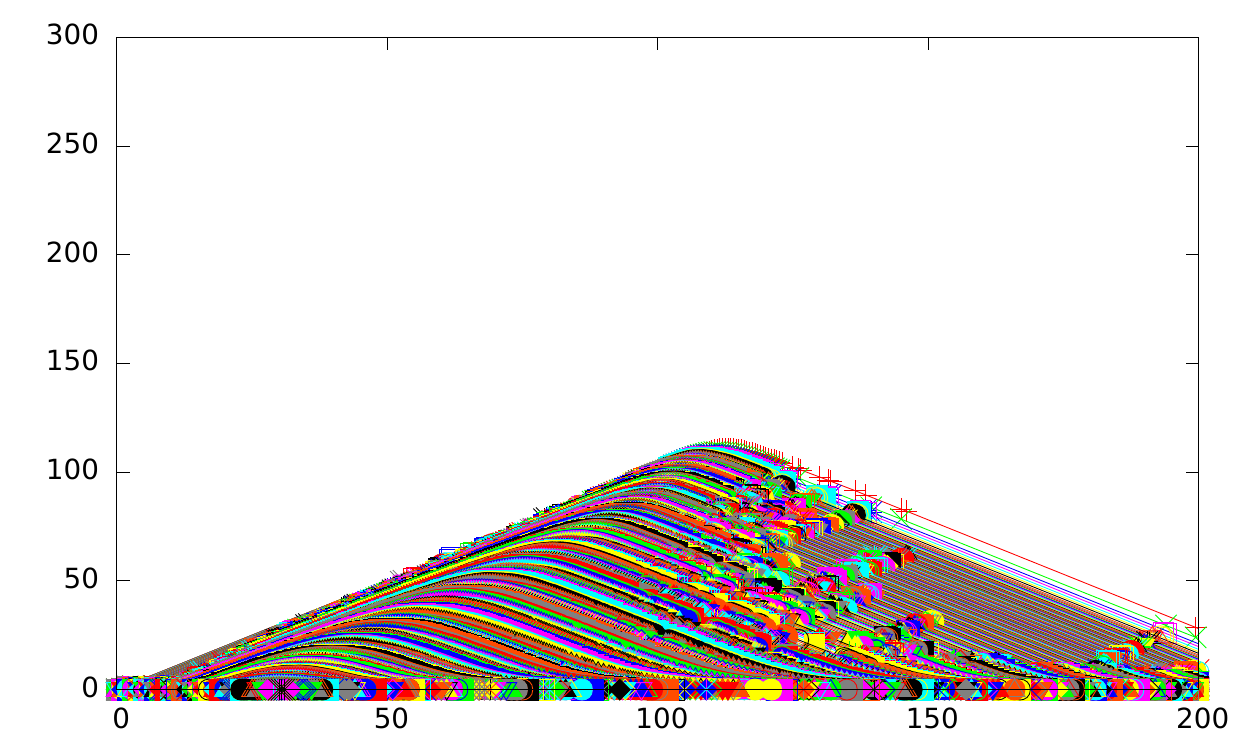}}
    \put(11.0,0.0){
      \includegraphics[width=5cm]{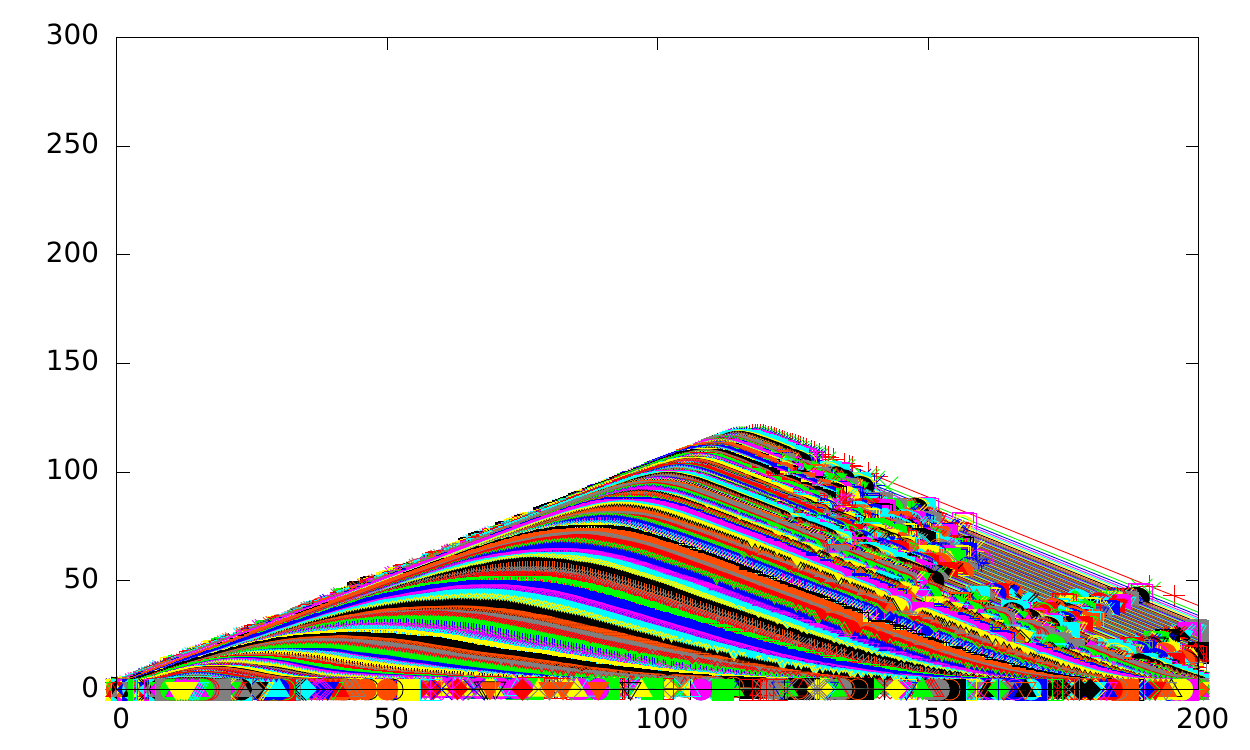}}
  \end{picture}
  \caption{Sample persistence landscapes in dimension zero for microstructures
           created by the Cahn-Hilliard-Cook model~(\ref{chc}). From top
           to bottom, the rows correspond to mass $\mu = 0$, $0.01$,
           and~$0.02$, respectively. From left to right, the respective columns
           use solution snapshots at times $t = 0.2 T_e$, $0.6 T_e$, and~$T_e$.}
  \label{fig:pldchdim0}
\end{figure}
\begin{figure} \centering
  \setlength{\unitlength}{1 cm}
  \begin{picture}(16.0,9.9)
    \put(0.0,6.6){
      \includegraphics[width=5cm]{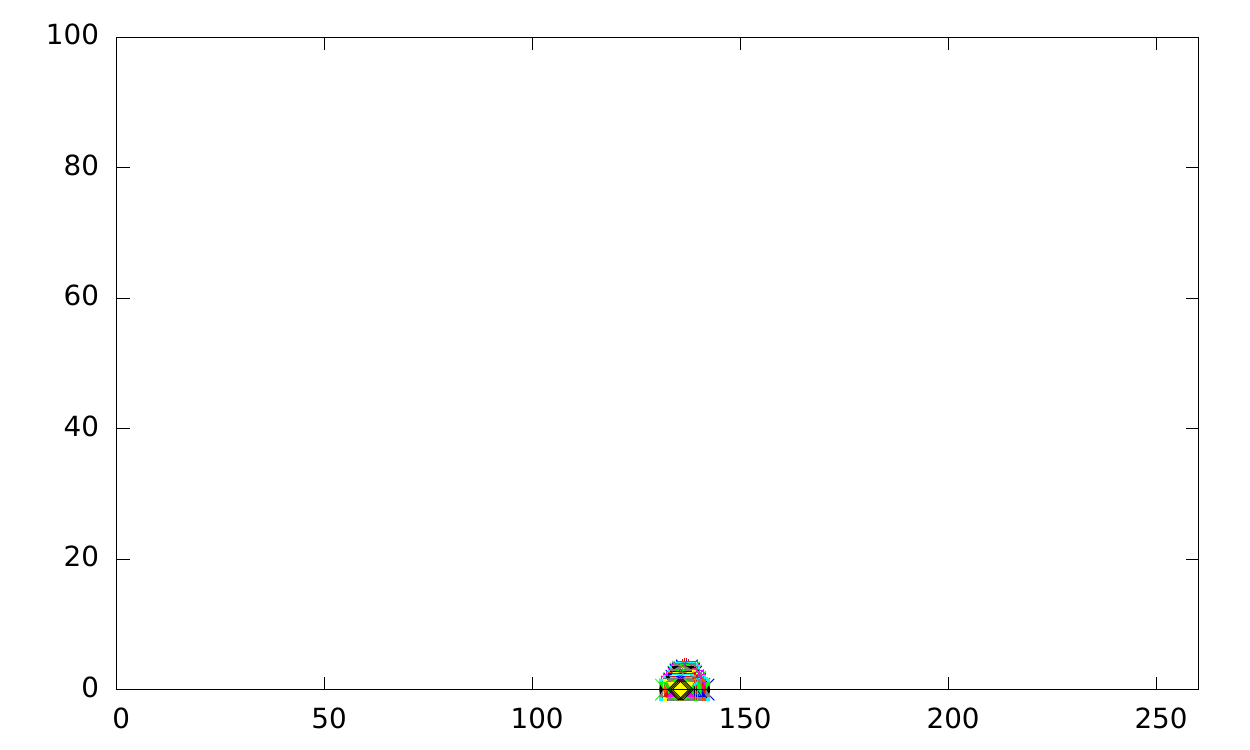}}
    \put(5.5,6.6){
      \includegraphics[width=5cm]{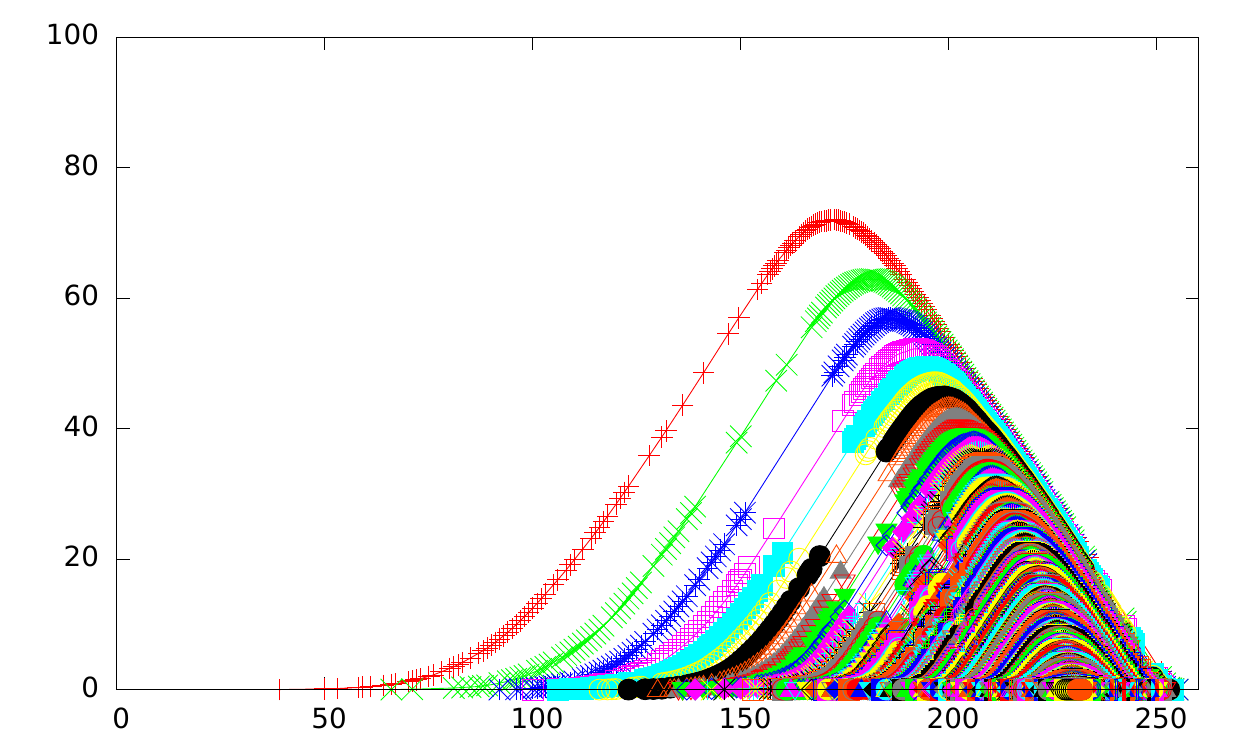}}
    \put(11.0,6.6){
      \includegraphics[width=5cm]{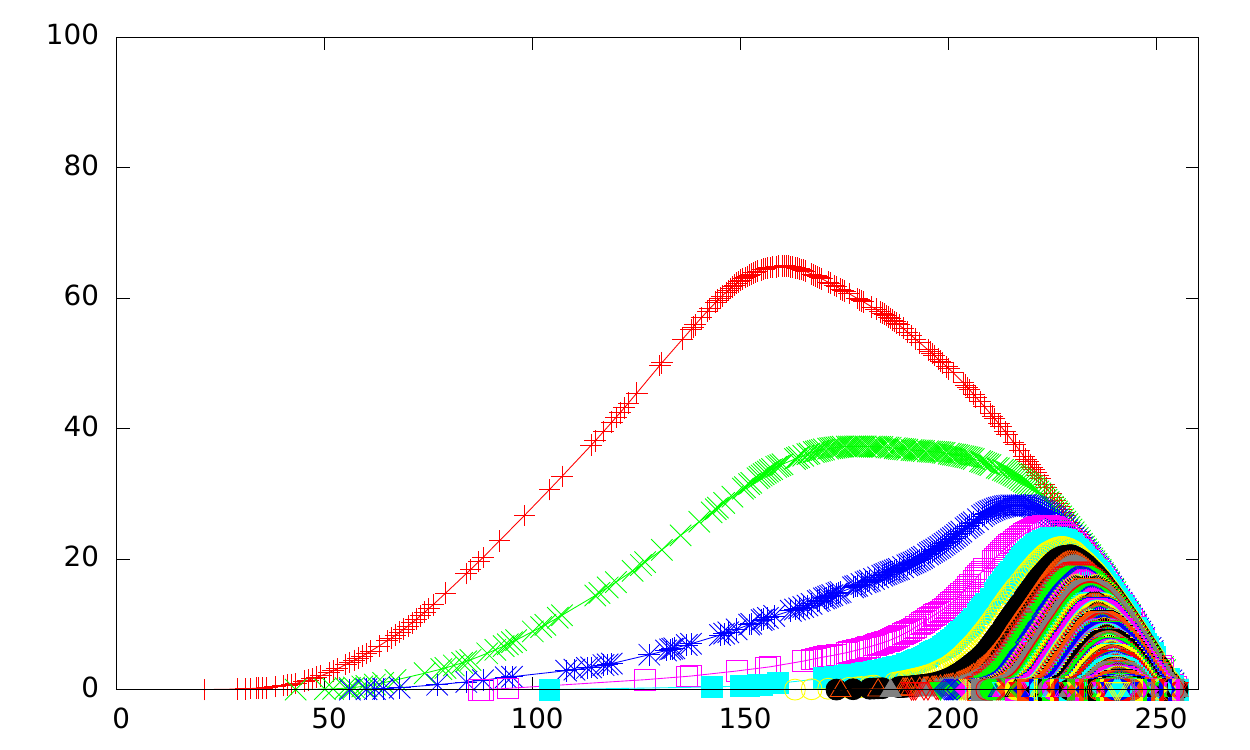}}
    \put(0.0,3.3){
      \includegraphics[width=5cm]{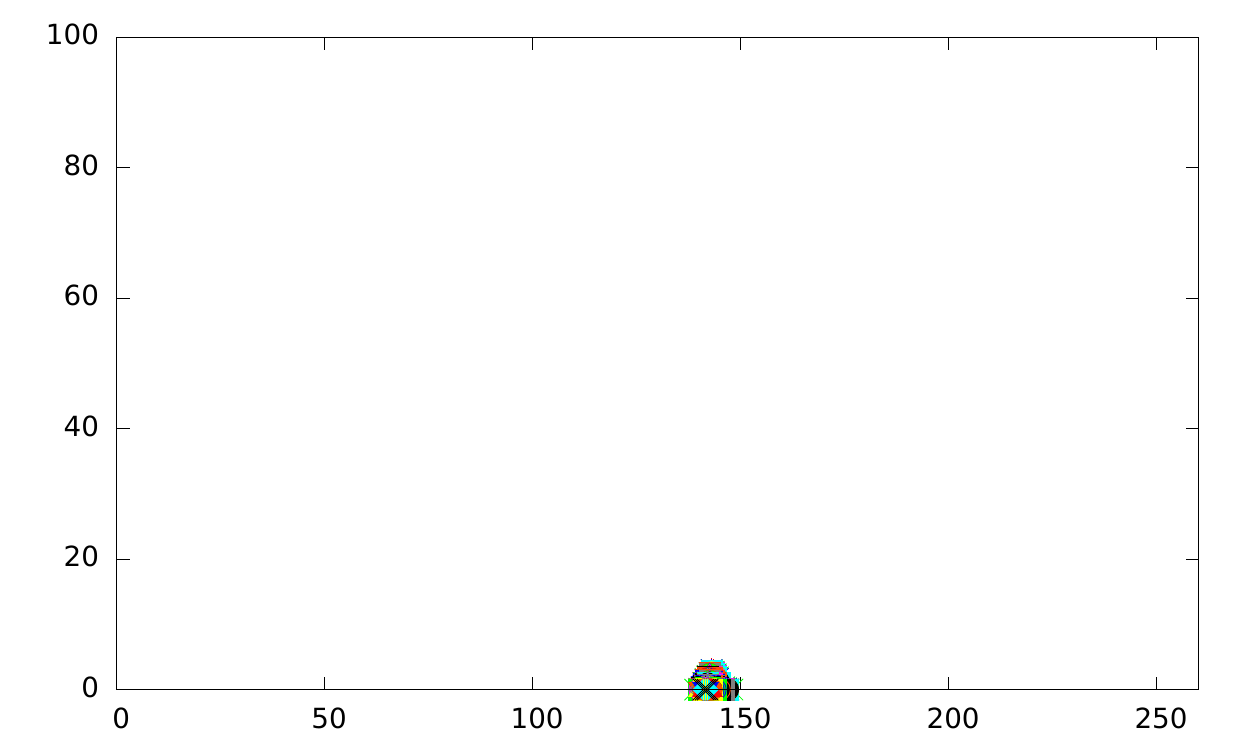}}
    \put(5.5,3.3){
      \includegraphics[width=5cm]{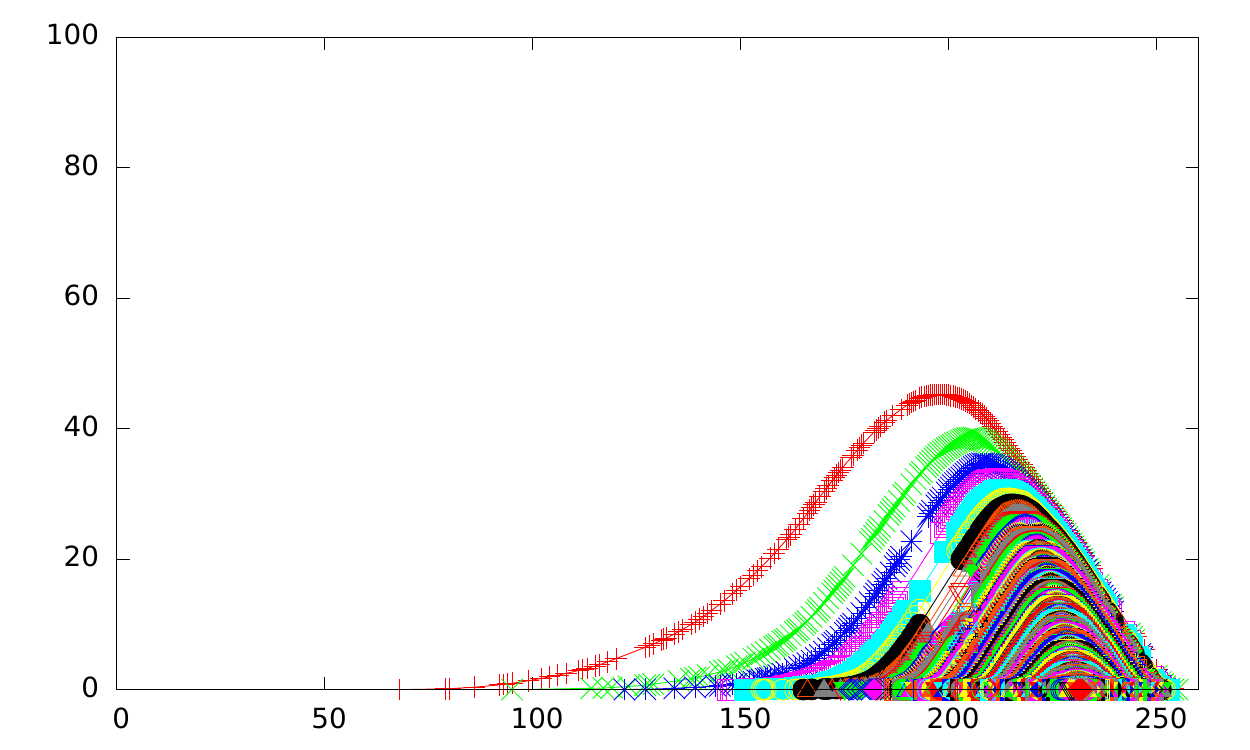}}
    \put(11.0,3.3){
      \includegraphics[width=5cm]{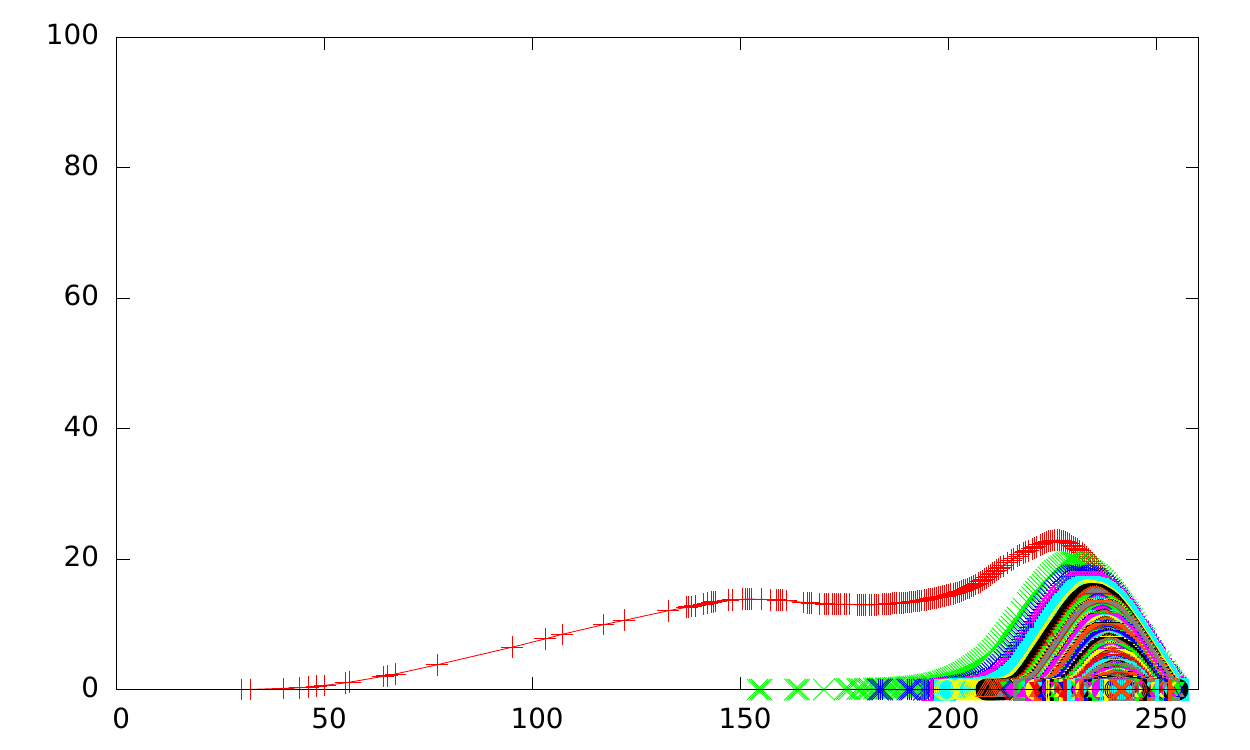}}
    \put(0.0,0.0){
      \includegraphics[width=5cm]{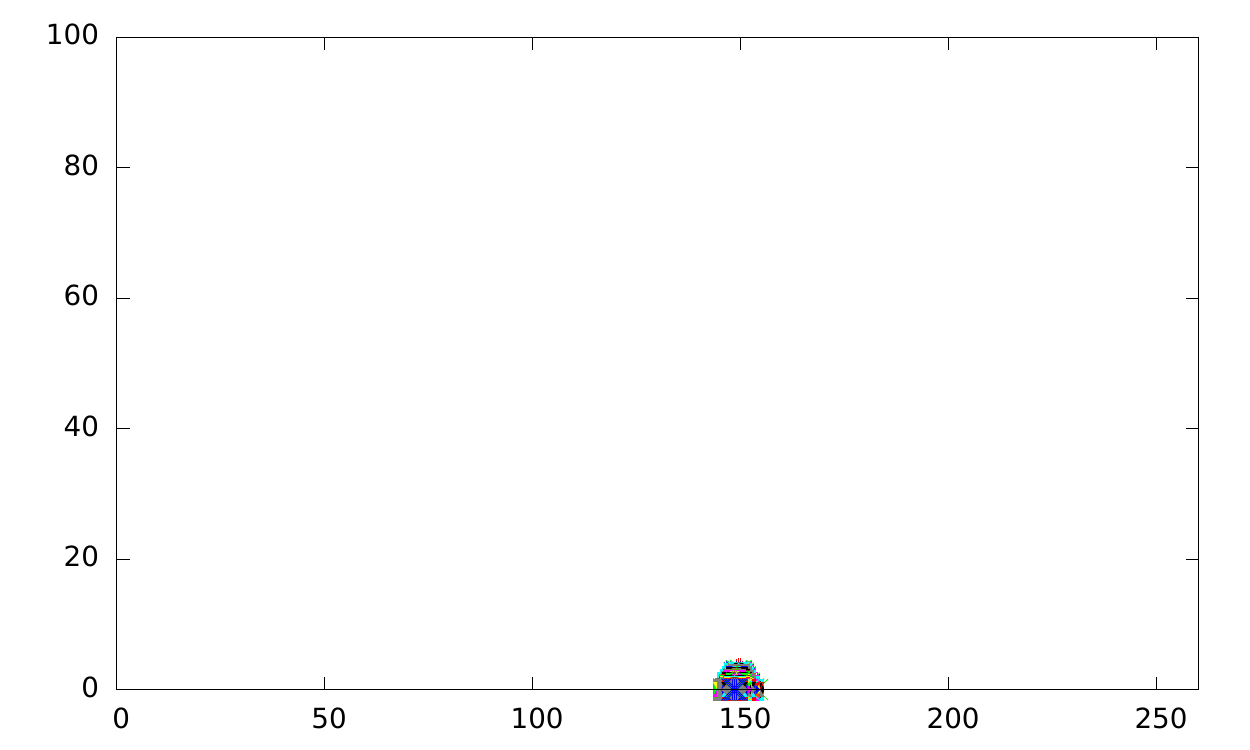}}
    \put(5.5,0.0){
      \includegraphics[width=5cm]{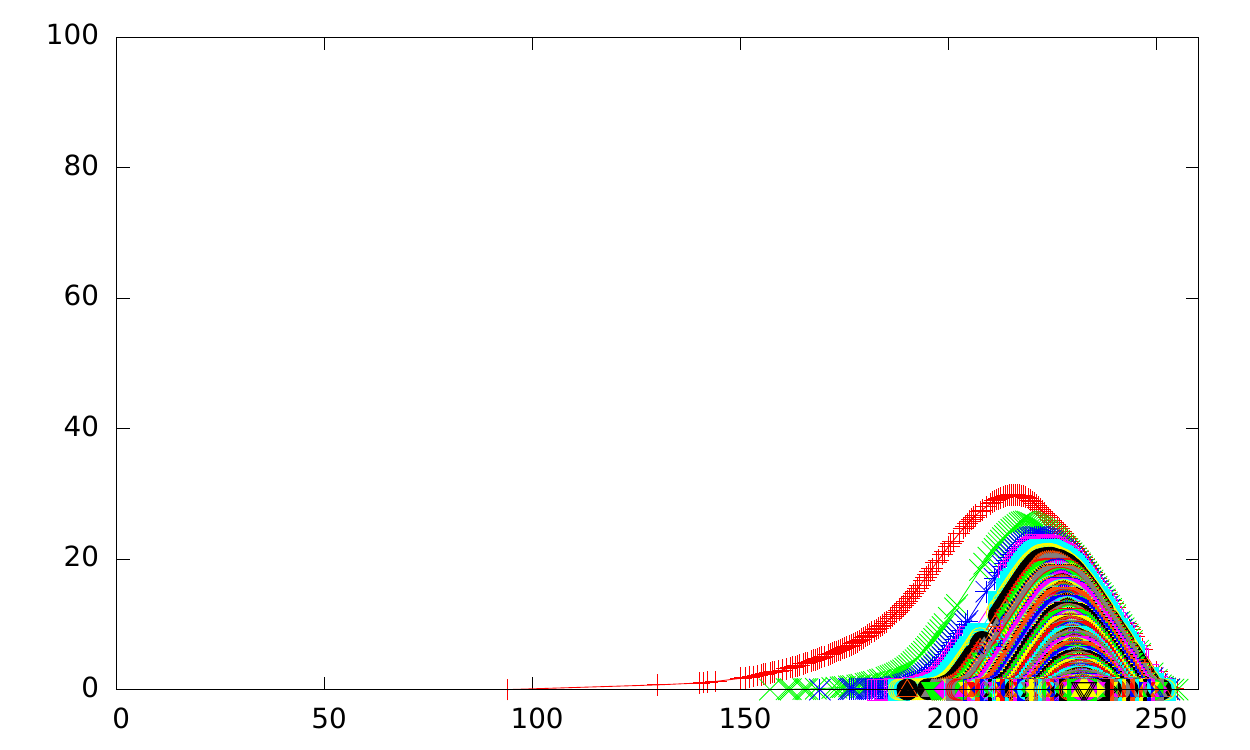}}
    \put(11.0,0.0){
      \includegraphics[width=5cm]{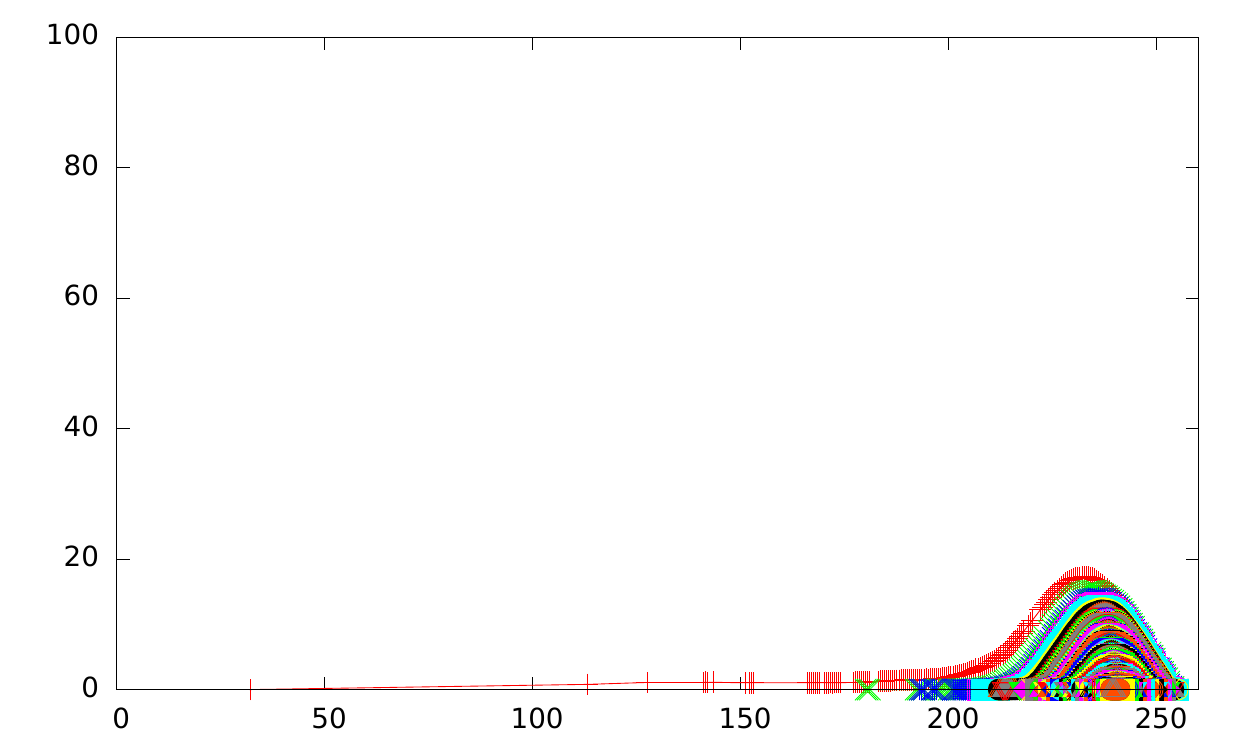}}
  \end{picture}
  \caption{Sample persistence landscapes in dimension one for microstructures
           created by the Cahn-Hilliard-Cook model~(\ref{chc}). From top
           to bottom, the rows correspond to mass $\mu = 0$, $0.01$,
           and~$0.02$, respectively. From left to right, the respective columns
           use solution snapshots at times $t = 0.2 T_e$, $0.6 T_e$, and~$T_e$.}
  \label{fig:pldchdim1}
\end{figure}

The solution snapshots~$u(t_*,\cdot)$ obtained through this procedure are
then analyzed using persistence landscapes. As a first step, the spatial 
domain is discretized into~$512^2$ equal-size subsquares, and on each square
the function~$u$ is approximated using its function value at the center 
of the subsquare. In this way, we obtain a pixelated image version of the
microstructure, as shown in Figure~\ref{fig:chcpatt}, where the function
values of~$u$ take the role of the gray-scales in the previous section. Next,
the solution range~$[-1,1]$ is discretized into~$256$ equal-sized intervals,
and the persistence diagrams are computed by considering sublevel sets
of~$u(t_*,\cdot)$ as the threshold parameter increases in~$256$ steps
from~$-1$ to~$1$. Notice that since the solution~$u$ could take values
slightly larger than one, we actually consider the function~$\min\{ u, 1\}$
instead of~$u$. In the language of Section~\ref{sec:intropers}, we consider
the cubical complex~$\mathcal{K}$ which consists of~$512^2$ equal
squares to form the unit square~$\Omega$ (where one unit in the cubical
complex representation corresponds to~$1/512$ units in the $\Omega$-scale),
and the filtering function~$\min\{ u, 1\} : \Omega \to (-\infty,1]$ evaluated
at the centers of the subsquares. Together with~(\ref{def:filteringfcomplex2})
this leads to a filtration~$\mathcal{K}_1 \subset \ldots \subset \mathcal{K}_{256}$,
where the thresholds~$a_i$ are chosen uniformly between~$-1$ and~$+1$. The
persistence diagrams are computed using the software {\tt PHAT}~\cite{phat},
together with an implementation of cubical complexes. Those computations can
also be performed with {\tt Perseus}~\cite{perseus}. The associated persistence
landscapes are then determined using~\cite{plt}. Some of the resulting persistence
landscapes can be found in Figures~\ref{fig:pldchdim0} and~\ref{fig:pldchdim1}
for dimensions zero and one, respectively. As we pointed out in the last
section, none of the patterns discussed in our setting will have
non-trivial two-dimensional homology, as the underlying two-dimensional
cubical complexes are flat.

Due to the stochastic nature of the Cahn-Hilliard-Cook model, the above
simulations are repeated many times in a Monte-Carlo type fashion. In this
way, we created three data sets that will then be analyzed further:
\begin{itemize}
\item[(D1)] The first data set considers the~$51$ mass values
$\mu = 0.01 \, n$ for $n = 0,\ldots,50$, and in each case five different
solution paths of~(\ref{chc}), i.e., five different initial
conditions~$u_0$. Along the solution, snapshots are saved at
times $t = 0.002 \, \ell \, T_e$ for $\ell = 1,\ldots,500$.
\item[(D2)] The second data set considers the~$6$ mass values
$\mu = 0.1 \, n$ for $n = 0,\ldots,5$, and in each case fifty different
solution paths of~(\ref{chc}). Along all solution paths, snapshots are
saved at times $t = 0.002 \, \ell \, T_e$ for $\ell = 1,\ldots,500$.
\item[(D3)] The third and last data set considers the~$51$ mass values
$\mu = 0.01 \, n$ for $n = 0,\ldots,50$, and in each case~$100$ different
solution paths of~(\ref{chc}). Along each solution, snapshots are saved at
times $t = 0.04 \, \ell \, T_e$ for $\ell = 1,\ldots,25$.
\end{itemize}
Thus, in each of the three data sets we consider~$N$ different mass
values~$\mu$, where $N = 51, 6, 51$, respectively, and for each simulation
we obtain solution snapshots at~$M$ different times, where~$M = 500, 500, 25$,
respectively. Every one of these solution snapshots leads to two associated
persistence landscapes, one each in dimensions zero and one.

We would like to briefly comment on the fact that all of our data sets
consider non-negative mass values~$\mu \ge 0$. Based on our thresholding
process described above, we study sublevel sets of the function~$u$ as
the threshold increases from~$-1$ to~$+1$. In Figure~\ref{fig:chcpatt},
these sublevel sets for the threshold~$\theta = 0$ are shown in red.
Thus, for mass values~$0 \ll \mu \le 0.5$ the sublevel sets always have 
nontrivial $1$-dimensional homology, which in fact measures the bulk
or interior behavior of the material, as described in~\cite{desi:etal:11a,
gameiro:etal:05a}. In contrast, one can easily see that for mass values
$-0.5 \le \mu \ll 0$ the sublevel sets usually have trivial first homology,
i.e., the first homology dimension is useless for topological classification.

Based on the above topological data, our goal in the next two sections
is a statistical study to solve inverse problems for the total mass~$\mu$
and the decomposition stage. While this could be done in a number of ways,
we concentrate on one statistical classification method which has been 
implemented in~\cite{plt} and which we recall now. Suppose we are given~$N$
different classes~$C_1, \ldots, C_N$ of persistence landscapes. Assume further
that each class~$C_n$ consists of~$S$ different landscapes~$L_{n,j}$, for
$j = 1,\ldots,S$. More precisely, in every class~$C_n$ we have~$S$ different
persistence landscapes~$L_{n,j}^k$, $j = 1,\ldots,S$, in dimension~$k$,
where~$k \in \mathbb{N}_0$. Our classification scheme will be based on the 
previously mentioned $L^p$-distances between persistence landscapes, and
therefore we choose and fix a real number~$p \in [1,\infty]$. Suppose we
are given a new persistence landscape~$L$ and would like to decide which of
the classes~$C_1, \ldots, C_N$ it belongs to, based only on the persistence
information in dimension $k \in \mathbb{N}_0$. Then we use the following
classification scheme:
\begin{itemize}
\item[(Ck)] Suppose~$N$ classes of persistence landscapes are given as
above. Then for each $n = 1,\ldots,N$ we define the \emph{average
classifier}~$\overline{C}_n^k$ of the~$n^{\rm th}$ class in
dimension~$k$ via
\begin{equation} \label{aveclassifierK}
  \overline{C}_n^k \; = \;
  \frac{1}{S} \cdot \sum_{j=1}^S L_{n,j}^k \; .
\end{equation}
Now let~$L$ be any other persistence landscape which we would like to
classify. Then we say that \emph{$L$ has been classified to belong to
class~$C_n$ using dimension~$k$}, if there exists a unique index
$n \in \{ 1, \ldots, N \}$ such that
\begin{equation} \label{classifierK}
  \left\| L^k - \overline{C}_n^k
    \right\|_{L^p}
  \; \le \;
  \left\| L^k - \overline{C}_j^k
    \right\|_{L^p}
  \qquad\mbox{ for all }\qquad
  j = 1, \ldots, N \; .
\end{equation}
If no such unique index~$n$ exists, then we say that
\emph{the classifier fails to classify~$L$}.
\end{itemize}
Intuitively, the classification scheme~(Ck) is easy to understand.
Based on the different persistence landscapes available in a given
class, we compute their average to determine the ``typical'' persistence
behavior in this class. Given a new persistence landscape~$L$, one then
uses the average classifier which lies closest to~$L$ to determine the
class which most likely contains~$L$.

The above classification scheme~(Ck) in fact describes a whole
family of classification schemes, indexed by the homology dimension~$k$.
In our situation, only the cases~$k = 0$ and~$k = 1$ lead to nontrivial
persistence landscapes, and therefore we will use only schemes~(C0) and~(C1)
in the following. Note, however, that in some situations one might like to
use all available topological information, i.e., all dimensions~$k$, at the
same time to determine the most likely class which contains a given persistence
landscape~$L$. This can be achieved by the following simple scheme.
\begin{itemize}
\item[(CA)] Suppose~$N$ classes of persistence landscapes are given as
above, and that the average classifier~$\overline{C}_n^k$ of the~$n^{\rm th}$
class in dimension~$k$ has been defined as in~(\ref{aveclassifierK}),
for $n = 1,\ldots,N$ and $k \in \mathbb{N}_0$. As before, let~$L$ be any
other persistence landscape which we would like to classify.

For each dimension~$k \in \mathbb{N}_0$, one can then order the~$N$
possible classification classes according to monotone increasing values of the
distances~$\| L^k - \overline{C}_n^k \|_{L^p}$, and we let~$\mathcal{S}_j^k$
denote the first~$j$ class indices in this sequence, for $j = 1,\ldots,N$.
Note that $\mathcal{S}_1^k$ is a singleton set which contains a class
whose averaged classifier~$\overline{C}_n^k$ is closest to~$L^k$, while
we have $\mathcal{S}_N^k = \{ 1,\ldots,N \}$ for all $k \in \mathbb{N}_0$.
Furthermore, let~$\ell = 1,\ldots,N$ be the smallest number such that
$\cap_{k=0}^\infty \mathcal{S}_\ell^k \neq \emptyset$.

Then we say that \emph{$L$ has been classified to belong to
class~$C_n$ using all dimensions}, if any ordering of the classes
as above leads to
\begin{equation} \label{classifierKA}
  \bigcap\limits_{k=0}^\infty \mathcal{S}_\ell^k
  \; = \; 
  \left\{ n \right \} \; .
\end{equation}
If this intersection contains more than one element, or if
there are ties in the above orderings due to equidistant average
classifiers, then we say that \emph{the classifier fails to
classify~$L$}.
\end{itemize}
This classification scheme can be useful if no clear choice of
dimension~$k$ in~(Ck) is evident a-priori, and we will use it also
in our further studies.

For our application to the Cahn-Hilliard-Cook equation in Section~\ref{sec:chcmass}
we are using a slight extension of the above classification schemes. For each
simulation, our input data consists of a sequence of~$M$ different persistence
landscapes which correspond to microstructure patterns at times~$t_1, \ldots, t_M$.
Such a sequence of persistence landscapes will be called a \emph{topological
process}. If we define the average of topological processes $P_1, \ldots, P_S$
as the topological process whose $i^{\rm th}$ persistence landscape in
dimension~$k$ is the average of all the $i^{\rm th}$ persistence landscapes
in dimension~$k$ of the processes $P_1, \ldots, P_S$, then we can easily talk
about an averaged classifier as in~(Ck). Furthermore, if we define the
$L^p$-distance between topological processes~$P$ and~$Q$ as the sum of
the $L^p$-distances between the persistence landscapes in the~$i^{\rm th}$
positions of the processes~$P$ and~$Q$, for $i = 1,\ldots, M$, then one can 
easily reformulate~(\ref{classifierK}) for the case of topological processes.
The total classification scheme~(CA) can be extended analogously. The
resulting schemes will be used for the classification of the Cahn-Hilliard-Cook
simulations described in~(D1), (D2), and~(D3).
\subsection{Identifying the Mass via Topology}
\label{sec:chcmass}
As described in the last section, our simulations of the Cahn-Hilliard-Cook
model~(\ref{chc}) provide snapshots of a solution realization along~$M$
equidistant timesteps~$t_1, \ldots, t_M$. For each time~$t_\ell$, computational
topology is then used to create the associated persistence landscapes in dimensions~$0$
and~$1$. In other words, for each simulation of~(\ref{chc}) we obtain a topological
process which captures the persistence evolution of the created microstructure
over the time interval~$(0,T_e]$. In this sense, our topological data is a 
refinement of the topology evolution curves which were considered
in~\cite{gameiro:etal:05a}.

What does this refined topology evolution capture? In the present section,
we demonstrate that it is in fact possible to determine the total mass~$\mu$
of the underlying simulation extremely accurately. This implies the somewhat
surprising observation that the total mass~$\mu$ has quite a significant 
influence on the pattern formation process --- to the extent that persistence
information alone, which of course only depends on $k$-dimensional connectivity
measurements, but not on size information, suffices to accurately infer~$\mu$.

We demonstrate this mass-topology dependence using the three data sets~(D1),
(D2), and~(D3) described in Section~\ref{sec:basicmethod}. The simulations
in each of these data sets give rise to topological processes of length~$M$,
and they can be divided into~$N$ different mass classes. While in~(D1) and~(D3)
we considered the $N = 51$ mass values $\mu = 0.01 \, n$ for $n=0,\ldots,50$, the
data set~(D2) considers the $N = 6$ mass values $\mu = 0.1 \, n$, where now
we have $n=0,\ldots,5$. Furthermore, the length of the topological processes
is given by $M = 500$ for~(D1) and~(D2), and $M = 25$ for data set~(D3).
For each mass value~$\mu$, a data set contains~$R$ stochastically independent
repetitions of the simulation. For~(D1), (D2), and~(D3) we have $R = 5$,
$R = 50$, and $R = 100$, respectively.

Based on the composition of the data sets, we apply the classification
schemes~(CA), as well as~(C0) and~(C1), from the previous section for
mass value classes~$C_1, \ldots, C_N$, which correspond to the chosen~$N$
mass values~$\mu$ in increasing order. Then we proceed as follows.
\begin{itemize}
\item Based on~$R_T < R$ training runs each from the~$N$ considered mass
values, we determine the averaged classifier~$\overline{C}_n^k$ of
the~$n^{\rm th}$ class in dimension~$k$ as defined in~(\ref{aveclassifierK})
with $S = R_T$, but extended for the case of topological processes.
\item We then try to classify each of the remaining~$(R - R_T) N$ simulations
in the data set, by either using the classification scheme~(CA), or one of
the schemes~(C0) or~(C1).
\end{itemize}
In each case, we note the index~$n$ obtained from the classification scheme,
and compare it to the actual mass index~$n_{true}$ of the underlying simulation.
If~$n = n_{true}$, then we have obtained a ``Hit'', otherwise we say that the
classification ``missed by~$|n - n_{true}|$''. Note that the integer
value~$|n - n_{true}|$ is directly proportional to the mass difference
incurred by the classification. For the data sets~(D1) and~(D3) the
associated proportionality factor is~$0.01$, while for~(D2) it is~$0.1$.
\begin{table}[tb]
  \centering
  \begin{tabular}{||c|c|c|c|c|c|c||}
  \hline\hline
  \multicolumn{7}{||c||}{{\bf $L^1$-norm}} \\
  \hline
  Batch No.\ & Hits & Missed by $1$ & Missed by $2$ & Missed by $3$ & Missed by $4$
    & Wrong \\
  \hline
  $1$ & $123$ & $29$ & $0$ & $0$ & $1$ & $0$ \\
  \hline
  $2$ & $128$ & $23$ & $1$ & $0$ & $1$ & $0$ \\
  \hline\hline
  \multicolumn{7}{||c||}{{\bf $L^2$-norm}} \\
  \hline
  Batch No.\ & Hits & Missed by $1$ & Missed by $2$ & Missed by $3$ & Missed by $4$
    & Wrong \\
  \hline
  $1$ & $112$ & $39$ & $1$ & $0$ & $1$ & $0$ \\
  \hline
  $2$ & $110$ & $41$ & $1$ & $0$ & $1$ & $0$ \\
  \hline\hline
  \multicolumn{7}{||c||}{{\bf $L^\infty$-norm}} \\
  \hline
  Batch No.\ & Hits & Missed by $1$ & Missed by $2$ & Missed by $3$ & Missed by $4$
    & Wrong  \\
  \hline
  $1$ & $78$ & $58$ & $14$ & $1$ & $1$ & $1$ \\
  \hline
  $2$ & $82$ & $58$ & $12$ & $0$ & $0$ & $1$ \\
  \hline\hline
  \end{tabular}
  \caption{Classification results for the data set~(D1) using the classification
           scheme~(CA). The two training set batches consist of $R_T = 2$
           simulations each. Once the averaged classifier processes~$\overline{C}_n^k$
           for $n = 1,\ldots,N$, where $N = 51$, have been computed, the
           remaining~$(R-R_T)N = 153$ solution snapshot processes are classified
           using the $L^1$-, $L^2$-, and~$L^\infty$-norms. The table contains
           the number of classification hits as well as detailed information
           on the index difference observed for the misclassifications. The
           column labeled ``Wrong'' contains classifications which missed by
           more than~$4$ indices.}
  \label{tableD1CA}
\end{table}

We now turn to the first data set~(D1), which consists of $R = 5$ simulations
each for $N = 51$ different mass values. In this case, we took the first
$R_T = 2$ data sets for each mass value as training sets, and their average
constitutes the averaged classifier in each case. Classifying the remaining
$(R-R_T)N = 153$ topological process using scheme~(CA) as described above then
leads to the results shown in Table~\ref{tableD1CA}, in the rows denoted
``Batch 1''. If instead we use the third and fourth simulations as two-element
training set, one obtains the lines labeled as ``Batch 2''. For the
classifications, we used the $L^p$-norm with $p = 1, 2, \infty$.

The results in the table are remarkable. If one uses the $L^1$-norm for
the classification, then in over~80\% of the cases the classification scheme
finds the correct mass value class. If one allows for an error of the form
$\mu \pm 0.01$, then the classification succeeds in at least 98\% of the
cases. While the behavior for the $L^2$-norm is only slightly worse, the
$L^\infty$-norm does not seem to be quite as accurate an estimator. The latter
norm gives correct classifications only in at least 50\% of the cases, errors
of the form $\mu \pm 0.01$ in at least~88\%, as well as errors $\mu \pm 0.02$
in at least 98\% of the cases. In addition, the two classifications appearing
in the ``Wrong'' column have index differences of~$21$, i.e., they are significant
misclassifications. Despite the worse performance of the $L^\infty$-norm estimator,
these results clearly indicate how powerful the above topological classification
scheme is. We would like to point out that one cannot expect the classifier to
have~100\% accuracy. On the one hand, these classifications used training sets
which consist of \emph{two} topological processes, which clearly is an extremely
small batch size. Moreover, we are classifying pattern evolutions created by a
stochastic partial differential equation which originate at random initial
conditions. Given these uncertainties, the presented topological classification
method is highly effective in distinguishing the pattern evolutions from the
Cahn-Hilliard-Cook equation~(\ref{chc}).
\begin{table}[tb]
  \centering
  \begin{tabular}{||c|c|c|c|c|c|c|c||}
  \hline\hline
  \multicolumn{8}{||c||}{{\bf Dimension $k = 0$}} \\
  \hline
  Norm & Hits & Missed by $1$ & Missed by $2$ & Missed by $3$ & Missed by $4$
    & Missed by $5$ & Wrong \\
  \hline
  $L^1$ & $112$ & $40$ & $0$ & $0$ & $0$ & $0$ & $1$ \\
  \hline
  $L^2$ & $100$ & $48$ & $4$ & $0$ & $0$ & $0$ & $1$ \\
  \hline
  $L^\infty$ & $81$ & $53$ & $11$ & $3$ & $1$ & $0$ & $4$ \\
  \hline\hline
  \multicolumn{8}{||c||}{{\bf Dimension $k = 1$}} \\
  \hline
  Norm & Hits & Missed by $1$ & Missed by $2$ & Missed by $3$ & Missed by $4$
    & Missed by $5$ & Wrong \\
  \hline
  $L^1$ & $129$ & $22$ & $1$ & $0$ & $1$ & $0$ & $0$ \\
  \hline
  $L^2$ & $120$ & $29$ & $3$ & $0$ & $1$ & $0$ & $0$ \\
  \hline
  $L^\infty$ & $68$ & $59$ & $18$ & $4$ & $3$ & $1$ & $0$ \\
  \hline\hline
  \end{tabular}
  \caption{Classification results for the data set~(D1) using the classification
           schemes~(C0) and~(C1). The training set consists of two simulations,
           the remaining three simulations for each of the~$N = 51$ mass values
           are then classified using the $L^1$-, $L^2$-, and~$L^\infty$-norms.
           The table contains the number of classification hits as well as
           detailed information on the index difference observed for the
           misclassifications. The column labeled ``Wrong'' contains
           classifications which missed by more than~$5$ indices.}
  \label{tableD1C01}
\end{table}

For the above classifications we have used the classification scheme~(CA)
which simultaneously considers persistence information in dimensions~$k = 0$
and~$k = 1$. How do the dimension-dependent classification schemes~(C0)
and~(C1) fare? The results for ``Batch 1'' as above are collected in
Table~\ref{tableD1C01}. As before, both the $L^1$- and the $L^2$-norm
classifications lead to precise mass value classifications, while the
classification based on the $L^\infty$-norm is worse. Note, however, that
the performance of the~(C1) scheme is consistently better than the one
of the~(C0) scheme. We believe that this can be explained as follows.
In our simulations, we considered mass values~$\mu \ge 0$. This implies
that in the droplet forming regime~$0 \ll \mu \le 0.5$, most threshold
values lead to a connected cubical complex with a number of holes, which
correspond to the droplets of one material forming in the background
matrix consisting of the second material. In this situation, the first
Betti number only counts the \emph{interior droplets}, i.e., the ones
that do not touch the boundary. In other words, in our setup, the 
$1$-dimensional Betti number of the considered sublevel set is equal
to the number of interior components of the minority phase. In contrast,
the $0$-dimensional Betti number counts both the interior \emph{and}
the boundary components of the majority phase. Thus, it seems reasonable
to attribute the better performance of the scheme~(C1) to its exclusive
focus on the \emph{bulk behavior} of the material, while~(C0) measures the
combination of both bulk and boundary phenomena. As was pointed out
in~\cite{desi:etal:11a}, these phenomena exhibit different scaling 
behaviors, and can therefore affect the results for certain values of
the parameter $\epsilon > 0$.
\begin{figure} \centering
  \setlength{\unitlength}{1 cm}
  \begin{picture}(16.75,10.5)
    \put(0.0,5.5){
      \includegraphics[height=5cm]{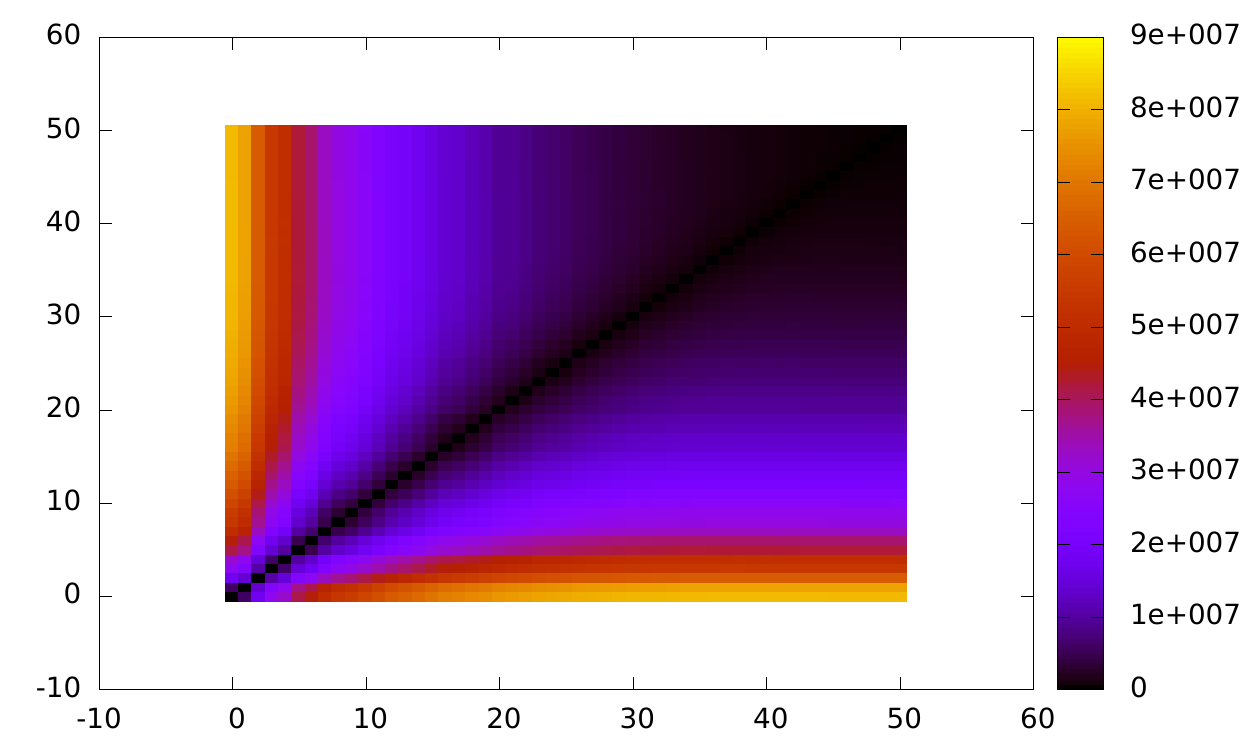}}
    \put(8.75,5.5){
      \includegraphics[height=5cm]{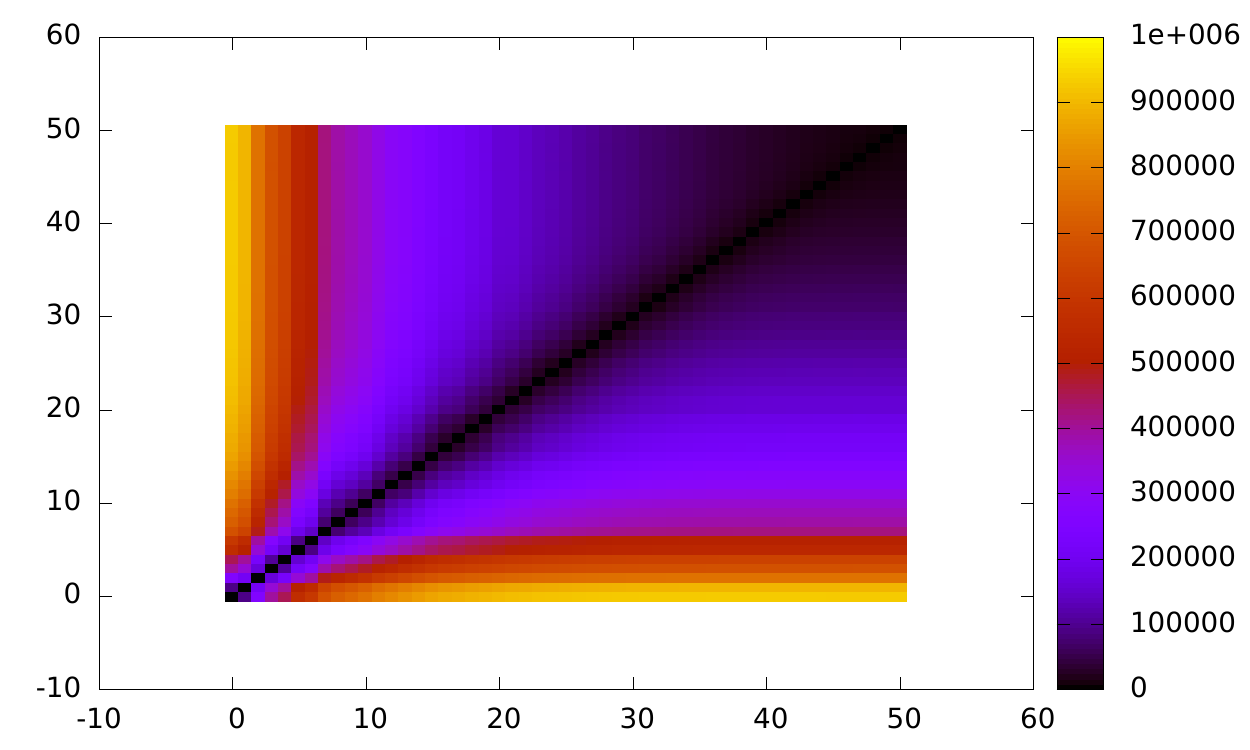}}
    \put(4.375,0.0){
      \includegraphics[height=5cm]{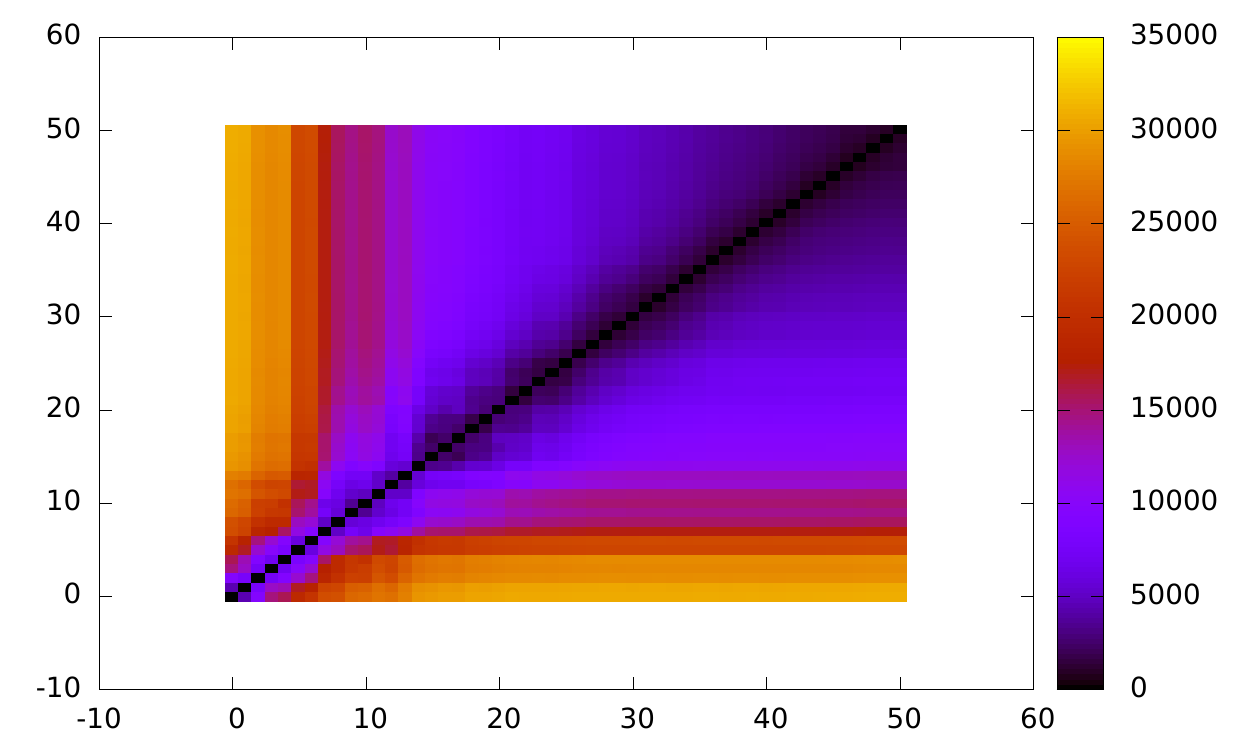}}
  \end{picture}
  \caption{Distances between the averaged classifiers~$\overline{C}_n^k$ of
           the~$n^{\rm th}$ mass value class for the data set~(D1), where
           $n = 1,\ldots,51$. The horizontal and vertical axes correspond
           to the mass value index~$n$, and the sum of the distances for
           dimensions $k = 0$ and $k = 1$, as defined in~(\ref{def:heatmapdist}),
           are color-coded as shown in the colorbars. From top left to bottom middle
           the images are for the $L^1$-, $L^2$-, and $L^\infty$-norm, respectively.}
  \label{fig:DistancesHeatMap}
\end{figure}

Before we move on to the other two data sets, we would like to elaborate
a bit more on the relative behavior of the $L^p$-classifiers for $p = 1,2,\infty$
in the above setting. One would expect good classification results only if
the averaged classifiers~$\overline{C}_n^k$ defined in~(\ref{aveclassifierK}),
but extended to the case of topological processes,
are suitably far apart with respect to the considered distance. More precisely,
one would expect that the distances
\begin{equation} \label{def:heatmapdist}
  \left\| \overline{C}_n^0 - \overline{C}_m^0 \right\|_{L^p} +
    \left\| \overline{C}_n^1 - \overline{C}_m^1 \right\|_{L^p}
  \quad\mbox{ for }\quad
  m,n = 1,\ldots,51
\end{equation}
should be smoothly increasing as the distance~$|m - n|$ becomes
larger. The distances defined in~(\ref{def:heatmapdist}) for $p = 1,2,\infty$
are shown as heat maps in Figure~\ref{fig:DistancesHeatMap}, from top left
to bottom middle. Notice that while in all three cases this ``monotonicity''
can be observed, the level of smoothness of the heat maps decreases with
increasing~$p$. The first of these observations is responsible for the fact
that all three choices of~$p$ lead to acceptable classification schemes, 
while the second observation explains the relative difference in their
performances.

We now turn our attention to the second data set~(D2). In this case, we
consider $R = 50$ simulations each for $N = 6$ different mass values,
and the length of each topological process is still given by $M = 500$.
We increase the size of the training sets to $R_T = 10$, and classify
the remaining $(R-R_T)N = 240$ simulations using the classification
schemes~(CA), (C0), and~(C1). Regardless of the choice of $L^p$-norm,
where again we use $p = 1,2,\infty$, the classification schemes yields 100\%
accuracy. We have also performed cross validations by choosing different
training sets of size ten, which again leads to perfect classifications. 
These results certainly owe to the fact that the discrete mass values~$\mu$
considered in the data set are spaced further apart. However, while in the
case~(D1) the $L^\infty$-classification led to isolated misclassifications 
with large deviations, this was not observed in the data set~(D2). We 
believe that these outliers are controlled due to the increased size of
the training sets.
\begin{table}[tb]
  \centering
  \begin{tabular}{||c|c|c|c|c|c|c||}
  \hline\hline
  \multicolumn{7}{||c||}{{\bf $L^1$-norm}} \\
  \hline
  Batch No.\ & Hits & Missed by $1$ & Missed by $2$ & Missed by $3$ & Missed by $4$
    & Wrong \\
  \hline
  $1$ & $2340$ & $210$ & $0$ & $0$ & $0$ & $0$ \\
  \hline
  $2$ & $2328$ & $222$ & $0$ & $0$ & $0$ & $0$ \\
  \hline\hline
  \multicolumn{7}{||c||}{{\bf $L^2$-norm}} \\
  \hline
  Batch No.\ & Hits & Missed by $1$ & Missed by $2$ & Missed by $3$ & Missed by $4$
    & Wrong \\
  \hline
  $1$ & $2188$ & $361$ & $1$ & $0$ & $0$ & $0$ \\
  \hline
  $2$ & $2163$ & $386$ & $1$ & $0$ & $0$ & $0$ \\
  \hline\hline
  \multicolumn{7}{||c||}{{\bf $L^\infty$-norm}} \\
  \hline
  Batch No.\ & Hits & Missed by $1$ & Missed by $2$ & Missed by $3$ & Missed by $4$
    & Wrong \\
  \hline
  $1$ & $1655$ & $769$ & $113$ & $11$ & $2$ & $0$ \\
  \hline
  $2$ & $1699$ & $721$ & $117$ & $11$ & $2$ & $0$ \\
  \hline\hline
  \end{tabular}
  \caption{Classification results for the data set~(D3) using the classification
           scheme~(CA). The two training set batches consist of~$R_T = 50$
           simulations each. Once the averaged classifier processes~$\overline{C}_n^k$
           for $n = 1,\ldots,N$, where $N = 51$, have been computed, the
           remaining~$(R-R_T)N = 2550$ solution snapshot processes are classified
           using the $L^1$-, $L^2$-, and~$L^\infty$-norms. The table contains
           the number of classification hits as well as detailed information
           on the index difference observed for the misclassifications. The
           column labeled ``Wrong'' contains classifications which missed by
           more than~$4$ indices.}
  \label{tableD3aCA}
\end{table}
\begin{table}[tb]
  \centering
  \begin{tabular}{||c|c|c|c|c|c|c|c||}
  \hline\hline
  \multicolumn{8}{||c||}{{\bf Dimension $k = 0$}} \\
  \hline
  Norm & Hits & Missed by $1$ & Missed by $2$ & Missed by $3$ & Missed by $4$
    & Missed by $5$ & Wrong \\
  \hline
  $L^1$ & $2113$ & $433$ & $4$ & $0$ & $0$ & $0$ & $0$ \\
  \hline
  $L^2$ & $1857$ & $669$ & $24$ & $0$ & $0$ & $0$ & $0$ \\
  \hline
  $L^\infty$ & $1581$ & $781$ & $142$ & $32$ & $7$ & $3$ & $4$ \\
  \hline\hline
  \multicolumn{8}{||c||}{{\bf Dimension $k = 1$}} \\
  \hline
  Norm & Hits & Missed by $1$ & Missed by $2$ & Missed by $3$ & Missed by $4$
    & Missed by $5$ & Wrong \\
  \hline
  $L^1$ & $2308$ & $242$ & $0$ & $0$ & $0$ & $0$ & $0$ \\
  \hline
  $L^2$ & $2210$ & $315$ & $25$ & $0$ & $0$ & $0$ & $0$ \\
  \hline
  $L^\infty$ & $1547$ & $643$ & $235$ & $80$ & $28$ & $10$ & $7$ \\
  \hline\hline
  \end{tabular}
  \caption{Classification results for the data set~(D3) using the classification
           schemes~(C0) and~(C1). The training set consists of~$R_T = 50$ simulations,
           the remaining~$50$ simulations for each of the~$N = 51$ mass values
           are then classified using the $L^1$-, $L^2$-, and~$L^\infty$-norms.
           The table contains the number of classification hits as well as
           detailed information on the index difference observed for the
           misclassifications. The column labeled ``Wrong'' contains
           classifications which missed by more than~$5$ indices.}
  \label{tableD3aC01}
\end{table}

The results so far have used fairly fine sampling in the temporal
direction, which has led to two data sets with topological processes
of length~$M = 500$. It is natural to wonder whether the classification
schemes work as well if this sampling size is decreased. To study this,
we consider the data set~(D3), which consists of $R = 100$ simulations
each for $N = 51$ different mass values, but only saves $M = 25$ solution
snapshots in the interval~$(0,T_e]$. As a first experiment, we use two
different training sets of size~$R_T = 50$, and in each case classify
the remaining $(R-R_T)N = 2550$ simulations. The two different sets
are referred to as ``Batch No.~1'' and ``Batch No.~2'', respectively,
the results for the classifier~(CA) are collected in Table~\ref{tableD3aCA},
while the results for~(C0) and~(C1) can be found in Table~\ref{tableD3aC01}.
\begin{table}[tb]
  \centering
  \begin{tabular}{||c|c|c|c|c|c|c|c||}
  \hline\hline
  \multicolumn{8}{||c||}{{\bf $L^1$-norm}} \\
  \hline
  Batch No.\ & Hits & Missed by $1$ & Missed by $2$ & Missed by $3$ & Missed by $4$
    & Missed by $5$ & Wrong \\
  \hline
  $1$ & $4174$ & $416$ & $0$ & $0$ & $0$ & $0$ & $0$ \\
  \hline
  $2$ & $4179$ & $411$ & $0$ & $0$ & $0$ & $0$ & $0$ \\
  \hline
  $3$ & $4163$ & $427$ & $0$ & $0$ & $0$ & $0$ & $0$ \\
  \hline
  $4$ & $4181$ & $409$ & $0$ & $0$ & $0$ & $0$ & $0$ \\
  \hline
  $5$ & $4143$ & $447$ & $0$ & $0$ & $0$ & $0$ & $0$ \\
  \hline
  $6$ & $4146$ & $444$ & $0$ & $0$ & $0$ & $0$ & $0$ \\
  \hline
  $7$ & $4149$ & $441$ & $0$ & $0$ & $0$ & $0$ & $0$ \\
  \hline
  $8$ & $4198$ & $392$ & $0$ & $0$ & $0$ & $0$ & $0$ \\
  \hline
  $9$ & $4173$ & $416$ & $1$ & $0$ & $0$ & $0$ & $0$ \\
  \hline
  $10$ & $4175$ & $415$ & $0$ & $0$ & $0$ & $0$ & $0$ \\
  \hline\hline
  \multicolumn{8}{||c||}{{\bf $L^2$-norm}} \\
  \hline
  Batch No.\ & Hits & Missed by $1$ & Missed by $2$ & Missed by $3$ & Missed by $4$
    & Missed by $5$ & Wrong \\
  \hline
  $1$ & $3859$ & $719$ & $12$ & $0$ & $0$ & $0$ & $0$ \\
  \hline
  $2$ & $3860$ & $724$ & $6$ & $0$ & $0$ & $0$ & $0$ \\
  \hline
  $3$ & $3809$ & $771$ & $10$ & $0$ & $0$ & $0$ & $0$ \\
  \hline
  $4$ & $3871$ & $713$ & $6$ & $0$ & $0$ & $0$ & $0$ \\
  \hline
  $5$ & $3747$ & $835$ & $8$ & $0$ & $0$ & $0$ & $0$ \\
  \hline
  $6$ & $3852$ & $734$ & $4$ & $0$ & $0$ & $0$ & $0$ \\
  \hline
  $7$ & $3846$ & $741$ & $3$ & $0$ & $0$ & $0$ & $0$ \\
  \hline
  $8$ & $3866$ & $721$ & $3$ & $0$ & $0$ & $0$ & $0$ \\
  \hline
  $9$ & $3843$ & $743$ & $4$ & $0$ & $0$ & $0$ & $0$ \\
  \hline
  $10$ & $3858$ & $729$ & $3$ & $0$ & $0$ & $0$ & $0$ \\
  \hline\hline
  \multicolumn{8}{||c||}{{\bf $L^\infty$-norm}} \\
  \hline
  Batch No.\ & Hits & Missed by $1$ & Missed by $2$ & Missed by $3$ & Missed by $4$
    & Missed by $5$ & Wrong \\
  \hline
  $1$ & $3019$ & $1350$ & $205$ & $11$ & $4$  & $1$ & $0$ \\
  \hline
  $2$ & $3004$ & $1373$ & $185$ & $23$ & $5$  & $0$ & $0$ \\
  \hline
  $3$ & $2972$ & $1395$ & $207$ & $14$ & $2$  & $0$ & $0$ \\
  \hline
  $4$ & $2899$ & $1417$ & $255$ & $14$ & $5$  & $0$ & $0$ \\
  \hline
  $5$ & $2903$ & $1430$ & $236$ & $17$ & $4$  & $0$ & $0$ \\
  \hline
  $6$ & $2954$ & $1412$ & $195$ & $23$ & $5$  & $1$ & $0$ \\
  \hline
  $7$ & $2934$ & $1403$ & $234$ & $16$ & $3$  & $0$ & $0$ \\
  \hline
  $8$ & $2895$ & $1465$ & $213$ & $14$ & $3$  & $0$ & $0$ \\
  \hline
  $9$ & $2908$ & $1432$ & $231$ & $17$ & $2$  & $0$ & $0$ \\
  \hline
  $10$ & $2977$ & $1382$ & $203$ & $23$ & $4$  & $1$ & $0$ \\
  \hline\hline
  \end{tabular}
  \caption{Classification results for the data set~(D3) using the classification
           scheme~(CA). The ten training set batches consist of~$10$ simulations
           each. Once the averaged classifier processes~$\overline{C}_n^k$ for
           $n = 1,\ldots,N$, where $N = 51$, have been computed, the
           remaining~$90N = 4590$ solution snapshot processes are classified
           using the $L^1$-, $L^2$-, and~$L^\infty$-norms. The table contains
           the number of classification hits as well as detailed information
           on the index difference observed for the misclassifications. The
           column labeled ``Wrong'' contains classifications which missed by
           more than~$5$ indices.}
  \label{tableD3bCA}
\end{table}

For the classification scheme~(CA), the simulation results are convincing.
If one uses the $L^1$-norm for the classification, then in over~91\% of the
cases the classification scheme finds the correct mass value class, and all
classifications are correct up to an error of the form $\mu \pm 0.01$. As before,
the $L^2$-norm classifier performs slightly worse, with correct classifications
in over 84\% of the cases, and errors of~$\mu \pm 0.01$ in over 99\% of the
cases. The maximally observed errors where of the form $\mu \pm 0.02$.
While the integral norms perform extremely well, the maximum norm
$L^\infty$-estimator again does not seem to be quite as accurate. It leads
to correct classifications only in at least 64\% of the cases, errors
of the form $\mu \pm 0.01$ in at least~94\%, as well as errors $\mu \pm 0.02$
in at least 99\% of the cases. We would like to point out, however, that
in all of the above situations, the performance exceeded the one observed
in Table~\ref{tableD1CA} for the data set~(D1), despite the fact that now
our topological processes have only length $M = 25$. This increased 
accuracy is due to the larger training set size $R_T = 50$.

How do the results change if we use the dimension-dependent classification
schemes~(C0) or~(C1)? The corresponding results are presented in
Table~\ref{tableD3aC01}, and they largely confirm the observations made
in the context of Table~\ref{tableD1C01}. Using only one dimension for the
classification yields still acceptable, but definitely worse results. While
now the $L^1$-classification is clearly the most accurate, the $L^2$- and
$L^\infty$-norm classifications come in second and third, respectively.
Finally, the classification using dimension $k = 1$ performs better
than the one in dimension zero, at least for the case of $L^1$-classifications.
These results show, however, that in general it is preferable to use all
available topological information for the classification.

As a final experiment, we again consider data set~(D3), but this time
the $R = 100$ simulations for each mass value~$\mu$ are partitioned into
ten different training sets of size~$R_T = 10$ each. These training sets
are referred to as ``Batch No.~$\ell$'' for $\ell = 1,\ldots,10$, and for
each training set the remaining $(R-R_T)N = 4590$ topological processes
are then classified using scheme~(CA). The results of these computations
are shown in Table~\ref{tableD3bCA}, and they confirm our previous
observations. For classifications using the $L^1$-, $L^2$-, and $L^\infty$-norm
one obtains exact hits in at least 90\%, 81\%, and 63\% of the cases, maximal
errors of the form $\mu \pm 0.01$ are observed in at least 99\%, 99\%,
and 94\% of the cases, and finally maximal errors $\mu \pm 0.02$ occur
in at least 100\%, 100\%, and 99\% of the cases, respectively. For the
$L^1$-, $L^2$-, and $L^\infty$-norm classifications, the maximal observed
classification errors where $\mu \pm 0.02$, $\mu \pm 0.02$, and $\mu \pm 0.05$,
respectively. While all of these numbers are slightly decreased from the ones
shown in Table~\ref{tableD3aCA}, this is to be expected due to the
smaller training set size $R_T = 10$. Nevertheless, especially the
$L^1$-norm classification performs amazingly accurate.
\subsection{Recovering the Snapshot Time via Topology}
\label{sec:chctime}
We have seen in the last section that if one considers all the topological
information encoded in the microstructure evolution of a solution path of the
Cahn-Hilliard-Cook model~(\ref{chc}), then one can use the classification
schemes~(CA) and~(Ck) to recover the total mass~$\mu$ with high accuracy.
Of course the total mass is only one of several parameters in~(\ref{chc})
which affect the topology, and it would be interesting to also determine
how strong the influence of these other parameters on the topology is.

Rather than pursuing this line of inquiry, we now consider a different
question. For this, consider again the simulations of the Cahn-Hilliard-Cook
model, but now for a fixed and a-priori known mass value~$\mu$. Then a
simulation of~(\ref{chc}) produces solution snapshots along~$M$ equidistant
timesteps~$t_1, \ldots, t_M$, and for each time~$t_\ell$, computational
topology can be used to create the associated persistence landscapes in
dimensions~$0$ and~$1$. Suppose now that we pick up one of these persistence
landscapes at random. Can we recover the actual time~$t_\ell$ which corresponds
to this landscape purely from the encoded topological information? In other
words, for fixed mass value~$\mu$, is it possible to determine the \emph{time
at which a solution snapshot has been taken\/}, purely based on
topological information? We will see in the following that this can
in fact be done, and in many cases with surprising accuracy.

As we saw in the previous section, our classification technique can only
recover mass information up to a certain tolerance, and the situation is
no different in the snapshot time case. Therefore, we only use the data
set~(D3) in the present section, which produces~$M = 25$ times~$t_\ell$
in the interval~$(0,T_e]$. Recall that this data set considers the $N = 51$
mass values $\mu = 0.01 \, n$ for $n=0,\ldots,50$, and for each mass value
one has~$R = 100$ stochastically independent repetitions of the simulation.

Based on the composition of the data set~(D3), we apply the classification
schemes~(C0) and~(C1) from Section~\ref{sec:basicmethod} for every fixed
mass value~$\mu$ to the snapshot time classes~$C_1, \ldots, C_M$, which
correspond to the~$M$ snapshot times~$t_1, \ldots, t_M$ in increasing
order. Then we proceed as follows.
\begin{itemize}
\item Based on~$R_T < R = 100$ training runs each from the~$M = 25$ considered
snapshot times~$t_\ell$, we determine the averaged classifier~$\overline{C}_m^k$
of the~$m^{\rm th}$ class in dimension~$k$ as defined in~(\ref{aveclassifierK})
with $S = R_T$ and $n = m = 1,\ldots,M$.
\item We then try to classify each of the remaining~$(R - R_T) M$ simulations
in the data set for fixed mass value~$\mu$, by either using the classification
scheme~(C0) or~(C1).
\end{itemize}
As before, in each case we note the index~$m$ obtained from the classification
scheme, and compare it to the snapshot time index~$m_{true}$ of the underlying
simulation. If~$m = m_{true}$, then we have obtained a ``Hit'', otherwise the
classification ``missed by~$|m - m_{true}|$''. Clearly the integer
value~$|m - m_{true}|$ is directly proportional to the time difference
incurred by the classification, with the proportionality factor being
given by~$T_e / M = 0.04 \, T_e$.
\begin{figure} \centering
  \setlength{\unitlength}{1 cm}
  \begin{picture}(16.75,16.0)
    \put(0.0,11.0){
      \includegraphics[height=5cm]{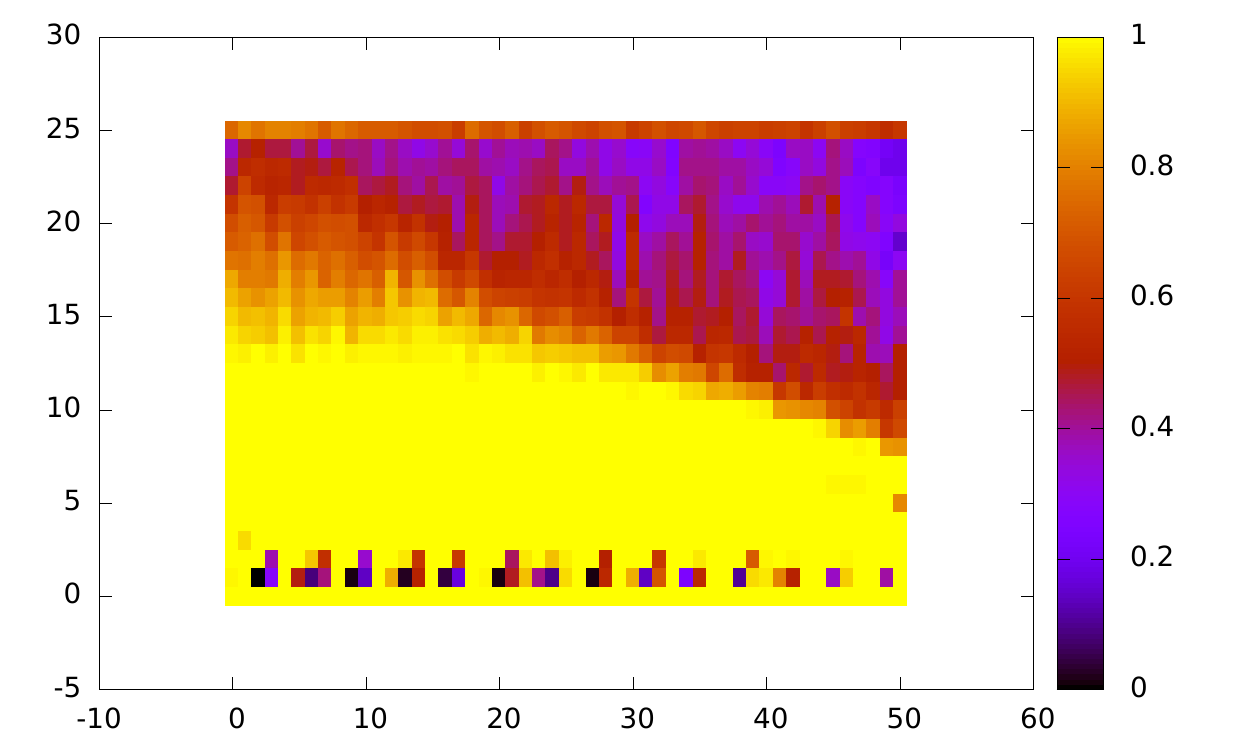}}
    \put(8.75,11.0){
      \includegraphics[height=5cm]{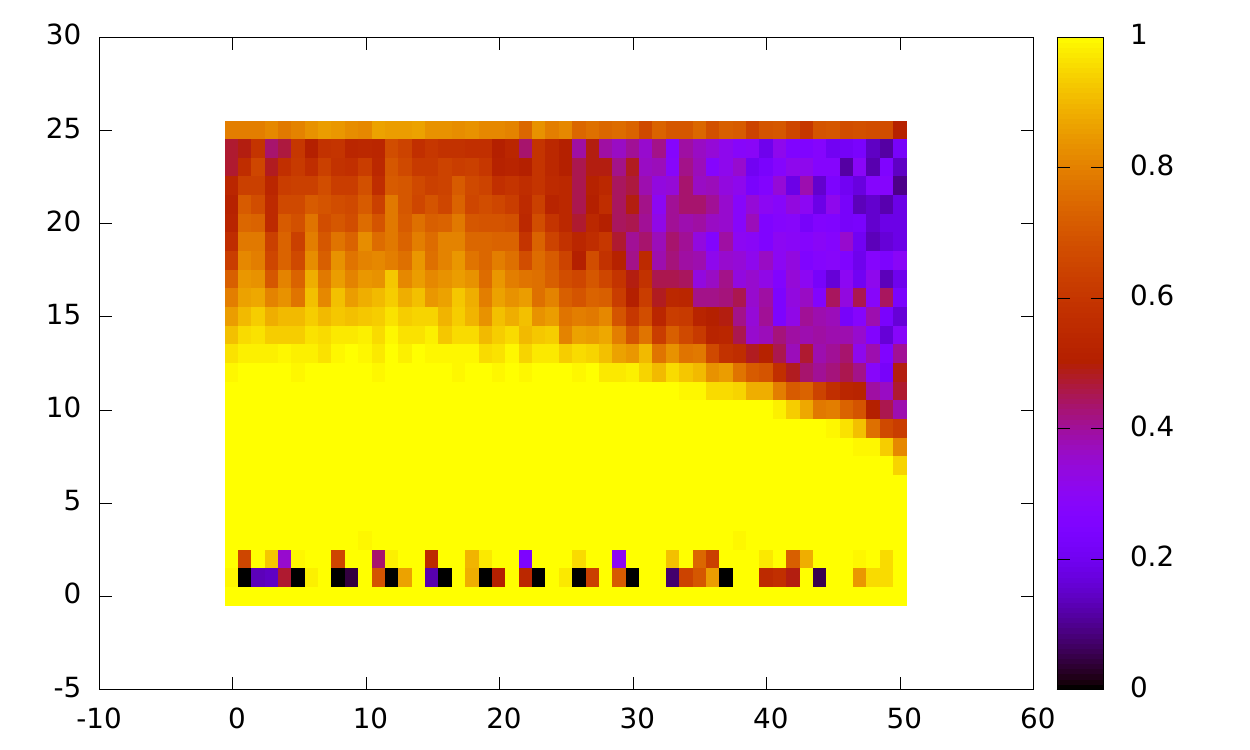}}
    \put(0.0,5.5){
      \includegraphics[height=5cm]{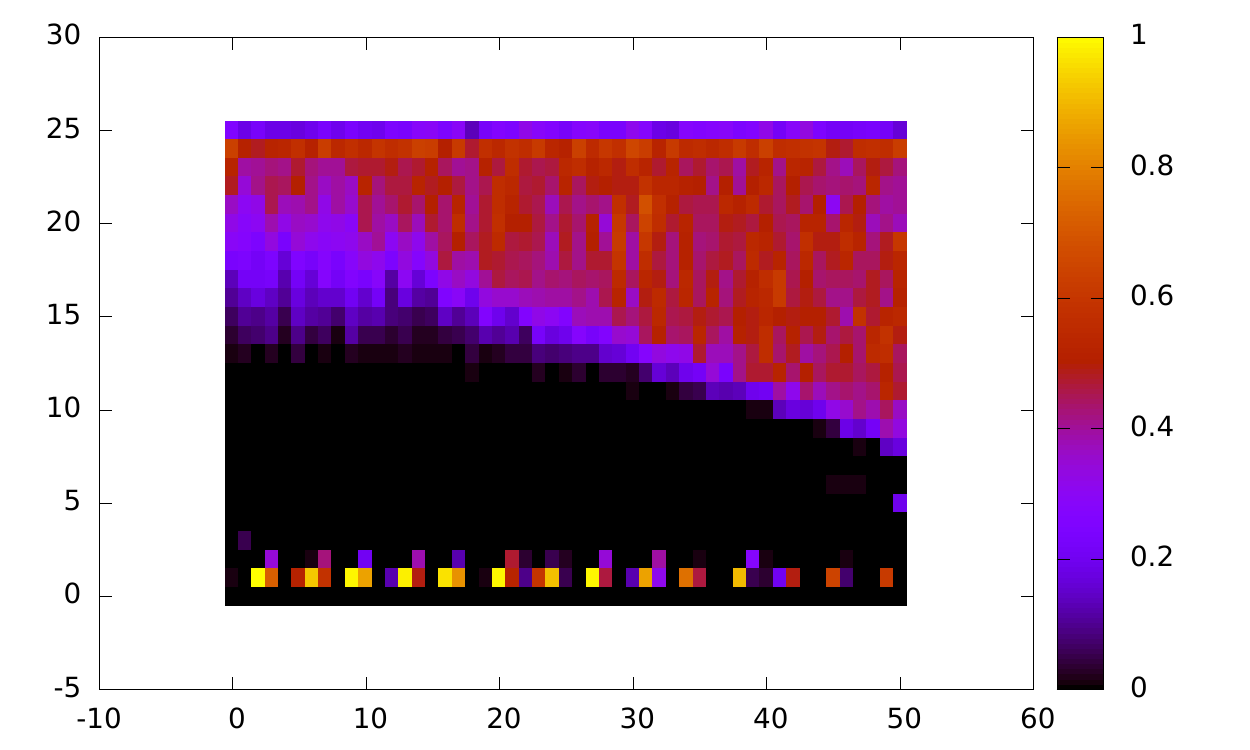}}
    \put(8.75,5.5){
      \includegraphics[height=5cm]{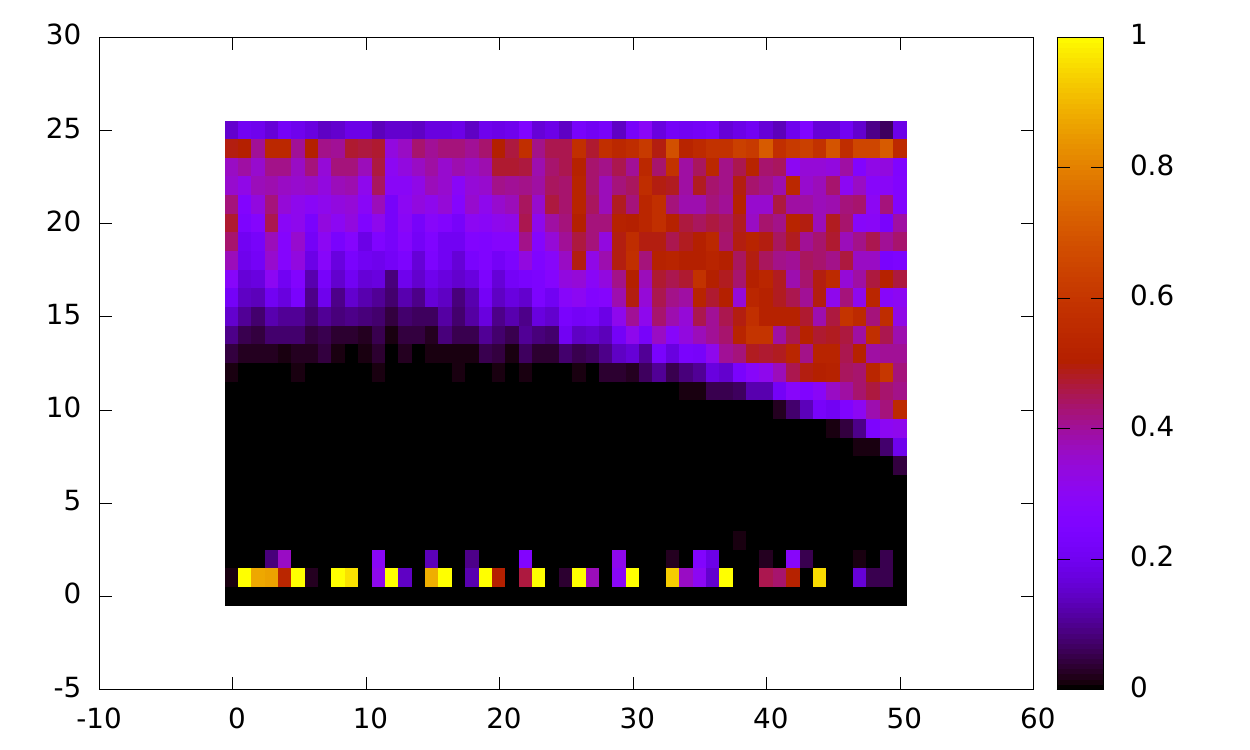}}
    \put(0.0,0.0){
      \includegraphics[height=5cm]{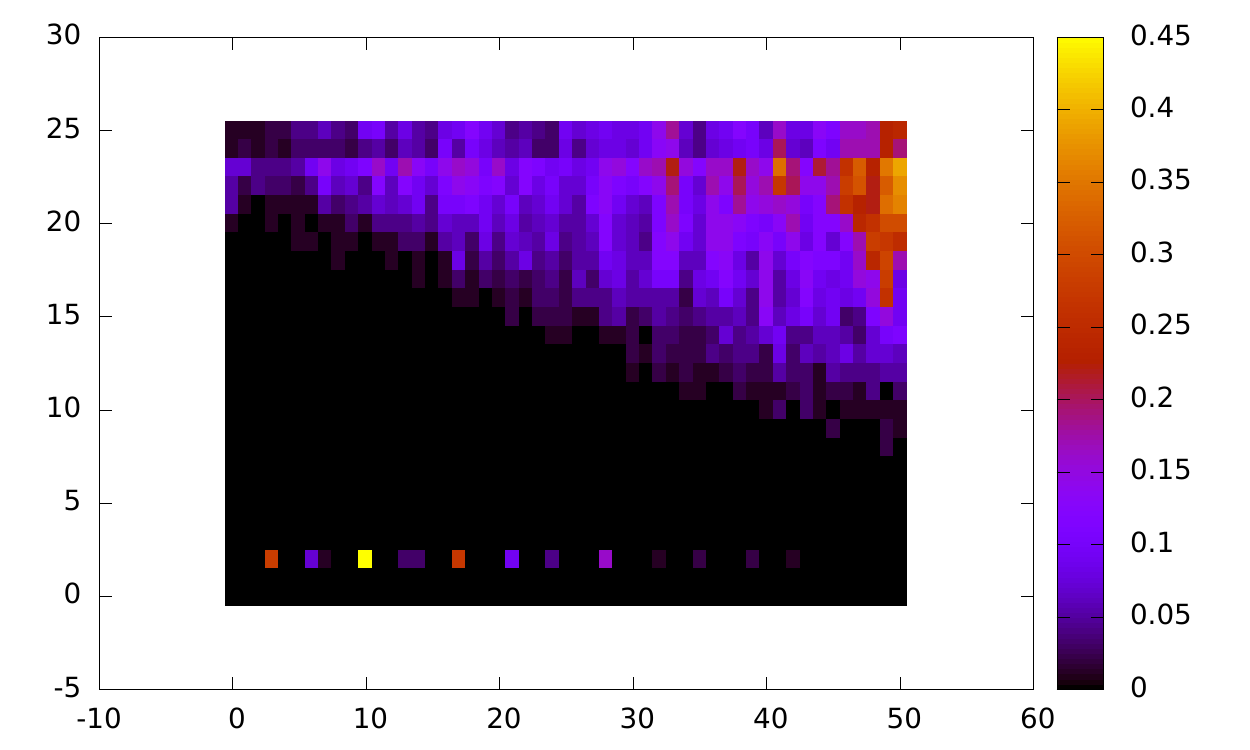}}
    \put(8.75,0.0){
      \includegraphics[height=5cm]{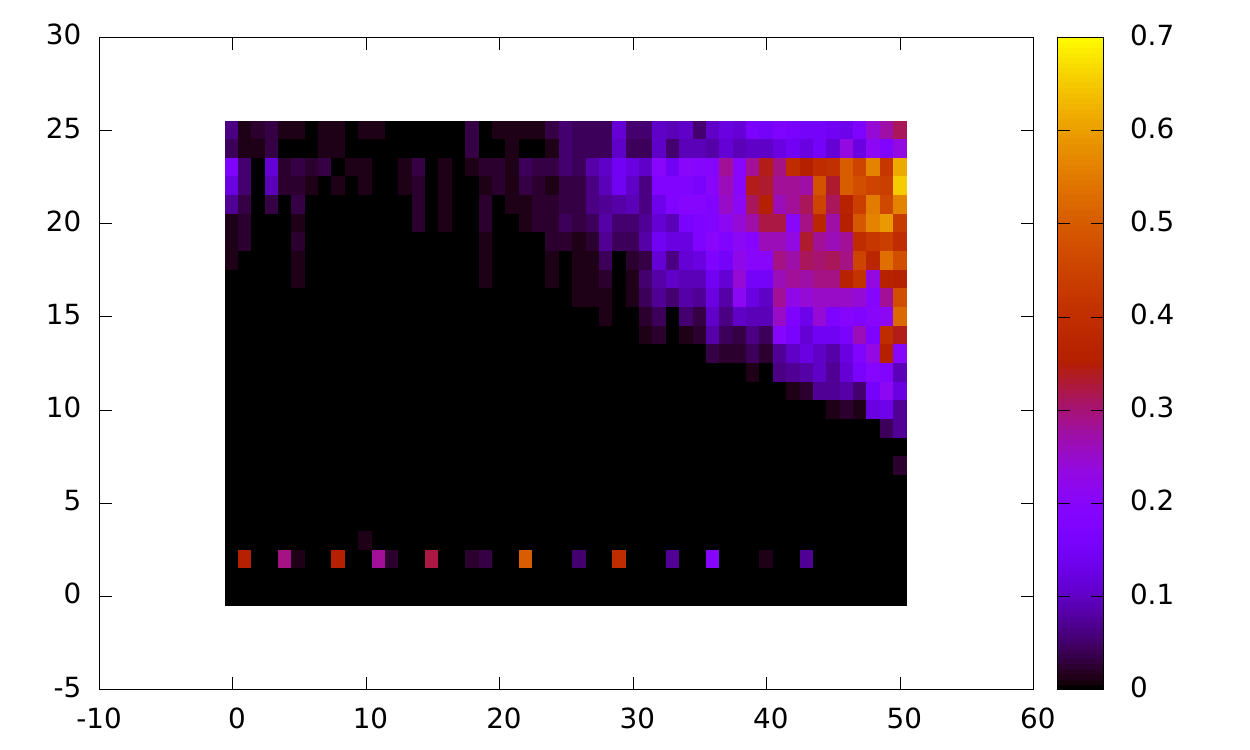}}
  \end{picture}
  \caption{Recovering the snapshot time from persistence information
           using the $L^1$-norm and training set size~$R_T = 50$. Each
           panel contains a heat map, whose horizontal axis corresponds
           to the~$N = 51$ mass values~$\mu$, and the vertical axis
           to the~$M = 25$ snapshot times, in increasing order. Color
           indicates frequency of observation, as indicated by the colorbar.
           Panels in the left column are for scheme~(C0), the right
           column corresponds to~(C1). From top to bottom, the three
           rows show the likelihoods of exact hits, misses by exactly one,
           and misses by at least two.}
  \label{fig:time1}
\end{figure}
\begin{figure} \centering
  \setlength{\unitlength}{1 cm}
  \begin{picture}(16.75,16.0)
    \put(0.0,11.0){
      \includegraphics[height=5cm]{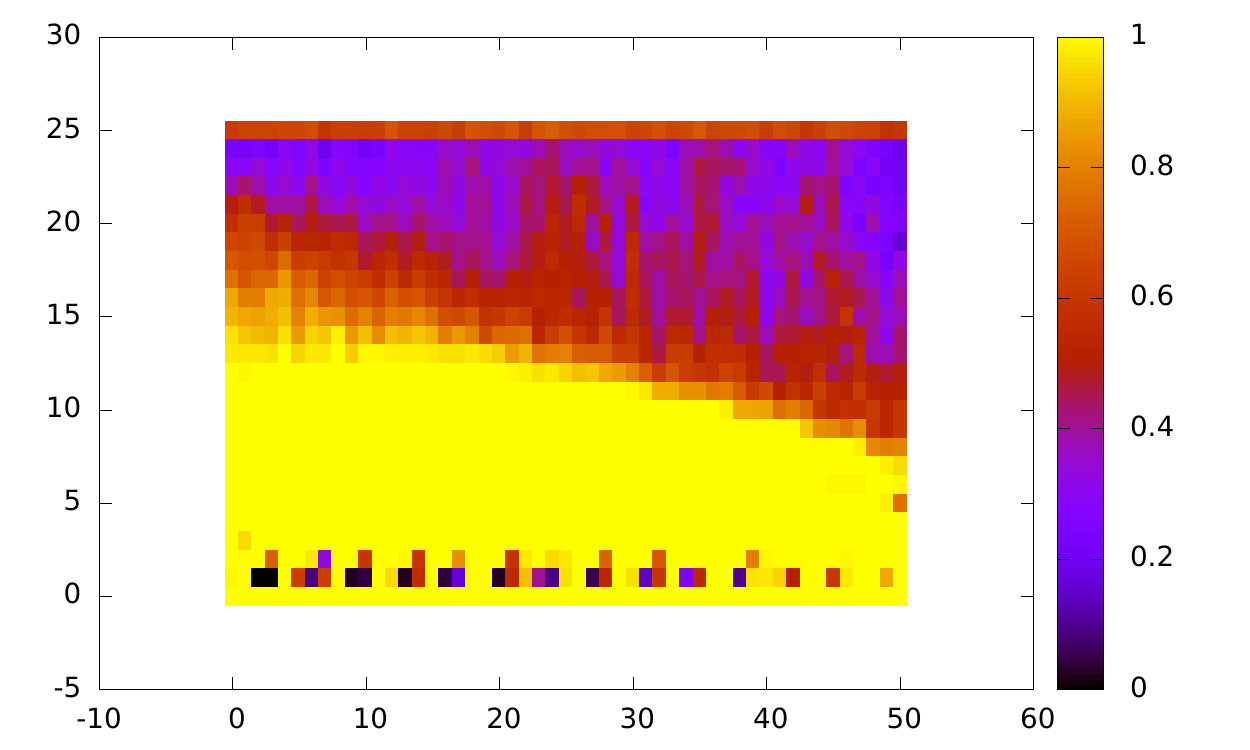}}
    \put(8.75,11.0){
      \includegraphics[height=5cm]{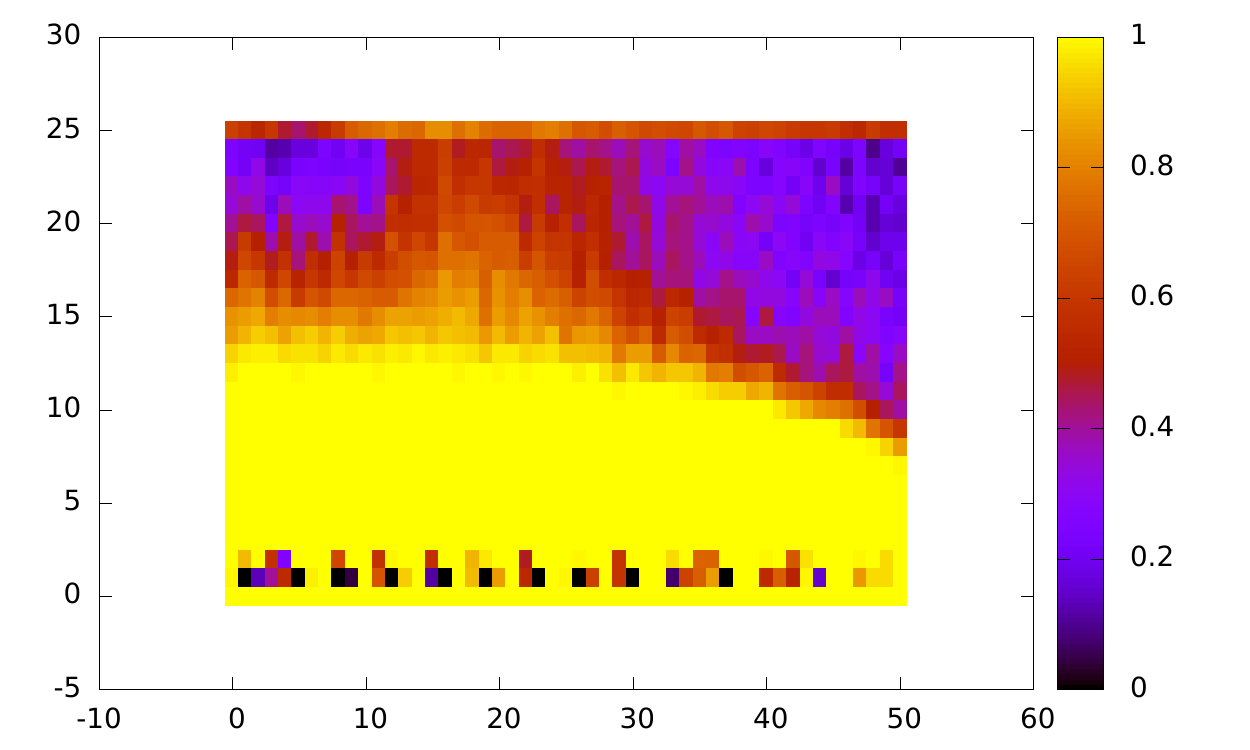}}
    \put(0.0,5.5){
      \includegraphics[height=5cm]{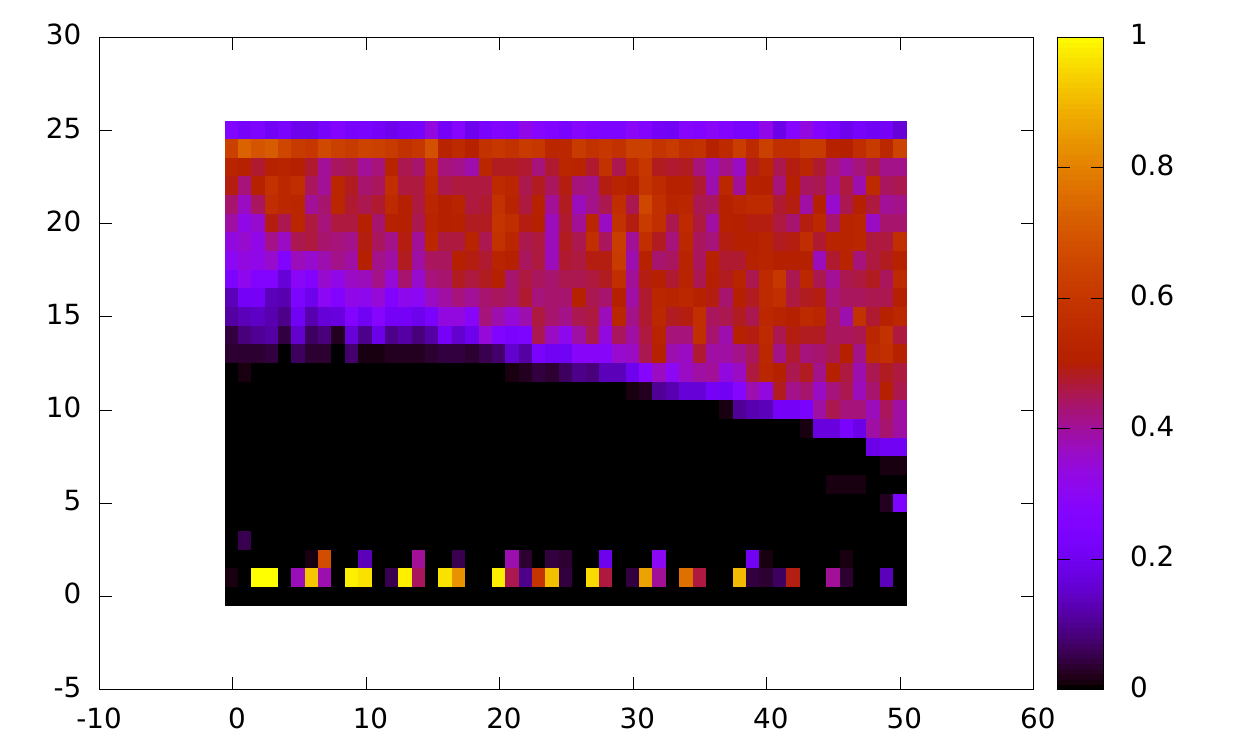}}
    \put(8.75,5.5){
      \includegraphics[height=5cm]{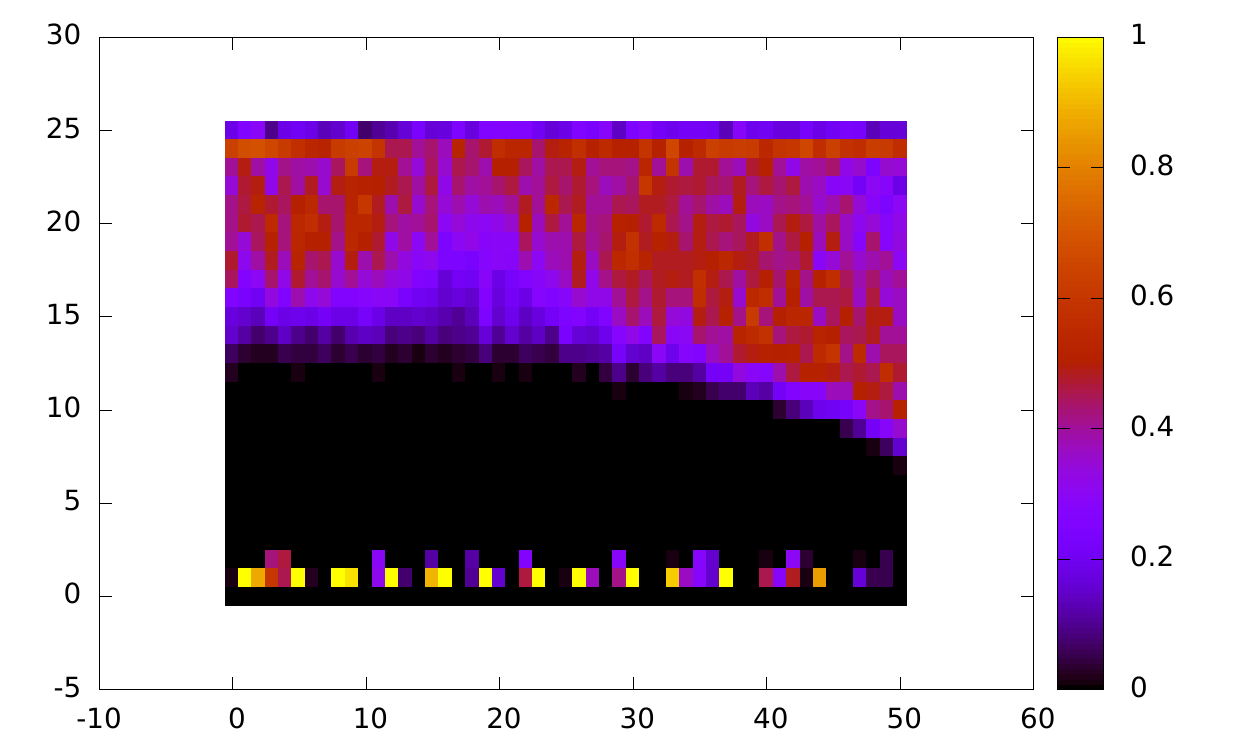}}
    \put(0.0,0.0){
      \includegraphics[height=5cm]{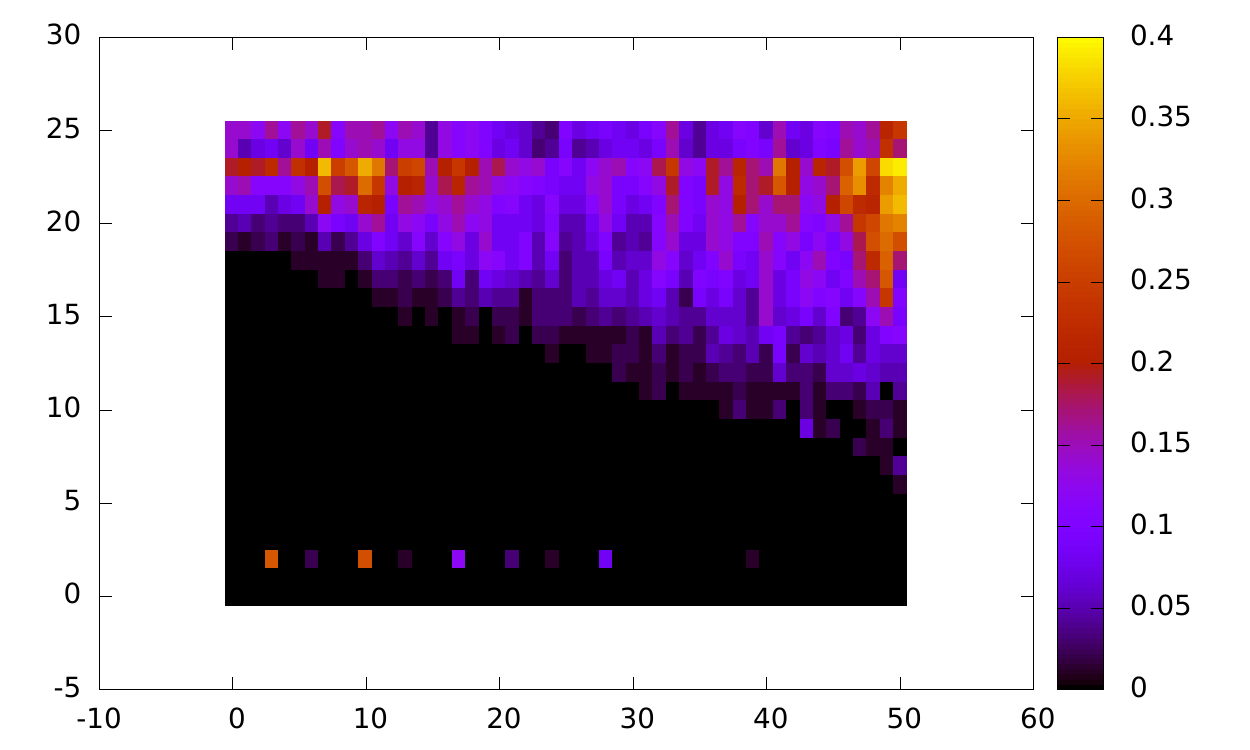}}
    \put(8.75,0.0){
      \includegraphics[height=5cm]{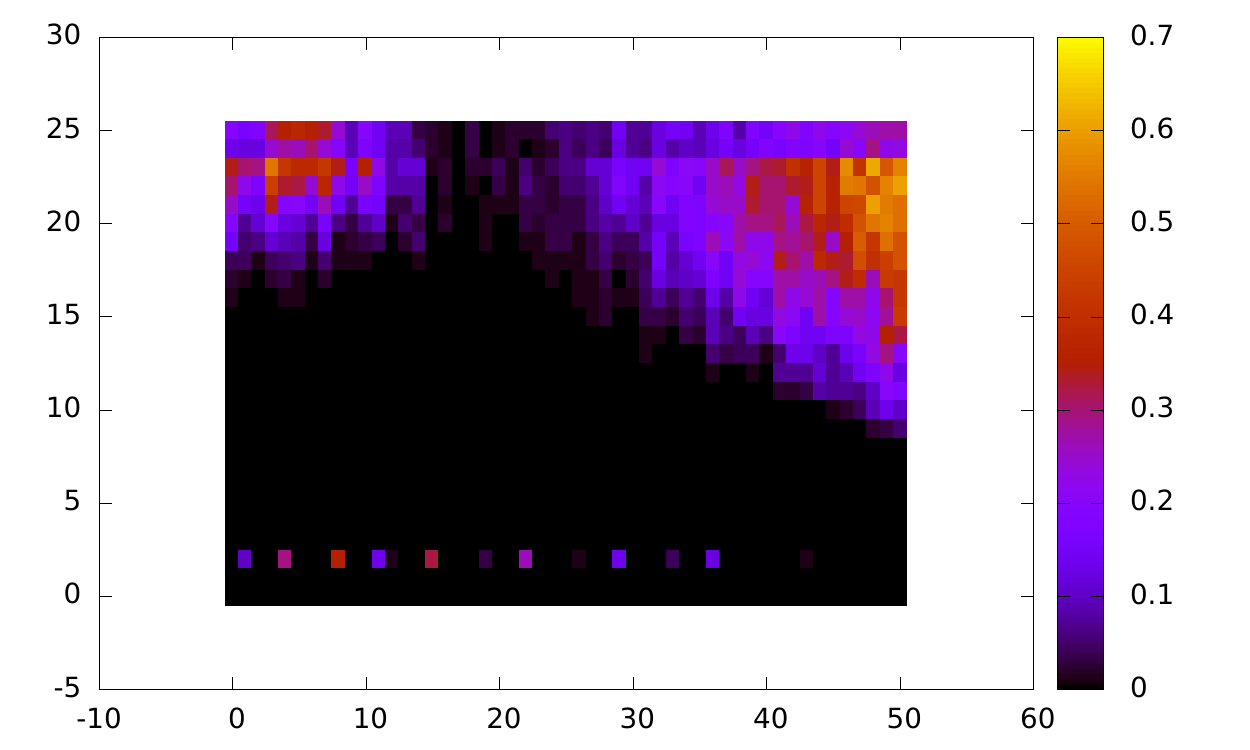}}
  \end{picture}
  \caption{Recovering the snapshot time from persistence information
           using the $L^2$-norm and training set size~$R_T = 50$. Each
           panel contains a heat map, whose horizontal axis corresponds
           to the~$N = 51$ mass values~$\mu$, and the vertical axis
           to the~$M = 25$ snapshot times, in increasing order. Color
           indicates frequency of observation, as indicated by the colorbar.
           Panels in the left column are for scheme~(C0), the right
           column corresponds to~(C1). From top to bottom, the three
           rows show the likelihoods of exact hits, misses by exactly one,
           and misses by at least two.}
  \label{fig:time2}
\end{figure}

Rather than performing an exhaustive study as in the previous section,
we use the already obtained insight to focus our approach. Since the
$L^\infty$-norm classification consistently performed significantly worse
than the $L^1$- and $L^2$-norm classifications, we only run simulations 
for the latter two. However, in contrast to the experiments presented in
Section~\ref{sec:chcmass}, it turns out that for recovering the snapshot
time both~(C0) and~(C1) are superior to~(CA). For this reason, we only
consider the dimension-dependent classification schemes below.

For the case of $L^1$-norm classifications using the schemes~(C0)
and~(C1), and for training set size~$R_T= 50$, the results can be found in
Figure~\ref{fig:time1}. For the sake of brevity, we have presented the
data in the form of heat maps, rather than providing precise tables.
Each panel in the figure represents frequencies which are indicated 
by colors, and which correspond to the ranges of the adjacent colorbars.
The horizontal axes correspond to the~$N = 51$ mass values~$\mu$, while
the vertical axes are for the~$M = 25$ snapshot times, both in increasing
order. Panels in the left column are for classification scheme~(C0), while
the right column corresponds to~(C1). From top to bottom, the three rows
show the likelihoods of \emph{exact hits\/}, \emph{misses by exactly one\/},
and \emph{misses by at least two\/}. Analogous results for the case of the
$L^2$-norm are depicted in Figure~\ref{fig:time2}. In each of the cases we
have cross-validated the results by choosing different training sets of the
same size, and the results are quantitatively similar. Notice that in this
setting, for each mass value~$\mu$ we are classifying $(R-R_T) M = 1250$
microstructures according to their snapshot time class.

A first glance at the heat maps in Figures~\ref{fig:time1} and~\ref{fig:time2},
especially the ones in the respective bottom rows, confirms our statement from
the beginning of this section. Averaged persistence landscapes can indeed be
used to accurately deduce the snapshot time at which a given microstructure
occurs. In most cases, the match is an exact hit, while in a few cases one
is off by one index. In more detail, one can draw the following conclusions
from the data in these figures:
\begin{itemize}
\item For every considered mass value~$\mu$, there exists an interval
of the form~$[\alpha_\mu, \beta_\mu] \subset (0,T_e]$ over which the
snapshot time can be determined exactly with probability close to one.
\item The times~$\alpha_\mu > 0$ are small, and more or less
constant with~$\mu$. Only at the beginning of the simulations does
it seem to be difficult to pin down the exact snapshot time. The volatile
timeframe~$(0,\alpha_\mu)$ corresponds to the initial rapid smoothing regime
in the parabolic stochastic partial differential equation~(\ref{chc}),
during with the high-frequency terms in the random initial conditions
decay due to dissipativity.
\item The value of~$\beta_\mu$ is close to about~60\% of the
endtime~$T_e$ for mass values~$\mu$ close to zero, and decreases
as~$\mu$ increases. At the final mass value $\mu = 0.5$ one has
$\beta_\mu \approx 0.3 \, T_e$. For larger simulation
times, i.e., times in the interval~$(\beta_\mu, T_e]$, it
appears to be more difficult to precisely locate the snapshot
time, and this becomes more pronounced as the mass~$\mu$ increases.
This timeframe is well into the coarsening regime, where microstructures,
particularly of droplet type, only move slightly, and no longer change
their topology frequently. Thus, the classification results are not
surprising.
\end{itemize}
Notice that in all panels of the two figures, the probabilities seem
to jump significantly along the topmost line. In other words, the
classification results for the time~$t_M$ are significantly different
from the one for~$t_{M-1}$. This is due to the fact that for all 
interior times~$t_\ell < T_e$ one can miss by one by either classifying
the time as~$t_{\ell-1}$ or~$t_{\ell+1}$, while at the last time, anything
that would have been misclassified as a larger time will most likely be
classified as~$t_M$.
\begin{figure} \centering
  \setlength{\unitlength}{1 cm}
  \begin{picture}(16.75,16.0)
    \put(0.0,11.0){
      \includegraphics[height=5cm]{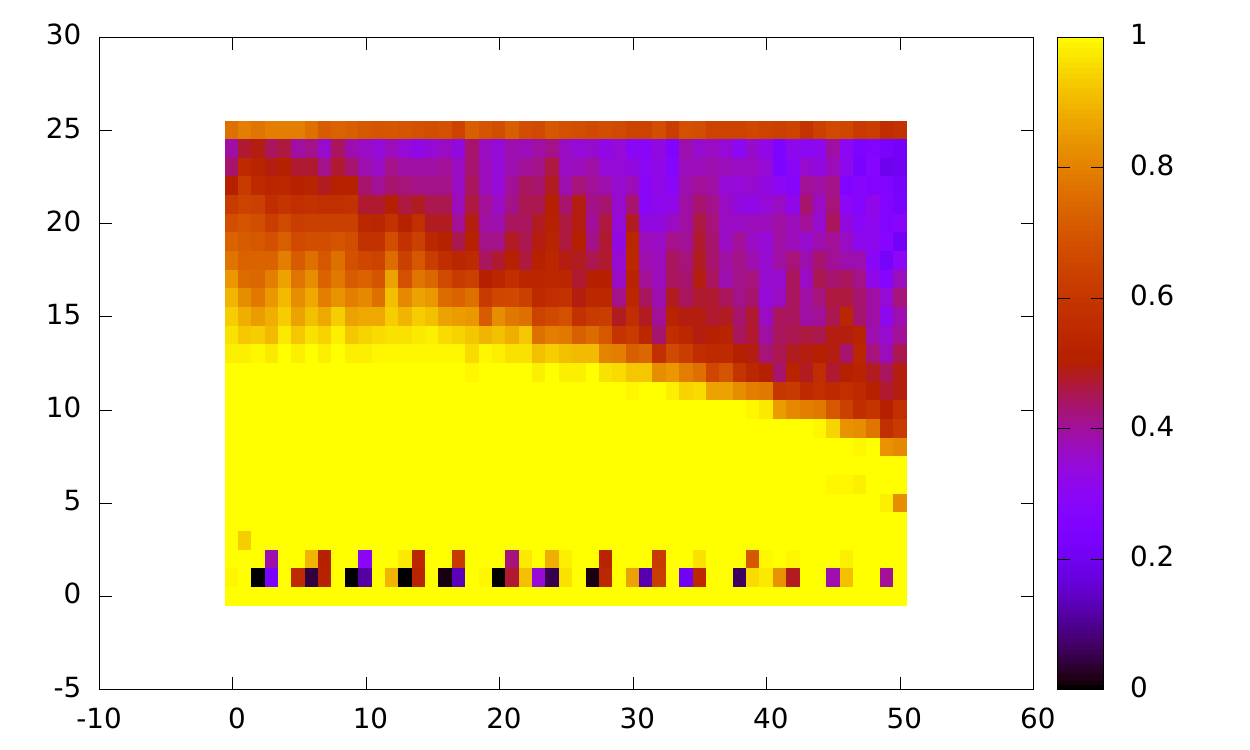}}
    \put(8.75,11.0){
      \includegraphics[height=5cm]{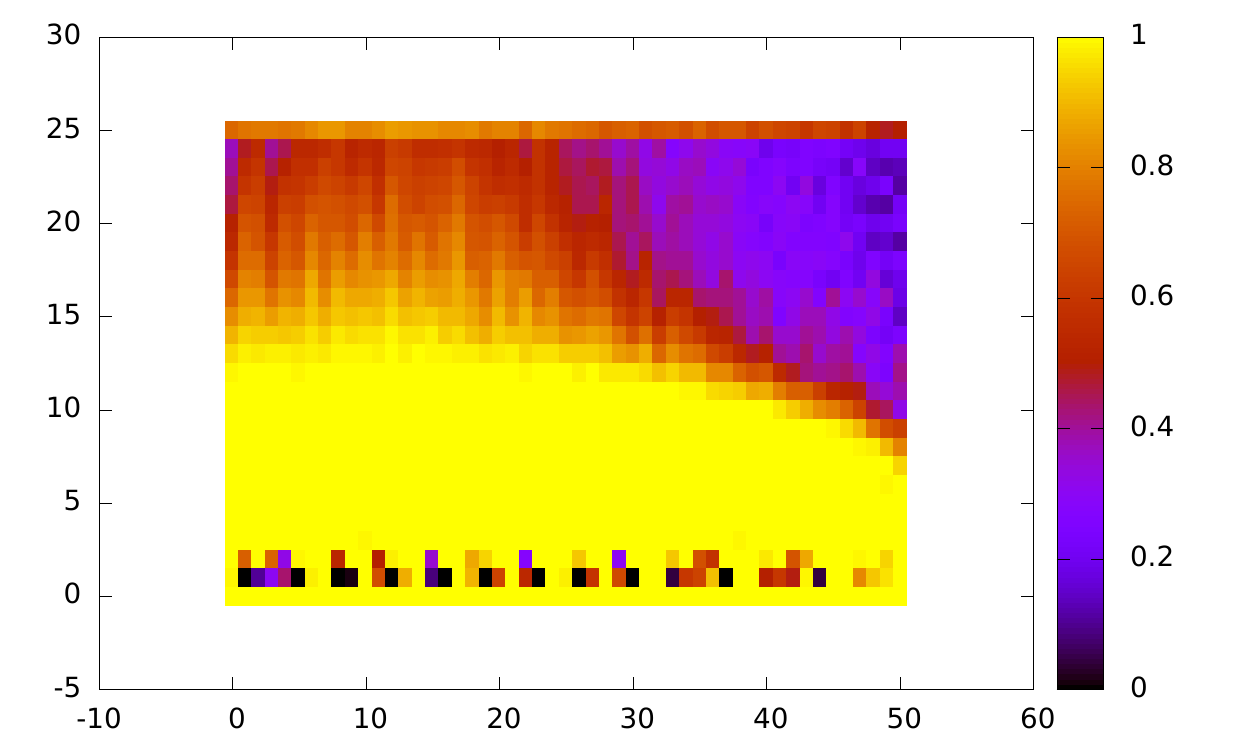}}
    \put(0.0,5.5){
      \includegraphics[height=5cm]{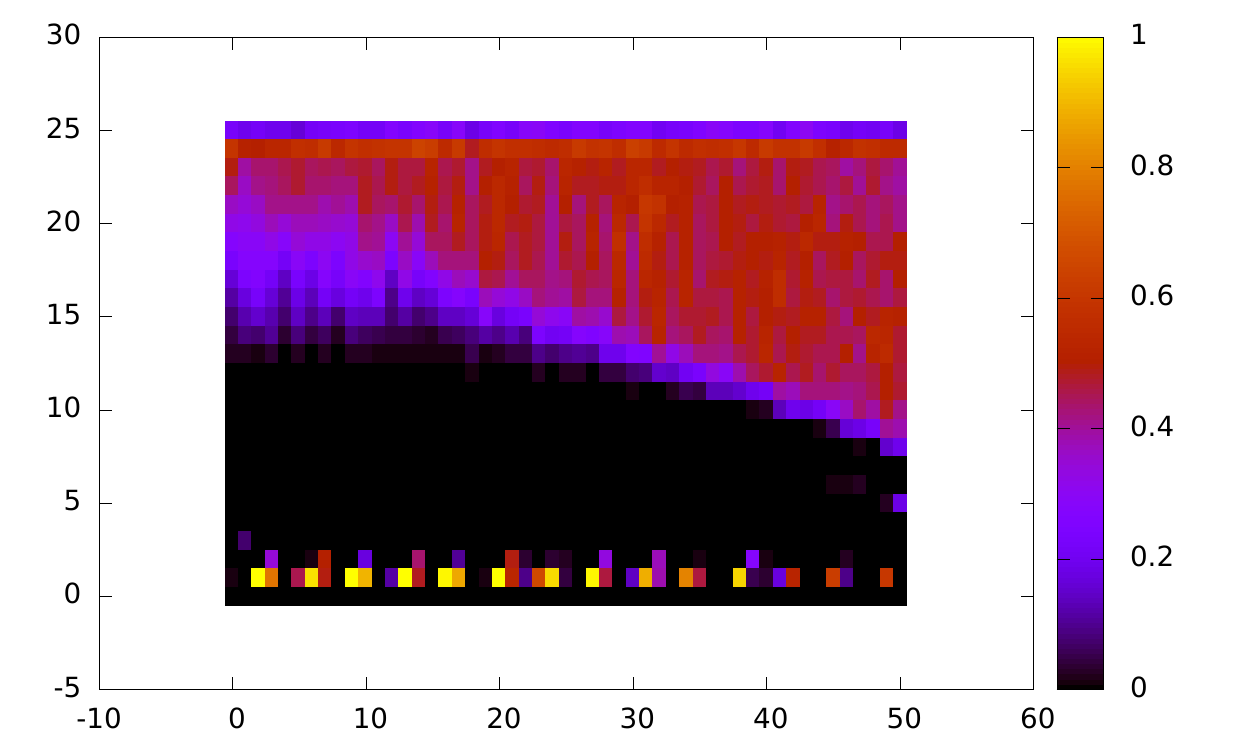}}
    \put(8.75,5.5){
      \includegraphics[height=5cm]{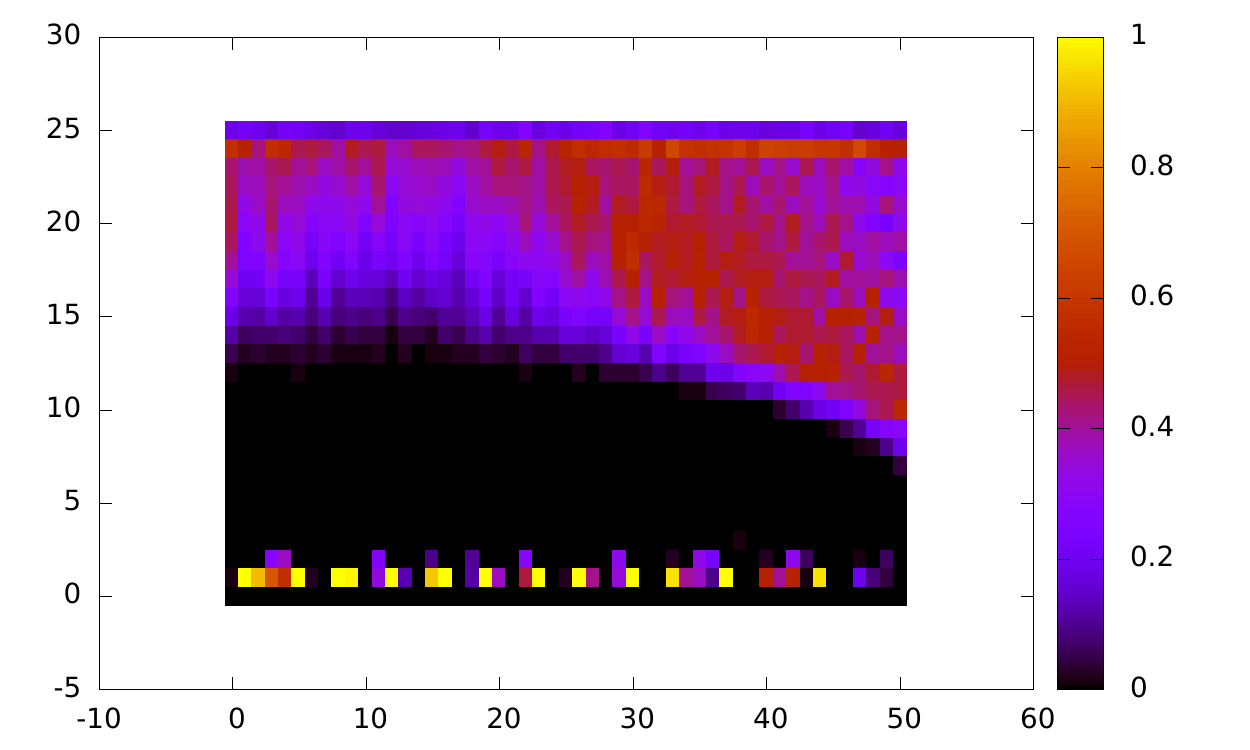}}
    \put(0.0,0.0){
      \includegraphics[height=5cm]{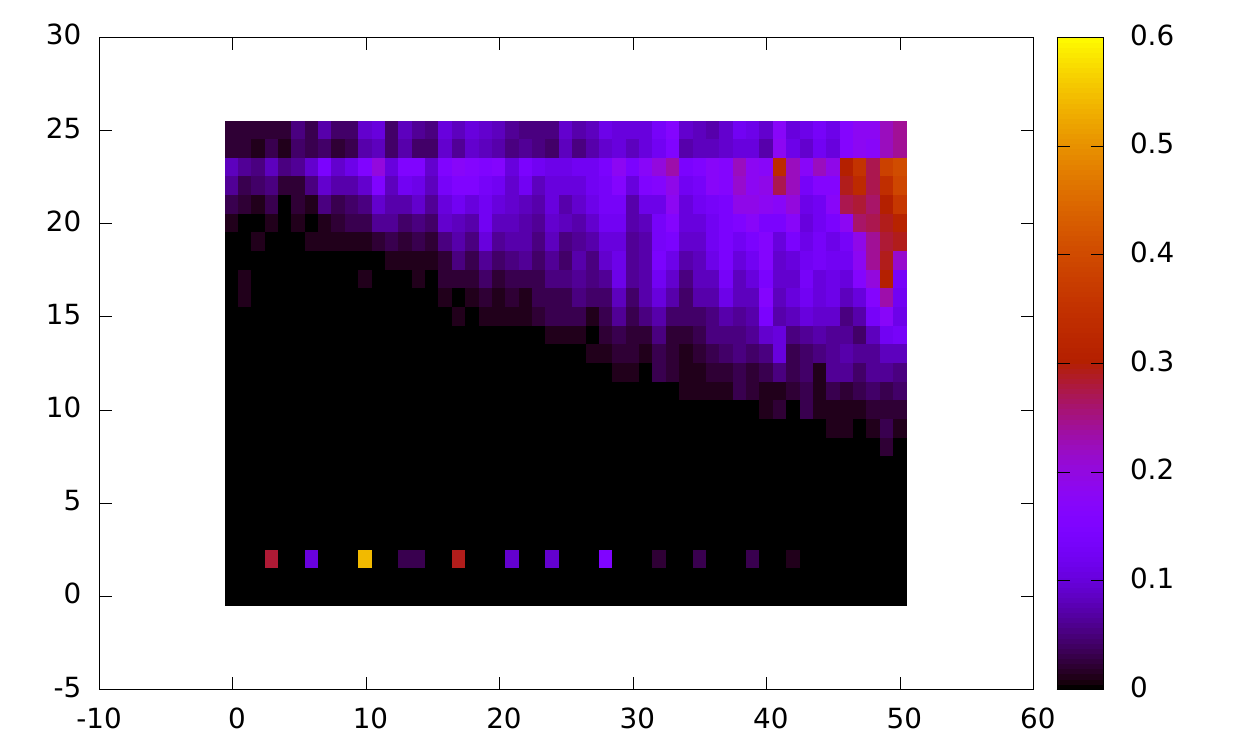}}
    \put(8.75,0.0){
      \includegraphics[height=5cm]{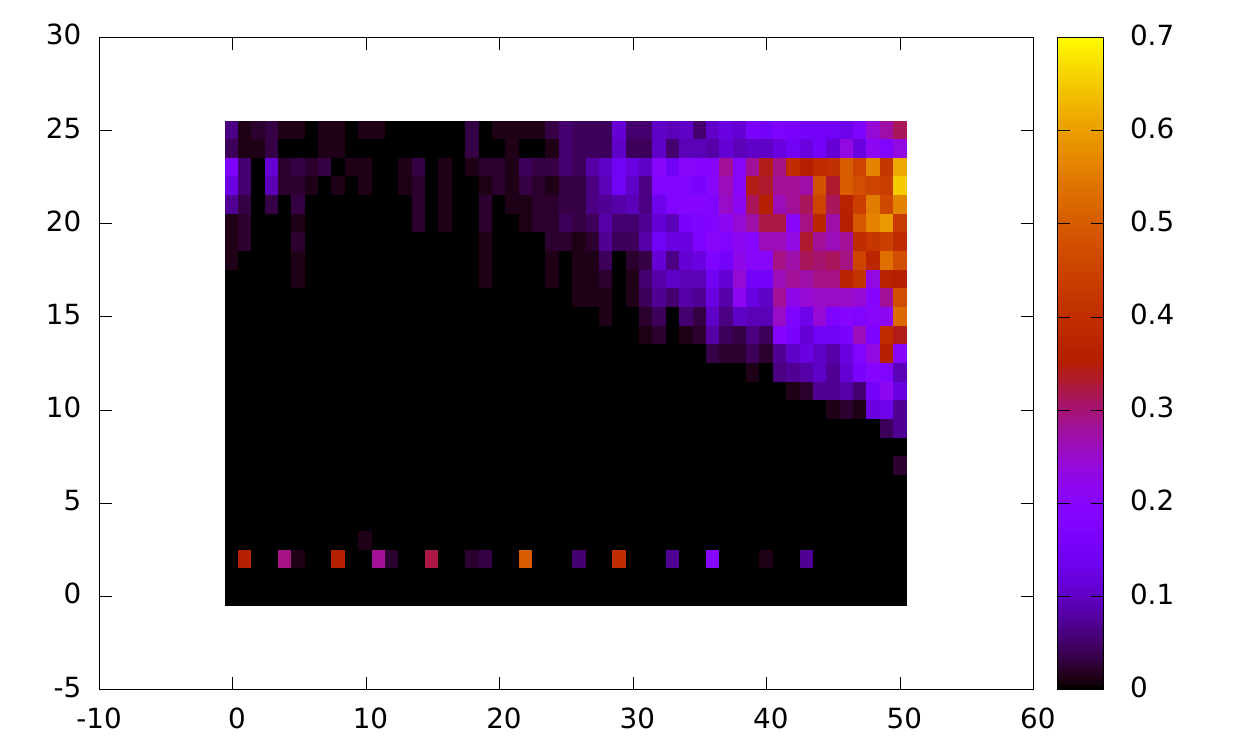}}
  \end{picture}
  \caption{Recovering the snapshot time from persistence information
           using the $L^1$-norm and training set size~$R_T = 20$. Each
           panel contains a heat map, whose horizontal axis corresponds
           to the~$N = 51$ mass values~$\mu$, and the vertical axis
           to the~$M = 25$ snapshot times, in increasing order. Color
           indicates frequency of observation, as indicated by the colorbar.
           Panels in the left column are for scheme~(C0), the right
           column corresponds to~(C1). From top to bottom, the three
           rows show the likelihoods of exact hits, misses by exactly one,
           and misses by at least two.}
  \label{fig:time3}
\end{figure}
\begin{figure} \centering
  \setlength{\unitlength}{1 cm}
  \begin{picture}(16.75,16.0)
    \put(0.0,11.0){
      \includegraphics[height=5cm]{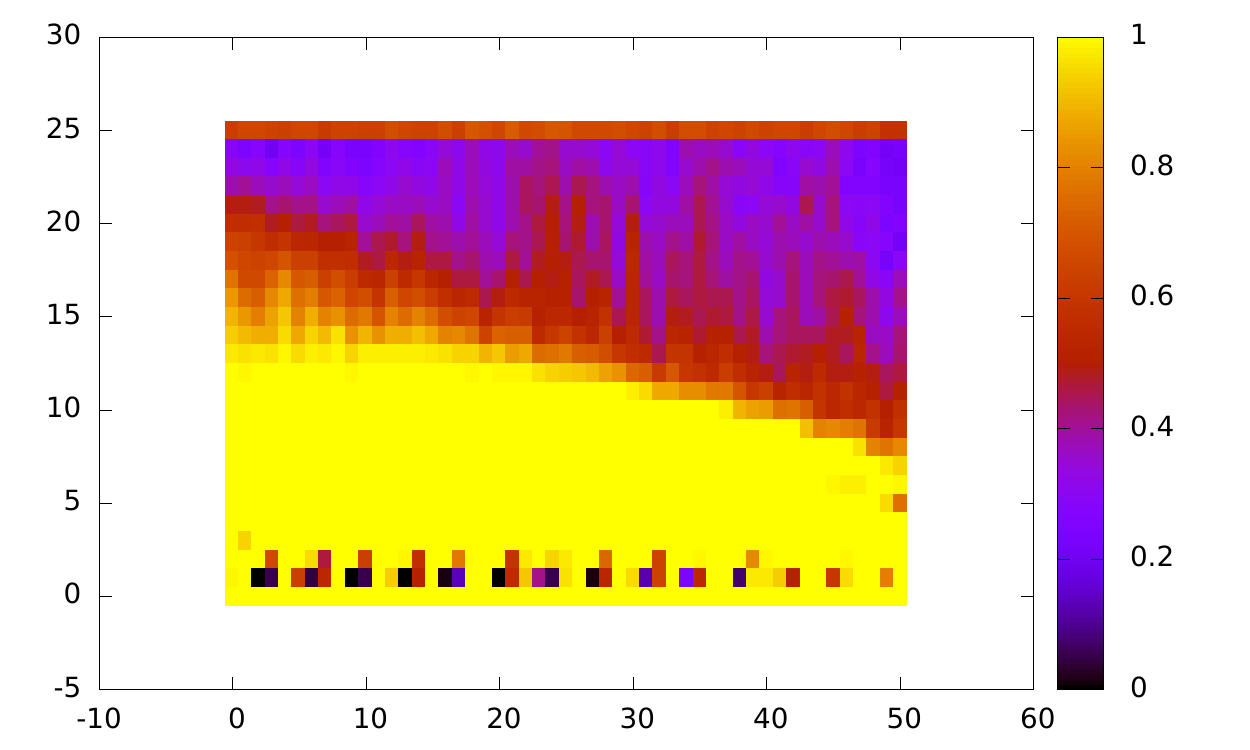}}
    \put(8.75,11.0){
      \includegraphics[height=5cm]{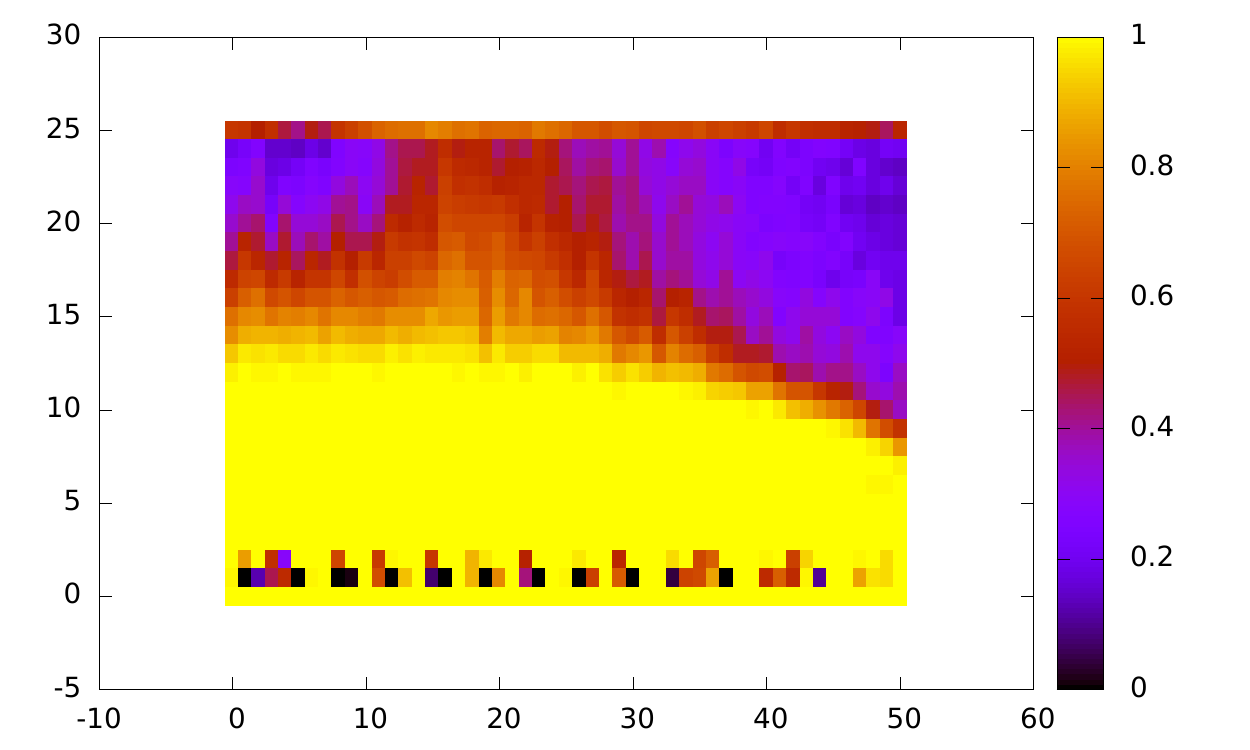}}
    \put(0.0,5.5){
      \includegraphics[height=5cm]{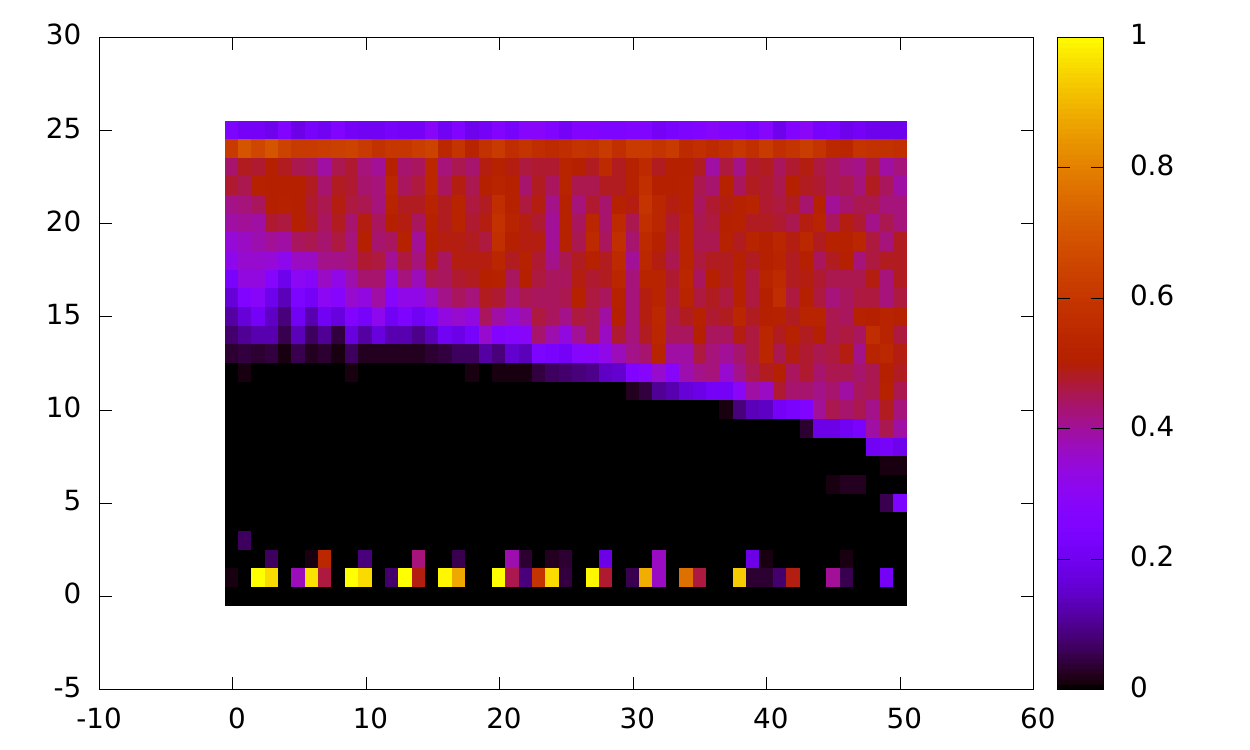}}
    \put(8.75,5.5){
      \includegraphics[height=5cm]{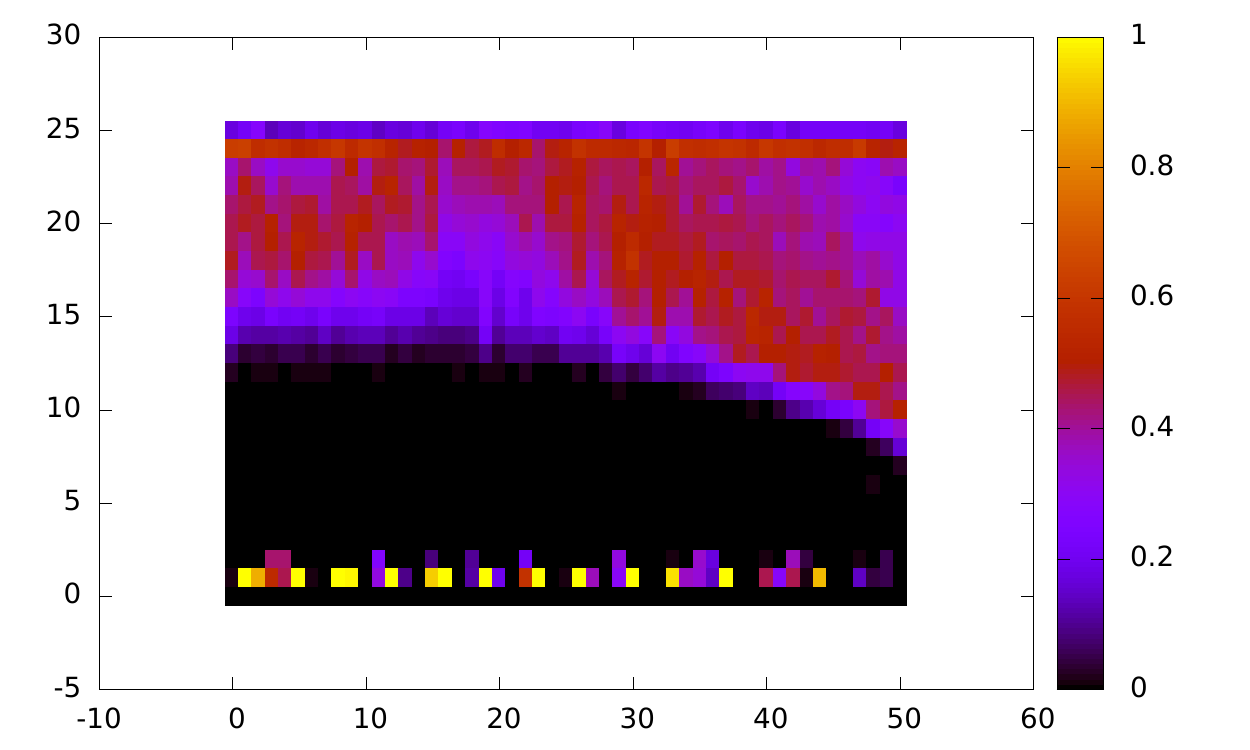}}
    \put(0.0,0.0){
      \includegraphics[height=5cm]{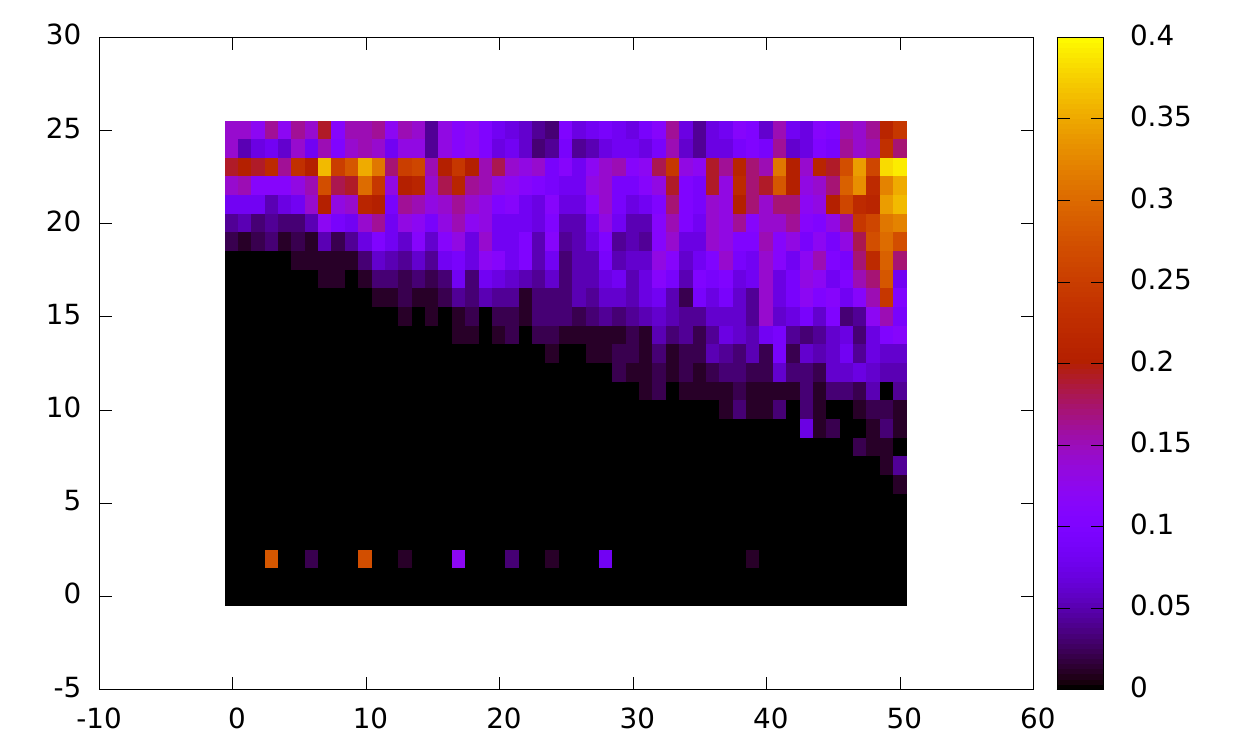}}
    \put(8.75,0.0){
      \includegraphics[height=5cm]{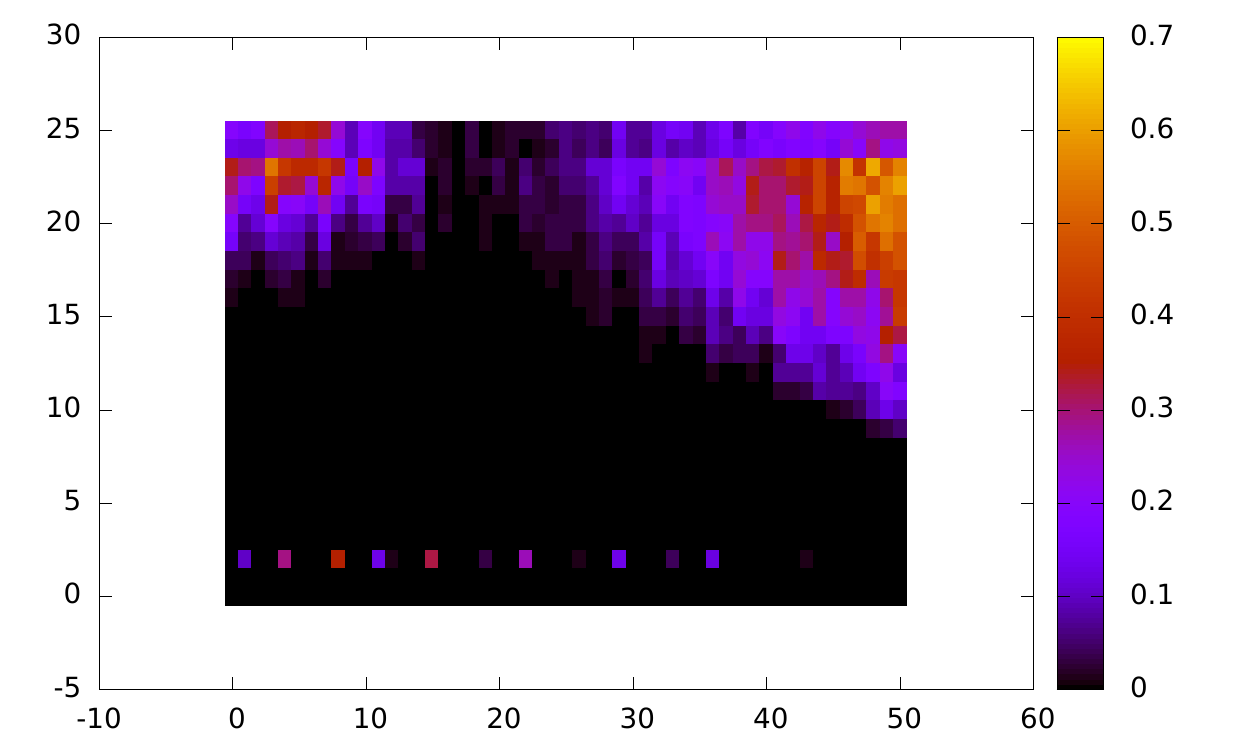}}
  \end{picture}
  \caption{Recovering the snapshot time from persistence information
           using the $L^2$-norm and training set size~$R_T = 20$. Each
           panel contains a heat map, whose horizontal axis corresponds
           to the~$N = 51$ mass values~$\mu$, and the vertical axis
           to the~$M = 25$ snapshot times, in increasing order. Color
           indicates frequency of observation, as indicated by the colorbar.
           Panels in the left column are for scheme~(C0), the right
           column corresponds to~(C1). From top to bottom, the three
           rows show the likelihoods of exact hits, misses by exactly one,
           and misses by at least two.}
  \label{fig:time4}
\end{figure}

For the above results we used different norms, but kept the training
set size fixed at~$R_T = 50$. In order to see how this size affects the
results, we have repeated the simulations, but this time with smaller
training sets of size~$R_T = 20$. The corresponding results are shown
in Figures~\ref{fig:time3} and~\ref{fig:time4}. Notice that despite the
smaller training set size, and the thereby induced larger simulation size
of $(R-R_T) M = 2000$ microstructures which have to be classified, the
results are almost the same. This indicates again that the topological
information which is captured by the random microstructures and their
evolution carries enormous information, enough to accurately detect the
snapshot time, and therefore the stage of material decomposition. The
only apparent exception to this is the case of large mass values~$\mu
\approx 0.5$ and time values well into the coarsening regime. In the
first case, the microstructures are of moving droplet type as shown
in the bottom row of Figure~\ref{fig:chcpatt}, and in the latter case
microstructure evolve slowly since there is less energy to dissipate,
and therefore there are fewer topology changes. But during the initial
spinodal decomposition regime, and well into the coarsening regime, 
topological information can be used to find the snapshot time.
\section{Summary and Conclusions}
\label{sec:future}
Homology is a way of quantifying complicated topological spaces
through a sequence of discrete objects, in fact, through a sequence
of integers in its most reduced form. This information is invariant
under continuous transformations of the underlying space, and in
particular, it does not contain any size or ``metric'' information
of the considered object. The only structures that ``count'', in the true
sense of the word, are $k$-dimensional holes, which cannot be contracted
within the space. In this form, homology has long been used in classical
mathematics to determine when two spaces are different within a certain
category, be it that they are not homeomorphic or not homotopy equivalent.
For data analysis, however, often the opposite question is of interest.
Can homology be used to say that two objects are similar, and if so,
how similar are they?

At its core, the above question is concerned with the amount of information
that is retained by homology. In the present paper, we study this problem 
in a well-defined and important situation, namely phase separation dynamics
in binary metal alloys as described by the stochastic Cahn-Hilliard-Cook 
model~(\ref{chc}). The equation can be used to create realistic evolutions
of phase separating microstructures, which due to their inherent complexity
cannot easily be classified using standard techniques. By using persistence
landscapes, we could obtain a sequence of discrete objects, which encapsulates
the topology evolution of a solution path. Since the Cahn-Hilliard-Cook model
is stochastic in nature, and in addition, since the initial alloy state is never
known precisely, one can only hope to capture typical behavior, and we do this
by averaging these persistence landscapes. It was demonstrated in this paper,
that these averaged landscapes encode enough information to make the following
decision problems solvable:
\begin{itemize}
\item Given only topological information about a specific microstructure
evolution, can we determine the total mass~$\mu$ of the underlying 
experiment?
\item Given only topological information about a specific microstructure
at a given mass value~$\mu$, can we determine the snapshot time during
the decomposition process at which this pattern was observed?
\end{itemize}
These questions can be answered by comparing the provided topological
information to the nearest averaged classifiers, which are created through
training sets.

We would like to stress that we did not set out to develop a method which
determines the total mass~$\mu$ from topological information --- that can 
obviously be done easier by just finding the average gray level of any solution
snapshot. Rather, we wanted to demonstrate that the seemingly stark reduction
in information, which happens during the passage from a complicated microstructure
to its persistent homology, still retains an incredible amount of information
about the underlying dynamical process. Yet, due to its inherently discrete
nature, this topological approach shows considerable promise for analyzing
and relating experimental data to numerical models, and this avenue has to
be explored further in future work.
\section*{Acknowledgements}
The first author was supported by the Advanced Grant of the European Research
Council GUDHI (Geometric Understanding in Higher Dimensions), by DARPA grant
FA9550-12-1-0416, by AFOSR grant FA9550-14-1-0012, and by a Foundation for
Polish Science Start Scholarship. The second author was partially supported
by National Science Foundation grants DMS-1114923 and DMS-1407087. The numerical
simulations for this work were run on ARGO, a research computing
cluster provided by the Office of Research Computing at George
Mason University, VA. (URL: http://orc.gmu.edu)

\end{document}